\documentclass{amsart}

\usepackage[english]{babel}
\usepackage{amsmath,amssymb,amsbsy,amsfonts,amsthm,latexsym,amsopn,amstext,amsxtra,euscript,amscd,paralist}
\usepackage{hyperref}
\usepackage{array}
\usepackage{ifthen}
\usepackage{url}

\usepackage{graphicx}%
\usepackage{multirow}%
\usepackage{amsmath,amssymb,amsfonts}%
\usepackage{amsthm}%
\usepackage{mathrsfs}%
\usepackage[title]{appendix}%
\usepackage{xcolor}%
\usepackage{textcomp}%
\usepackage{manyfoot}%
\usepackage{booktabs}%
\usepackage{algorithm}%
\usepackage{algorithmicx}%
\usepackage{algpseudocode}%
\usepackage{listings}%
\usepackage[english]{babel}
\usepackage{amsmath,amssymb,amsbsy,amsfonts,amsthm,latexsym,amsopn,amstext,amsxtra,euscript,amscd,paralist}
\usepackage{hyperref}
\usepackage{array}
\usepackage{ifthen}
\usepackage{url}

\newcommand{\GG}{\mathbb{G}}

\newcommand{\RR}{\mathbb{R}}
\newcommand{\ZZ}{\mathbb{Z}}

\newcommand{\Pc}{\mathcal{P}}

\newcommand{\Sc}{\mathcal{S}}
\newcommand{\T}{\mathcal{T}}

\newcommand{\set}[1]{\left\{ #1 \right\}}
\newcommand{\setb}[1]{\left( #1 \right)}
\newcommand{\abs}[1]{\left| #1 \right|}

\newtheorem{thm}{Theorem}

\newtheorem{lem}{Lemma}
\newtheorem{pro}{Proposition}

\newtheorem{defin}{Definition}

\newtheorem{remark}{Remark}

\newcommand{\nbf}{\mathbf{n}}
\newcommand{\xbf}{\mathbf{x}}
\newcommand{\bsa}{\boldsymbol{\alpha}}
\newcommand{\KK}{\mathbb{K}}
\newcommand{\G}{\mathcal{G}}

\newcommand{\ri}{m_i}

\theoremstyle{thmstyleone}%
\newtheorem{theorem}{Theorem}
\newtheorem{proposition}[theorem]{Proposition}%

\theoremstyle{thmstyletwo}%

\theoremstyle{thmstylethree}%

\begin{document}

\title[GCD for polynomials ring.]{On quantitative aspects of ideal membership certificate and GCD inequality  for polynomials evaluated at almost units.}

\author[A. Noubissie]{Armand Noubissie}
\address{Armand Noubissie\newline
	\indent Graz University of Technology\newline
	\indent Institute of Analysis and Number Theory\newline
	\indent M\"{u}nzgrabenstrasse 36/II \newline
	\indent A-8010 Graz, Austria.}
\email{armand.noubissie@tugraz.at}



\date{\today}
\begin{abstract}
In this paper, we obtain an effective bounds for  the ideal membership certificate in polynomial rings over almost all rings of $S$-integers, thereby generalizing results of Aschenbrenner and Evertse–Győry. We also establish an effective refinement of a result of Xiao, building on earlier work of Levin and Corvaja-Zannier\cite{CZ2}.  In particular, we provide explicit bounds for both the number and the total degrees of the exceptional hypersurfaces arising in an inequality that bounds the greatest common divisor of multivariable polynomials evaluated at almost S-units under consideration.

As an application of our results, we obtain an effective refinement of a theorem of Corvaja–Zannier \cite[Theorem 2]{CZ5}:  Let $(F_n)_n$ and $(G_n)_n$ be two linear recurrence sequences and $\mathcal{N}$ be the set of positive integers $n$ such that $F_n/G_n$ is an integer. Hadamard quotient theorem, conjectured by Pisot and proven by Van der Poorten states that  if the set  $\mathcal{N}$ coincides with $\mathbb{N}$ then $\left(F_n/G_n\right)_n$ is itself linear recurrence sequence. Corvaja-Zannier  proved that if  $\mathcal{N}$ is infinite set, then  there exists a nonzero polynomial $P(X)$ such that $\left( P(n)F_n/G_n\right)_n$ agrees with a linear recurrence sequence for all $n$ in a suitable arithmetic progression. We give an explicit bound $C$, depending explicitly on the parameters of the recurrence sequences  such that, the conclusion holds if $\mathcal{N}$ contains at least $C$ elements.
\end{abstract}

\maketitle

\section{Introduction and main results.}
We investigate two related questions: the complexity of the ideal membership problem for polynomial rings over the ring of $S$-integers and quantitative aspects of the greatest common divisor of polynomials evaluated at almost units.

\subsection{On the complexity of the ideal membership problem for polynomial rings.}

The study of the complexity of ideal membership in polynomial rings over a field was motivated by a project of Hentzelt and Noether \cite{HN}, which aimed to establish bounds for fundamental operations in commutative algebra.  In this direction, Hermann \cite{H} gave, in 1926, a bound for the complexity of the ideal membership problem for polynomial rings over a field by proving the following result:
{\it Given $f_0, \cdots, f_n \in F[x_1, \cdots, x_N]$ of total degree at most $d$ over a field $F$. If $f_0 \in <f_1, \cdots, f_n>$, then $$f_0=f_1g_1+ \cdots +f_ng_n$$ for certain $g_1, \cdots, g_n \in F[x_1, \cdots, x_N]$ whose degrees are bounded by an explicit constant $\beta(N, d)$ only depending on $N$ and $d$ (not on the polynomials $f_1, \cdots, f_n$ and the field $F$).\
}
Although Hermann's original proof was incomplete, Seidenberg \cite{Se} subsequently supplied a complete proof in 1974. Renschuch \cite{R} later showed that one may take $$\beta(N,d) = (2^d)^{2^N}.$$ Mayr and Meyer \cite{MM} proved that this bound is essentially optimal in the sense that the doubly exponential growth is unavoidable. However, Sharper bounds are available in special cases. When $f_0=1$, Kollár \cite{K}, using Hilbert's Nullstellensatz, obtained the bound $(2d)^{2^N}$ for $d>2$; for $d=2$  Sombra \cite{So} established the bound $2^{N+1}$.  Further improvements are possible under additional geometric hypotheses on the variety defined by $<f_1, \cdots, f_n>$,  for instance, hypotheses controlling its dimension.

A fundamental question raised by Kronecker (see \cite{Ed}) in the context of constructive mathematics is whether there exists an effective procedure for deciding ideal membership in polynomial rings over $\mathbb{Z}$. This problem has been studied extensively; see, for instance, \cite{Ay}, \cite{Ev}, and \cite{Ba}. In the broader context of deciding the truth of universal sentences in classes of commutative rings, Simmons \cite{Si} provided a nondeterministic algorithm (in other words, a procedure that is not primitive recursive, in the sense that it uses an unbounded loop such as "while") to decide whether a given polynomial $f_0 \in \mathbb{Z}[x_1, \cdots, x_N]$ belongs to a given ideal of $\mathbb{Z}[x_1, \cdots, x_N]$ generated by polynomials $f_1, \cdots, f_n$. The question of finding a primitive recursive procedure (i.e., a procedure that can be obtained using only the loop "for")  for this task is commonly referred to as the Kronecker problem. To address this problem, Gallo and Mishra \cite{GM} adapted Buchberger's algorithm for the construction of Gröbner bases and obtained a primitive recursive procedure for deciding the ideal membership problem for polynomials in $N$ variables over $\mathbb{Z}$, thereby answering Kronecker's problem. By analyzing their algorithm, they obtained the following bound on its complexity: $$W_{4N+8}(h(f_0) + \cdots + h(f_n) +N),$$ where $h(\cdot)$ is the logarithm Weil function that we will define in section \ref{sec1} and $W_k$ is often call the $k$th Wainer hierachy function (see \cite{Wa}). For the initial terms, we have $$w_0 (n)= n+1, ~~w_1 (n)= 2n+1, ~~w_2 (n)= 2^n $$ and $w_3$ is comparable to the $n$-times iterated exponential function. Clearly, this function grows extremely rapidly: depending on the indexing convention, the initial levels describe elementary growth rates, while higher levels eventually dominate iterated exponential functions of increasing height. Thus, although the procedure is primitive recursive, its complexity bound is very large and not explicit. 

Notice that a fundamental difference between the ideal membership problem over a field and that over $\mathbb{Z}$ is that, in the former setting, Hermann's degree bound is uniform in the coefficients of the input polynomials. Over $\mathbb{Z}$, by contrast, the  bound for some ideal membership certificates can not depend only on the number of variables and the degrees of the input polynomials. The following example, due to Aschenbrenner \cite{Asc}, illustrates the necessity of incorporating arithmetic information into such bounds.
{\it Let $m > 1 $ and $d \geq 1$ be integers. We have $$1 \in <1-mX, ~m^dX> $$ meaning that $1$ belongs to the ideal of $\mathbb{Z}[X]$ generated by $1-mX$ and $m^dX$ since $$1 = (1-mX)(1+ mX + \cdots + (mX)^{d-1}) + X^{d-1} m^dX.$$ Using the fact that $1-mX$ is unit in $(\mathbb{Z} / m^d\mathbb{Z})[X]$, it follows that, if $$1= (1-mX)g_1 +  m^dX g_2$$ then $\max \left \{ \deg g_1,  \deg g_2 \right \} \geq d-1$ which give us the desired result.}

To obtain explicit bounds for the degrees of ideal membership certificates over $\mathbb{Z}$, Aschenbrenner \cite{Asc} established the following result.

\begin{thm}\label{thmA}
Let $f_0, \cdots, f_n \in \mathbb{Z}[x_1, \cdots, x_N]$ be the polynomials with $$ f_0 \in <f_1, \cdots,f_n> $$ whose their total degrees are at most $d$ and logarithm Weil height of $f_1, \cdots, f_n$ (define in section \ref{sec1}) are at most $h$, then $$f_0 = g_1f_1 + \cdots + g_nf_n$$ for certain polynomials $g_1, \cdots, g_n \in \mathbb{Z}[x_1, \cdots, x_N]$ of degree at most $$\gamma (N,~d,~ h) = (2d)^{2^{O(N \log N +1)}}(h+1).$$   
\end{thm}

The constant hidden in the $O$-notation in $\gamma$ is explicit. Evertse and Gy\H{o}ry \cite[Theorem 6.1.2]{evertse-gyory-2022} combined Aschenbrenner's result (Theorem \ref{thmA}) with the result of Borosh et al. \cite{Bo} to obtain a deterministic algorithm for the ideal membership problem for polynomial rings over $\mathbb{Z}$, together with an explicit bound for the complexity.\\

In what follows, we extend Aschenbrenner's result to the following setting. Let $K$ be a number field of degree $d$ and $S$ be a finite set of places of $K$ containing all archimedean places and the nonarchimedean places corresponding to prime ideals dividing a set of integral ideals representing the ideal classes of the class group  of $K$. 
Let $\mathcal{P}_1, \cdots, \mathcal{P}_t$ denote the prime ideals corresponding to the finite places in $S$ and set $$Q_S = N(\mathcal{P}_1\cdots \mathcal{P}_t),$$ where $N(\cdot)$ denotes the norm function on $K$. Let $S_{\infty}$ be the set of all archimedean places of $K$. Finally, set $$c_4:=87d!\sqrt{d-1}d\log d.$$
From \cite[Proposition 5.3.6, Chap V]{bombieri-gubler-2006}, it is clear that $\mathcal{O}_S$ is a principal ideal domain. Let $m,n$ be two positive integers and $a=\max \{m, n\}$. We set $D= \mathcal{O}_S[X_1, \cdots, X_N]$. Let $A$ be an $m \times n$-matrix with entries $a_{ij} \in D$, and let $b \in D^m$ be a column vector. Let $U$ be the maximum of the total degrees of the entries $a_{ij}$ and the coordinates of $b$, and  put $h:= h^{hom}(1, a_{ij}, b_i~|~ ~ (i,j) \in \{1, \cdots, m\} \times \{1, \cdots, n\})$ be the homogeneous logarithm height of the vectors formed by all the coefficients of  $a_{ij}, b_i$ and $1$ (see section \ref{sec1}, for definition the homogeneous logarithm height).  We denote by $R_K$ the regulator of $K$ and $h_K$ the class number of $K$.

\begin{thm}\label{thmB}
If the system $AX = b$ has a solution in $D^n$, then it admits a solution whose coordinates have total degrees at most  $$C := \gamma_1 \cdot  3 \left(a \binom{N+ U+ \gamma}{N} \right)^3 \log \left(a\binom{N+ U+ \gamma}{N} \right) \left(h +T +1 \right),$$ where $$\gamma_1 = d 3^{2^{3N}}(2 (a+1) \cdot U)^{2^{3N+1}(N+1)^N} $$ and 
$$\gamma = (2 (a+1) \cdot U)^{2((N+1)^N-1)}$$ and $$T=c_4R_K + h_K \log Q_{S}.$$ If $S= S_{\infty}$, then we replace $T$ by $c_4R_K +1$.
\end{thm}

\begin{remark}
Aschenbrenner \cite[Theorem 6.5]{Asc} established a non-effective analogue of Theorem \ref{thmB} when  $\mathcal{O}_S$ is replaced by the ring of algebraic integers $\mathcal{O}_K$ of the number field $K$. The non-effectivity arises from his reliance on \cite[Lemma 1.4]{Asc}, whose proof uses a qualitative result of Lang \cite[Chapter V, Theorem 1]{L1}. 
\end{remark}

A related problem is to obtain explicit height bounds for solutions of systems of homogeneous linear equations over rings of integers. Siegel initiated this line of research by proving what is now known as Siegel's lemma in the case $\mathcal{O}_K = \mathbb{Z}$. Bombieri and Vaaler \cite{BV} subsequently extended this result to number fields, obtaining a basis constitute of elements in $\mathcal{O}_K^N$ of the subspace $V$ associated with the system, whose such elements are bounded in height in terms of the discriminant of $K$. Results of Roy and Thunder \cite{RT} and Masser \cite{Mas} show that the dependence on the discriminant in the Bombieri--Vaaler bounds is essentially optimal.\\

We extend this line of investigation in the following setting.  Let $m,n$ be two positive integers and $a=\max\{m, n\}$.  Let $A$ be an $m \times n$- matrix with entries  $a_{ij} \in \mathcal{O}_S$. Define  the $\mathcal{O}_S$-module of solutions of the associated homogeneous system by $${\rm Sol_{\mathcal{O}_S}}(A) = \{X \in \mathcal{O}_S^n~~|~~ AX=0\}$$ and   put $h= h^{hom}(1, a_{ij}~|~ ~ (i,j) \in \{1, \cdots, m\} \times \{1, \cdots, n\})$, where $h^{hom}(\cdot)$ is defined in section \ref{sec1}. Since  $\mathcal{O}_S$ is a Noetherian ring, ${\rm Sol_{\mathcal{O}_S}}(A)$ is a finite generated $\mathcal{O}_S$-module. 

\begin{thm}\label{thmC}
There exists a generator $\{w^{(1)}, \cdots, w^{(k)}\}$   of  $~~{\rm Sol_{\mathcal{O}_S}}(A)$ such that  the logarithm Weil height of any coordinate of the vectors $w^{(i)}$ is bounded above by $$ a \left( \frac{5}{2}a^2(\log a + h)+ c_4R_K + h_K \log Q_S \right),$$
where $R_K, h_K$ are respectively the regulator and the class number of $K$ and 
$$c_4:=87d!\sqrt{d-1}d\log d, \quad Q_S = N(\mathcal{P}_1\cdots \mathcal{P}_t).$$
 If $S= S_{\infty}$, then we replace $ h_K \log Q_S$ by $1$.
\end{thm}

\begin{remark}
 Theorem \ref{thmC} plays a crucial role for proving Theorem \ref{thmB}. In \cite[Corollary 5.4]{Asc}, Aschenbrenner established a non-effective version of Theorem \ref{thmC} when $\mathcal{O}_S$ is replaced by the ring of algebraic integers $\mathcal{O}_K$. The non-effectivity arises from his reliance on \cite[Lemma 1.4]{Asc}, whose proof uses a qualitative result of Lang \cite[Chapter V, Theorem 1]{L1}.
\end{remark}

Let $K$ be a number field of degree $d$ and $S$ be a finite set of places of $K$ containing all archimedean places and the nonarchimedean places corresponding to prime ideals dividing a set of integral ideals representing the ideal classes of the class group  of $K$. We denote the cardinality of $S$ by $s$.
Let $\mathcal{P}_1, \cdots, \mathcal{P}_t$ denote the prime ideals corresponding to the finite places in $S$. 
Given $f_1, \cdots, f_m \in \mathcal{O}_S[x_1, \cdots,x_n]$, we write  $ I_K = < f_1, \cdots, f_m >$  for the ideal they generate in $K[x_1, \cdots,x_n]$.  A natural problem is to determine an explicit set of generators for the contraction  $$<f_1, \cdots, f_m>  \bigcap  \mathcal{O}_S[x_1, \cdots,x_n]$$ that is, the ideal obtained by contracting $I_K$ to  $\mathcal{O}_S[x_1, \cdots,x_n]$. This problem plays an important role in the study of Diophantine equations over finitely generated domains; see \cite[Chap IX]{evertse-gyory-2022}.

Evertse and Győry \cite[Theorem 6.1.7]{evertse-gyory-2022} solved this problem in the case  $K= \mathbb{Q}$ and $\mathcal{O}_S = \mathbb{Z}$.  We extend their result, providing an explicit description of the contracted ideal in this more general setting. We denote by $R_K$ the regulator of $K$ and $h_K$ the class number of $K$. For the definition of the homogeneous logarithm height (see section \ref{sec1}).

\begin{thm}\label{thmD}
Let $f_1, \cdots, f_m \in \mathcal{O}_S[x_1, \cdots, x_N]$ be the polynomials of total degrees  at most $D$ and  $h = h^{hom}(f_1, \cdots , f_m)$ be the homogeneous logarithm height of the vector formed by their coefficients. Then $<f_1, \cdots, f_m>  \bigcap  \mathcal{O}_S[x_1, \cdots,x_N]$ is an ideal of 
$\mathcal{O}_S[x_1, \cdots,x_N]$ generated by polynomials of total degrees at most  $$ c :=  (2 m\cdot D)^{2((N+1)^N-1)}$$ and  the homogeneous logarithm height  at most $$C_2 := (h+1) \left (\binom{c+D+N}{N} |d \Delta_K|\right)+ (s-1)c_5 h_K R_K {\displaystyle \prod_{i=1}^{t}} \log N(\mathcal{P}_i),$$  with
 $$c_5 :=  \dfrac{((s-1)!)^2}{2^{s-2}d^{s-1}}(53 d^2 \log 6d)^{s-2}$$
and $\Delta_K$  is the discriminant of $K$.  If $S= S_{\infty}$, then we  replace $h_K {\displaystyle \prod_{i=1}^{t}} \log N(\mathcal{P}_i)$ by $1$.
\end{thm}  

\subsection{ On quantitative aspect of the greatest common divisor for polynomials.}

Let $a, b > 1$ be integers. If $a,~b$ are multiplicatively dependent, then there exist coprime positive integers $r, s$ and an integer $c>1$ such that  $a=c^s$ and $b=c^r$. Consequently, $c^n-1$ divide ${\rm gcd}(a^n-1, b^n-1)$ for every integer $n$. This makes multiplicative independence the natural hypothesis under which to study the size ${\rm gcd}(a^n-1, b^n-1)$. Bugeaud, Corvaja and Zannier \cite{BCZ} established the following fundamental result. 

\begin{thm}\label{thmE}
Let $a,b >1$ be multiplicatively independent integers, $\epsilon >0$ be a positive real number. Then, the inequality
\begin{equation} \label{eq:a'}
{\rm log } \gcd (a^n-1, b^n-1) > \epsilon n.
\end{equation}
has finitely many solutions in $n$.
\end{thm}

On the other hand, Addleman et al. \cite[Proposition 1]{A} implies that there exists a constant $c>0$ such that for all pair $a,b >1$ there is infinitely positive integers $n$ with $${\rm log } \gcd (a^n-1, b^n-1) > \exp(c\log n / \log \log n).$$ This implies that the relation \eqref{eq:a'} is closed to be the best possible. To formulate extensions of Theorem \ref{thmE} over number fields, we first introduce a natural generalization of the $\gcd$.  For nonzero integers $a$ and $b$, we have 
\begin{equation*}\label{eq:E}
\begin{split}
\log \gcd (a,b) &=  {\displaystyle \sum_{p~ prime}} \min\{ {\rm ord}_p(a), {\rm ord}_p(b) \} \log p\\
& =  -  {\displaystyle \sum_{v \in M_{\mathbb{Q}}^0}} \log \max \{|a|_v, |b|_v\} \\
& = -  {\displaystyle \sum_{v \in M_{\mathbb{Q}}}} \log^{-} \max \{|a|_v, |b|_v\}, 
\end{split}
\end{equation*}
where $M_{\mathbb{Q}}^0$ denotes the set of nonarchimedean places of $\mathbb{Q}$ and $\log^{-}z = \min\{0, \log z\}$. It is clear that the quantity $$\log \gcd (\alpha, \beta)  = - {\displaystyle \sum_{v \in M_{K}}} \log^{-} \max \{|\alpha|_v, |\beta|_v\}$$ is well defined for all $\alpha, \beta \in K.$ Note that  $$ \log \gcd(\alpha,\beta)= h^{hom}(1, \alpha,\beta) - h^{hom}( \alpha,\beta),$$ where the homogeneous logarithm Weil height $h^{hom}(\cdot)$ is defined in section \ref{sec1}. 

Corvaja and Zannier \cite{CZ1}, independently of Hernández and Luca \cite{HL}, extended Theorem \ref{thmE} to elements of finitely generated multiplicative groups of  $\mathbb{Q}.$ One formulation of their result is the following theorem.

\begin{thm}\label{thmF}
Let $p_1, \cdots, p_t \in \mathbb{Z}$ be prime numbers and $S=\{\infty, p_1, \cdots, p_t\}$. Then for every $\epsilon >0,$ 
$${\rm log } \gcd (u-1, v-1) < \epsilon \max \{\log |u|, \log|v|\}.$$
for all but finitely many multiplicatively independent $S$-unit integers $u, v \in \mathbb{Z}_{S}^*$.
\end{thm}

Note that Theorem \ref{thmE} and Theorem \ref{thmF} have also been studied in the function field setting. Indeed, Ailon and Rudnick \cite{AR} showed that, if we replace $a$ and $b$ in \eqref{eq:a'} by multiplicatively independent polynomials  $f, g \in \mathbb{C}[x]$, then there exists a polynomial $h$ such that  $  \gcd (f^n-1, g^n-1)| h $ 
for all $n >1.$  In particular, their result yields a uniform bound, independent of $n$, for the degree of the greatest common divisor of $f^n-1$ and $g^n-1$. Their proof relies on a theorem concerning the finiteness of torsion points on algebraic curves in the two-dimensional algebraic torus, a result arising from a conjecture of Lang \cite[Theorem 6.1, Chapter 8]{L1} and proved in \cite{La}, following work attributed to Ihara, Serre, and Tate. Subsequently, Beukers and Smyth \cite{BS} obtained an explicit bound for the number of torsion points that is independent of the genus of the curve, while Corvaja and Zannier \cite{CZa} bounded the maximal order of such points. Maurin \cite{Ma} later extended the relevant finiteness results. In 2020, Ostafe and Shparlinski \cite{OS} used Maurin's theorem to extend the result of Ailon and Rudnick to $\gcd(h_1(f^n, g^n), h_2(f^n, g^n) )$, where polynomials $h_1, h_2 \in \mathbb{C}[x,y]$. The function field analogue of Theorem \ref{thmF} was established in \cite{CZ4}, where a sharper bound was obtained.

Returning to the number field setting, the multiplicative independence of $u,v$ in Theorem \ref{thmF} is equivalent to the condition that the point $(u,v) \in \mathbb{G}_m^2$ does not lie in  a proper algebraic subgroup of $\mathbb{G}_m^2$.  Corvaja and Zannier \cite{CZ2} obtained a further generalization of Theorem \ref{thmF} by replacing  $u-1$ and $v-1$ with more general polynomial expressions evaluated at  $u$ and $v$. Their result provides the following theorem. 

\begin{thm}\label{thmG}
Let $\Gamma \subseteq \mathbb{G}_m^2(\overline{\mathbb{Q}}) $ be a finitely generated group. Let $f(x,y),g(x,y) \in \overline{\mathbb{Q}}[x,y]$ be non-constant coprime polynomials such that not both of them vanish at $(0,0)$. For all $\epsilon >0,$ there is a finite union $Z$ of translate of proper algebraic subgroup of $\mathbb{G}_m^2(\overline{\mathbb{Q}})$ such that 
$${\rm log } \gcd (f(u,v), g(u,v)) < \epsilon \max \{h(u), h(v)\}$$
for all $(u,v) \in \Gamma \setminus Z$, where $h(\cdot)$ denotes the logarithm Weil height function defined in section \ref{sec1}.
\end{thm}

Levin \cite{L} extended Theorem \ref{thmG} to higher dimensions by establishing $\gcd$ estimates for multivariate polynomials evaluated at $S$-unit points. More recently, Xiao \cite{X} generalized Levin's result to multivariate polynomials evaluated at almost $S$-unit points. Let $K$ be a number field of degree $d$

\begin{thm}\label{thmI}
Let $K$ be a number field of degree $d$ and $S$ be a finite set of places of $K$ containing all archimedean ones. Let $f,g \in K[x_1, \cdots, x_n]$ be non-constant coprime polynomials such that not both of them vanish at $(0, \cdots,0)$. For all $\epsilon >0,$ there exists $\delta >0$ and a finite union $Z$ of  hypersurfaces such that 
$${\rm log } \gcd (f(u_1, \cdots, u_n), g(u_1, \cdots, u_n)) < \epsilon \max \{h(u_1), \cdots, h(u_n)\}$$
for all $(u_1, \cdots, u_n) \in \mathbb{G}_m^n(K)_{S, \delta} \setminus Z$.  For the definition of  $\mathbb{G}_m^n(K)_{S, \delta}$ (see Definition \ref{def}).  
\end{thm}

Theorem \ref{thmI} is qualitative in that it provides no explicit bounds for either the number or the total degree of the exceptional hypersurfaces. We establish a quantitative refinement of this result by obtaining explicit bounds for both quantities.

\begin{thm}\label{thmJ}
Let $K$ be a number field of degree $d$ and $S$ be a finite set of places of $K$ containing all archimedean ones. We denote the cardinality of $S$ by $s$. Let $f, g \in K[x_1, \cdots, x_N]$ be non-constant coprime polynomials of total degree at most $D$ such that not both of them vanish at $(0, \cdots,0)$  and $ 0< \epsilon < 1 $ be a real number. For $$\eta < \dfrac{\epsilon^2}{72D^2N^4} $$ we consider the inequality 
\begin{equation}\label{eq:b'}
\log \gcd\left ( f(u_1, \cdots, u_N),~ g(u_1, \cdots, u_N) \right ) > \epsilon {\displaystyle \sum_{i=1}^{N}} h(u_i)
\end{equation}
with $(u_1, \cdots, u_N) \in \mathbb{G}_{m}^{N}(K)_{S, \eta}$. There are hypersurfaces $H_1, \cdots , H_{t}$ of $K^N$ of total degree at most $m: =  12D^2N^3\epsilon^{-1} $ with $$t \leq 2{\rm exp}\left( 75sNd(2m)^{2N}(1+\log \epsilon^{-1})\right)\left( \log (6m)+ 2h^{hom}(f, g, 1) \right)^{3Nd}$$ such that any solution $(u_1, \cdots, u_n) \in \mathbb{G}_{m}^{n}(K)_{S, \eta}$ of \eqref {eq:b'} lies in $$H_1 \bigcup  H_2 \bigcup \cdots \bigcup H_{t},$$ where  $h^{hom}(f,g,1)$ denotes the homogeneous logarithm  height of the vector formed by the coefficients of the polynomials $f, g,$ and $1$ as defined in section \ref{sec1}. For the definition of  $\mathbb{G}_m^n(K)_{S, \eta}$ (see Definition \ref{def}).  
\end{thm}

\begin{remark}
Combining Theorem \ref{thmJ}, applied with $f(X, Y)= X - 1$ and $g(X, Y)= Y - 1$  with Lemma \ref{lemJ}, we obtain an explicit upper bound for the number of positive integers $n$ satisfying \eqref{eq:a'}. More precisely, $$\# \{n~:~{\rm log } \gcd (a^n-1, b^n-1) > \epsilon n\} < \zeta  10^{6} {\rm exp}\left(4.2\cdot s^6\cdot 10^{11} \epsilon^{-5}\right)  \left(\max \{|a|, |b|\} \right)^{2},$$ where $$ s= 2 \log \max\{|a|, |b|\}, ~~\zeta = \exp((7k)^{8k}), \quad \mbox{with} \quad k=\binom{96\epsilon^{-1}+2}{2}.$$ Thus, our results make the finiteness assertion of Theorem \ref{thmE} effective by providing an explicit upper bound for the number of exceptional indices in terms of $a, b$, and $\epsilon$. 
\end{remark}

\subsection{Applications.}

Our results may be interpreted from the perspective of the quantitative form of Vojta’s conjecture, as emphasized by Levin \cite[Section 1.2]{L} and Xiao \cite{X}.

\subsubsection{Greatest common divisor from terms of multi-recurrence sequences.}

Linear recurrence sequences (LRS) are fundamental objects in number theory and dynamical systems. A  sequence of  numbers ${u_n}$ is a linear recurrence sequence over $\mathbb{Q}$ of order $l$ if it is defined by the recurrence relation
$$
u_{n+l}= a_1u_{n+l-1}+ a_2u_{n+l-2}+ \cdots + a_lu_{n}, \quad (n \in \mathbb{N}),
$$
where the coefficients $a_1, \cdots, a_l$ are complex numbers. Such sequences admit a polynomial-exponential representation, often called a Binet formula, of the form 
\begin{equation}\label{eq:c'}
u_n = \sum_{i=1}^{k} p_i(n) \alpha_i^n,
\end{equation}
where $\alpha_i$ are the roots of the characteristic polynomial of ${u_n}$, and the coefficients of the polynomials $p_i \in K$ with $K := \mathbb{Q}(\alpha_1, \cdots, \alpha_k)$.  The sequence is called non-degenerate if no ratios of two roots of the characteristic polynomial is a root of unity. It is worth pointing out that any LRS can be effectively decomposed as an interleaving of finitely many non-degenerate LRS. Thus, the non-degenerate case is central to the study of zeros of linear recurrence sequences. The Skolem–Mahler–Lech theorem asserts that the zero set of a linear recurrence sequence over a field of characteristic zero is a finite union of arithmetic progressions and a finite set. In particular, a non-degenerate linear recurrence sequence that is not identically zero has only finitely many zeros. However, the theorem does not provide an effective bound for the number of zeros or for their indices.

In \cite{BM1}, Brestel and Mignotte obtained explicit arithmetic progressions describing the zero set, leading to an algorithm for deciding whether a given linear recurrence sequence has infinitely many zeros. Nevertheless, the general problem of deciding whether a linear recurrence sequence has even a  zero or not  remains open and is known as the Skolem problem. This problem arises in several areas of theoretical computer science and applied mathematics, including the analysis of L-systems and population dynamics, software verification and termination of linear programs \cite{OW15}, stochastic model checking and reachability in Markov chains \cite{AAGT}, and threshold problems for quantum automata \cite{LOW}.

Using Baker's theory of linear forms in logarithms of algebraic numbers, Mignotte et al. \cite{MST} and Vereshchagin \cite{V} independently established the decidability of the Skolem problem for linear recurrence sequences of order at most four. Extending these methods to higher orders has proved difficult, owing to the limitations of current estimates for linear forms in logarithms in higher-dimensional settings.

\par A recent line of research on the Skolem problem has focused on the notion of a universal Skolem set. Rather than imposing restrictions on the recurrence sequence, such as bounds on its order or conditions on the moduli of its characteristic roots, this approach restricts the set of indices at which zeros are sought. A domain $\Gamma \subseteq \mathbb{N}$, for which there exists a recursive procedure to decide whether there exists an $n \in \Gamma $  with $u_n = 0$ for a given integer LRS $\{u_n\}$, is called a universal Skolem set.  In particular, the Skolem problem is decidable if and only if $\mathbb{N}$ is a Universal Skolem Set.

\par From this perspective, a natural objective is to construct universal Skolem sets with lower density as close to 1 as possible. Luca et al.\cite{LOW2} constructed the first infinite universal Skolem set, although its density is zero. Subsequently, Luca et al. \cite{LOW1} constructed a universal Skolem set of positive lower density for simple linear recurrence sequences, for which the algorithm applies under the additional simplicity hypothesis. Their construction relies on a result of Schmidt giving an explicit upper bound, depending only on the order of the recurrence, for the number of solutions of $u_n=0$. This result is based on the parametric Subspace Theorem. In \cite{LMNOW}, we combined this approach with classical methods from analytic number theory to construct a universal Skolem set of lower density at least 0.29 unconditionally, as well as a set of density 1 under the Bateman–Horn conjecture. More recently, Luca, Ouaknine, and Worrell established a conditional decidability result for the Skolem problem under Cramér's conjecture.

\par A related question concerns the rate of growth of linear recurrence sequences. Let $\alpha$ be a dominant root with absolute value greather that $1$ of characteristic polynomial of $\{u_n\}$. It's not difficult to see  there is an effectively computable constant $a$ such that, for all  $n \geq 1, ~~  \vert u_n \vert <  a n^l \vert \alpha  \vert^n$ for a suitable non-negative integer $l$.  In 1977, Loxton and van der Poorten \cite{LV} formulated a conjecture asserting, in essence, that a non-degenerate linear recurrence sequence cannot remain exponentially smaller than its natural maximal growth rate at arbitrarily large indices. More precisely,   for any $\epsilon >0,$ there is an effectively computable constant $C(\epsilon),$ such that , if $ \vert u_n \vert < (\max_i\{ \vert \alpha_i  \vert \})^{n(1-\epsilon)} $ then  $n<C(\epsilon)$. Fuchs and Heintze \cite{FH} established a noneffective version of this conjecture in 2021, building on results of Schmidt and Evertse. More recently, I obtained an explicit upper bound for the number of solutions of the corresponding inequality, using methods related to those employed by Evertse \cite{E} in his quantitative version of the Subspace Theorem; see \cite{N}. 

\par Questions concerning the growth of linear recurrence sequences are also closely related to two fundamental problems in theoretical computer science: the Positivity Problem, which asks whether every term of a given linear recurrence sequence is nonnegative, and the Ultimate Positivity Problem, which asks whether all sufficiently large terms are nonnegative. Ouaknine and Worrell \cite{OW3} proved that the Ultimate Positivity Problem is decidable for simple linear recurrence sequences. Their approach relies substantially on estimates for the growth of linear recurrence sequences.

\par The notion of a linear recurrence sequence admits a natural multivariate generalization. Namely, let $r,t$ be positive integers, and consider a polynomial-exponential function 
\begin{equation} \label{eq:A1}
G(n_1, \ldots, n_t) = \sum_{i=1}^{r} P_i(n_1, \ldots, n_t)\alpha_{i1}^{n_1}\cdots \alpha_{it}^{n_t}, 
 \end{equation}
 where the $P_i$ are polynomials in $t$ variables with coefficients $\textbf{a}_i =(a_{i1}, \ldots, a_{i\ri})$,  for $i=1,\ldots,r$, and $n_1, \ldots, n_t$ are non-negative integers. The polynomial-exponential function $G: \mathbb{Z}^t \rightarrow K$ with $\nbf=(n_1,\ldots,n_t)\mapsto G(\nbf)$, abbreviated by $(G(\nbf))$, is called a multi-recurrence sequence (or multi-recurrence) and equation \eqref{eq:A1} is called its power sum representation. We assume that $\alpha_{ij}$'s are algebraic integers over $K$ and the coefficients of the $P_i$ belong to the number field $K(\alpha_{11},\alpha_{12},\ldots,\alpha_{rt})$. We put $\bsa_i^{\nbf} := \alpha_{i1}^{n_1} \cdots \alpha_{it}^{n_t}$ for $i=1,\ldots,r$. Assume that $\alpha_{11}$ has maximal modulus with $\vert \alpha_{11} \vert >1$ among the $\alpha_{ij}$. We also put $\vert \nbf \vert = \max_i\{\vert n_i \vert \}$.  Schlickewei and van der Poorten \cite{VS} claimed the following result, which was subsequently proved by Fuchs and Heintze \cite{FH1}: 
{\it 
Fix $\epsilon >0$. Let $\mathcal{A}$ be the set of vectors $\nbf$ such that for any subset $I \subseteq \{1, \ldots, r\}$ with $1 \in I,$ we have $$\sum_{i \in I}^{} P_i(\nbf)\bsa_i^{\nbf} \neq 0.$$ Then, the inequality 
\begin{equation} \label{eq:2}
\vert G(\nbf)\vert < \vert P_1(\nbf) \bsa_1^{\nbf} \vert e^{-\epsilon (\vert n_1\vert+\cdots+\vert n_t\vert)}
 \end{equation} 
 has finitely many solutions in $\mathcal{A}$.\\
}

More recently, we obtained an explicit upper bound for the number of solutions of  \eqref{eq:2} in $\mathcal{A}$; see \cite[Theorem 1]{FN}. Let $S$ the set containing all the prime ideals above the $\alpha_{ij}$ and all the Archimedean places over $K$, and let $s$ the cardinality of $S$. Let $\Gamma = \{1, \ldots,r\}$ and $\G$ be the subset of $\mathbb{Z}^t$ consisting of \textbf{z} with $\bsa_l^{\textbf{z}} = \bsa_k^{\textbf{z}}$ for some $l,k \in \Gamma$ with $l \neq k$. We denote by $$A' = \sum_{i \in \Gamma}^{} \binom{t+\ri}{t} \mbox{ and } A= \max \{t, A'\}.$$
Notice that $A = \max\{t,r\}$ if all the polynomials $P_1, \ldots, P_r$ are constants, and $A=A'$ otherwise. We write $$\mathcal{S}_i = \{ \textbf{z} \in \mathbb{Z}^t: P_i(\textbf{z}) = 0\}$$ for $i=1,\ldots,r$. We set $$h := {\rm max}_{i,j} \left \{ h(a_{i,j}), ~ 1  \right \}$$ where $a_{i,j}$ range over the coefficients of the polynomials $P_1, \cdots, P_r.$ By restricting the condition  on the maximal modulus among the $\alpha_{ij}$ we obtain the following result. 

\begin{thm}\label{thmK}
Let $D$  be a positive real number. Assume that $\G = \{0\}$ and  $$H_i := \{(u_1, \cdots, u_t) \in \mathbb{Z}^t ~: ~P_i(u_1, \cdots, u_t) = 0\} $$ has at most $D$ elements and the $P_i$ have total degrees at most $\sigma \geq 1$.  Let $|\cdot |_v$ be an absolute value on K such that, for some $i,$ we have  $|\alpha_{i,j}|_v \geq 1$ for all $j\in \{1, \cdots, t\}$. For $0< \epsilon <1,$ there exists at most $$ \rho(\epsilon):= 4 \left( 2^{60r^2}\delta^{-7r}\right)^{s+1}\frac{r!}{2} 2^{35A^3}d^{6A^2} +D^r + L^t $$ with $$\eta =  \dfrac{\epsilon(r-1)c_6h^r}{(2r +1)(r + 6r^2t c_6h^{r+1})} \quad \mbox{and}\quad \delta =  \dfrac{\epsilon(r-1)c_6h^r}{r + 6r^2t c_6h^{r+1}}$$ and $$ L:= \max\left\{ (2t\sigma)^t\exp\left((2t\sigma)^th \right), (4 \eta^{-1}c_6h^r(\sigma+1)(r-1))^{3/2}\right\}$$ and 
$$c_6 =  \dfrac{58t!3^td^{t+1} \log d}{t^t}$$
values of $\nbf$ such that
\begin{equation}\label{eq:l_1}
\log |G(\nbf)|_v < -\epsilon |\nbf|,
\end{equation}
where $|\nbf| = {\rm max}_i\{n_i\}.$
\end{thm}

\begin{remark}
 In certain cases, the parameter $D$ can be chosen explicitly. If $t=1$ and  all the $P_i$ are constant polynomials, we may take $D=1$. In the case where $P_i$ is a polynomial in one indeterminate, we take $D= {\rm max}_{i} \deg P_i$. Finally, if $t=2$ and each $P_i$ defines an algebraic curve with at least three points at infinity, an explicit value of $D$ can be obtained by combining a quantitative version of the Northcott property (Lemma \ref{lemE}) with an effective version of Siegel's theorem \cite[Theorem 1]{CZ3}. 
\end{remark}

\begin{remark}
A stronger version of this result, under the assumptions that every  $\alpha_{i,j}$ is algebraic integer for all $i,j$ and ${\rm max}_{i,j} |\alpha_{i,j}|_v >1$ is established in \cite[Theorem 1]{FN}. In the case $t=1$, our result also yields effective versions of \cite[Lemma 5.2]{L}, \cite[Lemma 4.1]{X}, and \cite[Lemma 4.1]{GW}. A stronger result was obtained in \cite[Theorem 1]{N}, under the additional hypothesis that the sequence has a dominant root of absolute value greater than $1$ with respect to the place $v$.
\end{remark}

Let $K$ be a number field of degree $d$ and  $(F(\nbf)), ~(G(\nbf))$ be two multi-recurrences defined over $K$ by 
$$
F(n_1, \ldots, n_t) = \sum_{i=1}^{r} Q_i(n_1, \ldots, n_t)c_{i1}^{n_1}\cdots c_{it}^{n_t}, 
$$
and
$$
G(n_1, \ldots, n_t) = \sum_{i=1}^{r} P_i(n_1, \ldots, n_t)\alpha_{i1}^{n_1}\cdots \alpha_{it}^{n_t}, 
$$
 where $r$ and $ t$ are positive integers,  $P_i, Q_i \in K[x_1, \cdots,x_t]$ for $i=1,\ldots,r$, and $n_1, \ldots, n_t$ are non-negative integers and $\alpha_{i,j}, c_{i,j} \in K$ for all $(i,j) \in \{1, \cdots, r\} \times \{1, \cdots, t\}$. $\alpha_{i,j}$ and $ c_{i,j}$ are called "multi-roots" of multi-recurrences  $(F(\nbf)), ~(G(\nbf))$ respectively.
 
 Let $S$ be a set  containing all the archimedean places of $K$ and the nonarchimedean places corresponding to the prime ideals lying above  $\alpha_{ij}$ and $c_{i,j} $. Let $s$ be the number of elements of $S$. Put $\Gamma = \{1, \ldots,r\}$ and $\G$ be the subset of $\mathbb{Z}^t$ consisting of those  $\textbf{z}$ for which at least one of the following equalities holds for some distinct  $l,k \in \Gamma$ :  
 $$\bsa_l^{\textbf{z}} = \bsa_k^{\textbf{z}}, \quad
 \textbf{c}_l^{\textbf{z}} = \textbf{c}_k^{\textbf{z}}, \quad
 \bsa_l^{\textbf{z}} = \textbf{c}_k^{\textbf{z}}.$$ Finally, set $$h := h^{hom}(\textbf{V}, 1)$$ where $\textbf{V}$ is the vector consisting of the coefficients of $P_1, \cdots, P_r$ and $Q_1, \cdots, Q_r$.  The following results are the applications of Theorem \ref{thmJ} to greatest common divisor of terms of multi-recurrence sequences.

\begin{thm}\label{thmL}
Let $\epsilon >0$ and $$S_0 = \{v \in M_K~ :~ \max_{ijkl} \{|\alpha_{i,j}|_v, |c_{k,l}|_v \} <1\}.$$ Assume  $\G = \{0\}$.   
 We also suppose that $F$ and $G$ are  coprime and simple (i.e. $P_i$, $Q_i$ are constants) and the inequality
 \begin{equation}\label{eq:l_2}
 {\displaystyle \sum_{M_K \setminus S_0}} -\log^{-} \max \{|F(\nbf)|_v, |G(\nbf)|_v\} > \epsilon |\nbf | 
 \end{equation}
 has infinitely many solutions. There exist finitely many  proper subspaces $H_1, \cdots, H_l$ of $\mathbb{R}^t$; constants $C_1, \cdots C_l >0$ and finitely many tuples  $(\theta, \textbf{u}) \in \mathrm{GL}_t(\ZZ) \times \mathbb{Z}^t$  such that all but finitely many solutions of \eqref{eq:l_2}  satisfy  $$\nbf = \theta (\textbf{x}) + \textbf{u} \quad  \mbox{where} \quad \textbf{x} \in \Sc_i$$ with
\begin{equation*}
		\Sc_i = \set{\textbf{x} \in \ZZ^t : \abs{\textbf{x}^{\perp}} \leq C_i \log^+ \abs{\textbf{x}^{H_i}}},
\end{equation*}
for some $i \in \{1, \cdots, l \}$,  where $\textbf{x}^{\perp} $ lies in the orthogonal of $H_i$ and $\textbf{x}^{H_i} \in H_i$ with $\textbf{x} = \textbf{x}^{H_i} + \textbf{x}^{\perp}$.
\end{thm}

\begin{remark}
We emphasize that the conclusion of Theorem \ref{thmL} remains valid for multi-recurrence sequence $F$ and $G$ that are not necessarily simple.
\end{remark}

\begin{thm}\label{thmM}
Let $\epsilon >0$ and $$S_0 = \{v \in M_K~ :~ \max_{ijkl} \{|\alpha_{i,j}|_v, |c_{k,l}|_v \} <1\}.$$   
\begin{compactitem}
 \item[a)] If $F$ and $G$ are  simple (i.e. $P_i$, $Q_i$ are constants) and the $\alpha_{i,j}$ and $c_{i,j}$  are multiplicatively independent.  Then, the number of solutions of the inequality 
 \begin{equation}\label{eq:l_3}
 {\displaystyle \sum_{M_K \setminus S_0}} -\log^{-} \max \{|F(\nbf)|_v, |G(\nbf)|_v\} > \epsilon |\nbf |
 \end{equation}
does not exceed $$ T_1(\epsilon):= t_1 \left( \dfrac{\epsilon}{2t \left(\sum_{i,j} h(\alpha_{i,j}) + h(c_{i,j}) \right)}\right) \frac{|I|!}{2} 2^{35|I|^3}d^{6|I|^2} + s \rho\left( \frac{\epsilon}{2s}\right)$$
with $$|I| = \binom{2rt^2+96r^3t^8\epsilon^{-1}}{2rt^2},$$ where, the functions $\rho(\cdot)$ is defined in Theorem \ref{thmK}, with parameter ($D=0$) and the function $t_1(\cdot)$ is defined in Proposition \ref{pro3}, with the parameters ($h^{hom}(f,g,1)=h,~n = 2rt^2,~ D=t, ~ m= 96t^8r^3$).
  \item[b)] If $t=1$ (i.e. $F, G$ are linear recurrence sequences) and $F,~G$ are  coprime and $\G = \{0\}$, then the number of solutions of inequality \begin{equation}\label{eq:l_4}
  {\displaystyle \sum_{M_K \setminus S_0}} -\log^{-} \max \{|F(n)|_v, |G(n)|_v\} > \epsilon n
  \end{equation}
does not exceed $$ T_2(\epsilon) :=  L + t_1 \left( \dfrac{\epsilon}{3r^2 \max_{i,j}\{h(c_i), h(\alpha_j), 1\}}\right) \exp((7C^m)^{8C^m}) + s \rho\left( \frac{\epsilon}{2s}\right)$$
with $$D= \max \{\deg P_i, \deg Q_i\} + 2c_7 \left(r\max_{i,j} \{h(\alpha_{i}), h(c_j)\} \right)^{2r+1};~ L=  \left(\dfrac{288 D^2(2r+1)^4d^4}{ \epsilon^2} \right)^{3/2}$$ and $$ m= 12 D^2(2r+1)^3\epsilon^{-1}3r^2 \max_{i,j}\{h(c_i), h(\alpha_j), 1\},~C= \binom{2r+1+m}{2r+1}$$ and $$c_7 =  \dfrac{58(2r+1)!(3d)^{2r+2} \log d}{(2r+1)^{2r+1}}.$$ Here, the functions $\rho(\cdot)$ is defined in Theorem \ref{thmK}, with parameters ($D= \max \{\deg P_i, \deg Q_i\} + 1$ and $t=1$) and the function $t_1(\cdot)$ is  defined in Proposition \ref{pro3}, with the parameters ($h^{hom}(f,g,1)=h,~n = 2r+1$ and , $ D, ~ m$ as defined above).
\end{compactitem}
\end{thm}

Theorems \ref{thmL} and \ref{thmM} describe, in broad terms, the algebraic structure of pairs of multi-recurrence sequences  $(F(\nbf), G(\nbf))$ whose terms have a large greatest common divisor. They may be viewed as higher-dimensional analogues of \cite[Theorem 1.11]{L} and \cite[Theorem 1.7]{GW}.

\subsubsection{Effective version of Hadamard quotient theorem and beyond.}

 Let $(F_n)_n$ and $(G_n)_n$ be two integer linear recurrence sequences  and $\mathcal{N}$ be the set of natural integers such that $F_n/G_n$ is an integer. Hadamard quotient theorem, conjectured by Pisot and proven by Van der Poorten states that : if the set  $\mathcal{N}$ coincides with $\mathbb{N}$ then $\left(F_n/G_n\right)_n$ is itself linear recurrence sequence. Corvaja-Zannier  proved that if  $\mathcal{N}$ is infinite set, then there exists a non zero polynomial such that $\left( P(n)F_n/G_n\right)_n$ agrees with a linear recurrence sequence for all $n$ in a suitable arithmetic progression (see \cite[Theorem 2]{CZ5}). More precisely,  We recall the relevant result in a number-field setting. Let $K$ be a number field of degree $d$, and suppose that the sequences admit representations

\begin{equation}\label{eq}
G_n = \sum_{i=1}^{r} Q_i(n)c_{i}^{n}, \quad \mbox{and}\quad F_n = \sum_{i=1}^{r} P_i(n)\alpha_{i}^{n},
\end{equation}
where the $P_i, Q_i$ are polynomials  with coefficients respectively $a_{ij}, b_{ij}$ in $K$ and  $\alpha_{i}, c_i\in K$ for all $i \in \{1, \cdots, r\} $. Let $S$ be the set consisting of all archimedean places and all nonarchimedean places associated with the roots  $\alpha_{i}$ and $c_{i} $  and  $s$ the cardinality of $S$. We denote $h_1= {\rm max}_{i,j,k,l}\{ h( a_{ij}),  h( b_{kl}),1\}$.

\begin{thm}(\cite[Theorem 2]{CZ5})
Assume that the roots  of the characteristic polynomials of $(F_n)_n$ and  $(G_n)_n$ generate together a torsion-free multiplicative group denoted $\mathcal{A}$. We also suppose that for infinitely many $n\in \mathbb{N}$, we have $G(n) \neq 0$ and $F_n/G_n \in \mathcal{O}_S$. Then there exists a nonzero polynomial $P(X) \in K[X]$ such that both sequences $F_n/G_n$ and $G_n/Q(n)$ are linear recurrence sequences.
\end{thm}

Our results extend this theorem in a quantitative direction by providing explicit bounds for the relevant exceptional sets. 

\begin{thm}\label{thmN}
Assume that the roots  of the characteristic polynomials of $(F_n)_n$ and  $(G_n)_n$ generate together a torsion-free multiplicative group denoted $\mathcal{A}$. We also suppose that  the set $\mathcal{R}$ of $n$ such that $G(n) \neq 0$ and $F(n)/G(n) \in \mathcal{O}_S$ contains at least $$C = T_2 \left( \frac{1}{4s(s+1)d^4}\right) + \rho\left( \frac{1}{8sd^4}\right) + (s-1) \rho\left( \frac{1}{4s^2d^4}\right) +1$$ elements, where  $T_2(\cdot)$  is the function defined in Theorem \ref{thmM} b), with $r$ replaced by $l$ defined in the Lemma \ref{lemA"} and $\rho(\cdot)$ is the function defined in Lemma \ref{lemA"}. Then there exists a nonzero polynomial $P(X) \in K[X]$ such that both sequences $(P(n)F_n/G_n)_n$ and $(G_n/Q(n))_n$ are linear recurrence sequences.
\end{thm}

\begin{remark}
If the polynomials $Q_1, \cdots, Q_r$  are pairwise coprime, then the polynomial $P$ in Theorem \ref{thmN} can be taken to be constant. This follows from the uniqueness of the Binet representation \eqref{eq} of the sequence $(G_n)_n$.

The torsion-free hypothesis can be viewed as a normalization condition. Indeed let $q$ be  the order of the torsion in multiplicative group $\mathcal{A}$ generated by the roots of characteristic polynomials of  $(F_n)_n$ and  $(G_n)_n$, then, for each $r=0, \cdots, q-1,$ the sequences $(F_{nq+r})_n$ and $(G_{nq+r})_n$ have characteristic roots that generate a torsion-free multiplicative group. Consequently, Theorem \ref{thmN} can be applied to these subsequences, yielding the corresponding conclusion on a suitable arithmetic progression. 
\end{remark}

 For multi-recurrence sequences, a partial result in this direction was obtained in \cite[Proposition 1.10]{GW}. Adapting the argument used in the proof of Theorem \ref{thmN}, we show that an analogous statement holds for simple multi-recurrence sequences, subject to suitable hypotheses on the multiplicative group generated by their multi-roots. More precisely, we obtain the following result.

\begin{thm}\label{thmO}
Let $F(\nbf)$ and $G(\nbf)$ be two simple multi-recurrence sequences over $K$  their multi-roots $\alpha_{i,j}'$s and $c_{ij}'$s are multiplicatively independent.  Let $S$ the set containing all the prime ideals above their multi-roots and all the archimedean places over $K$, and let $s$ the cardinality of $S$. Assume that the set $\mathcal{R}$ of $\nbf$ such that $G(\nbf) \neq 0$ and $F(\nbf)/G(\nbf) \in \mathcal{O}_S$ contains at least $$\zeta = \rho \left( \frac{1}{8sd^4}\right)+ (s-1) \rho \left( \frac{1}{4s^2d^4}\right) + T_1\left( \frac{1}{4s(s+1)d^4}\right) +1$$ elements, where the functions $\rho(\cdot)$ and $T_1(\cdot)$ are defined in  Theorem \ref{thmK}, and Theorem \ref{thmM}  respectively.  Suppose that there exist  $i, j \in \{1, \cdots, r\}$  and a place $v$ of $K$ such that $$|c_{jk}|_v > |c_{ik}|_v \quad \mbox{for}~\mbox{all}\quad k \in \{1, \cdots, t\}.$$ Moreover we assume that there is $ k \in \{1, \cdots, t\}$ such that $$\dfrac{|c_{jk}|_v}{|c_{ik}|_v}  = \max_{w \in S} \left \{\dfrac{|c_{jk}|_w}{|c_{ik}|_w} \right \}. $$ Then the sequence  $F(\nbf)/G(\nbf)$ is a simple multi-recurrence.
\end{thm}

 \subsection*{Organization of the paper including main idea of the proof.} The paper is organized as follows. In Section \ref{sec1}, we introduce the notation used throughout and recall several fundamental results from effective Diophantine geometry. These include properties of the logarithmic Weil height over number fields, a quantitative version of the Northcott property, and a quantitative form of the Subspace Theorem. We conclude the section with several effective results concerning multi-recurrence sequences, which will be used in the subsequent arguments.

\par Section \ref{sec2} is devoted to the proofs of our results on ideal membership and the construction of generators for modules of solutions to linear systems over $\mathcal{O}_S$. We begin by introducing the notion of a lattice over $\mathcal{O}_S$, which is a principal ideal domain under our standing hypotheses. We then establish a procedure for constructing a basis of a lattice from a basis of a full-rank sublattice. This construction is a key ingredient in the proof of Theorem \ref{thmD}. \\
Our approach to constructing explicit generators for the $\mathcal{O}_S$-module ${\rm Sol{\mathcal{O}_S}}(A)$,
consisting of the solutions to the system of homogeneous linear equations associated to a matrix $A$ with coefficients in $\mathcal{O}_S[X_1, \cdots,X_N]$ is the following. We first use  Aschenbrenner's approach to show the existence of the generators on a localization of  $\mathcal{O}_S$ at a suitable prime ideal.  However, it does not provide sufficient quantitative control of their logarithmic Weil heights in our setting, since the argument relies essentially on a qualitative result of Lang \cite[Chap. V, Theorem 1]{L1}. Our method rely on two key ingredients. First, we establish an effective upper bound for the S-norms of the coefficients appearing in the coordinates of these generators. Second, combining this estimate with the relation between the $S$-norm and the logarithmic Weil height established in Lemma \ref{lemB}, we derive explicit height bounds for these coefficients. We then apply a generalized Bézout identity over $\mathcal{O}_S$ to pass from the generators constructed locally to an explicit generating set over the global ring $\mathcal{O}_S$. This completes the proof of Theorem \ref{thmC}. Finally, Theorem \ref{thmB} follows by combining Lemma \ref{lemI} with Theorem \ref{thmC}.

 Section \ref{sec3} is devoted to the proof of Theorem \ref{thmJ}. The proof proceeds by decomposing the sum on the left-hand side into contributions from places $v \notin S$ and places $v \in S$. We first describe the main ideas underlying the treatment of the former contribution. Our approach provides a quantitative refinement of Levin's method \cite{L}. We begin by dividing the solutions of the inequality under consideration into two classes according to their logarithmic Weil heights: those of small height and those of large height. For the former, an effective version of the Northcott property (Lemma \ref{lemE}) yields an explicit upper bound for the number of solutions. For the latter, we construct, for each relevant place $v\in S$, a basis of the quotient space associated with the ideal $(f,g)_m$ in the polynomial ring over the corresponding completion. This construction allows us to produce a collection of linear forms satisfying the key inequality \eqref{eq:p}. The principal difficulty is that Levin's method does not provide sufficient control over the logarithmic Weil heights of the coefficients of these linear forms. Our main innovation is to establish an explicit height bound for these coefficients. This bound enables us to apply a quantitative version of the Subspace Theorem (Lemma \ref{thm1}) and thereby obtain the desired estimate. The contribution from the places in $S$ is handled by a similar argument.

The last section (sec \ref{sec4}) of this paper is devoted to applications of the results established in section \ref{sec3} to the study of greatest common divisor from terms of multi-recurrence sequences. The key observation is that, for a finitely generated multiplicative subgroup $\Gamma$ of $K^*$, there is a one-to-one correspondence between the set of multi-recurrences with multi-roots in $\Gamma$ and the set of Laurent polynomials. Combining this correspondence with the results of Section \ref{sec3},  we show that the solutions of the inequalities \eqref{eq:l_2}, \eqref{eq:l_3}, and \eqref{eq:l_4} are contained in the solution sets of finitely many exponential Diophantine equations, with an explicit bound for the number of equations required.  We then invoke the quantitative results on multi-recurrence sequences established in Section \ref{sec1}, particularly Lemmas \ref{lemJ}, \ref{thm2}, and \ref{thm3}, to complete the proofs of Theorems \ref{thmL} and \ref{thmM}. Finally, combining Theorems \ref{thmK} and \ref{thmM} with tools from the geometry of numbers, including Minkowski's convex body theorem and results on multiplicative independence, we obtain an effective version of the Hadamard quotient theorem (see Theorem \ref{thmN}).

\section{Preliminaries}\label{sec1}
In this section, we recall several effective results over number fields that will be used throughout the paper.
Let $K$ be a number field of degree $d$ and let $M_{K}$ denote the set of places of $K$. For $v \in M_K$ and $x \in K$, we define the absolute value $\vert x \vert_v$ as follows.
\begin{itemize}
\item $\vert x \vert_v = \vert \sigma (x) \vert ^{1/d}$ if $v$ corresponds to the embedding $\sigma: K \rightarrow \mathbb{R}$,\\
\item $\vert x \vert_v = \vert \sigma (x) \vert^{2/d}$ if $v$ corresponds to the embedding $\sigma: K \rightarrow \mathbb{C}$,\\
\item $\vert x \vert_v = (N(\mathfrak{P}))^{-{\rm ord}_{\mathfrak{P}}(x)/d}$ if $v$ corresponds to the prime ideal $\mathfrak{P}$ of $\mathcal{O}_K$ and ${\rm ord}_{\mathfrak{P}}(x)$ the exponent of $\mathfrak{P}$ in the decomposition of the ideal generated by $x$. 
 \end{itemize}
A place $v$ (real resp. complex) is infinite if $v$ corresponds to an embedding (in $\mathbb{R}$ resp. $\mathbb{C}$) and finite if it corresponds to a non-zero prime ideal of $\mathcal{O}_K$. The definitions are such that the product formula $$ \prod_{ v \in M_K}^{} \vert x \vert_v =1, \mbox{ for } x \in K^*=K \setminus \{0\}$$ holds. Let $S \subset M_K$  be a finite set of cardinality $s$ containing all infinite places of $K$.  The ring of $S$-integers is defined by $$\mathcal{O}_S= \{x \in K: \vert x \vert_v \leq 1 \mbox{ for all } v\notin S\}$$ and its group of units is given by $$\mathcal{O}_S^*= \{x \in K: \vert x \vert_v = 1 \mbox{ for all } v\notin S\}.$$ For $v \in M_K$, we define the quantity $s(v)$ as follows.
\begin{equation*}
s(v) = \left\{
\begin{aligned} 1/d & \mbox{ if } v \mbox{ is real infinite},\\2/d & \mbox{ if } v \mbox{ is complex infinite},\\ 0 & \mbox{ if } v \mbox{ is finite}.
\end{aligned}
\right.
\end{equation*}
By the definition of $s(v)$, we have $$\sum_{v \in S}^{} s(v) = 1.$$ 
For a nonzero vector $\xbf=(x_1, \ldots, x_m) \in K^m$, we define
 \begin{equation*}
\vert \xbf \vert_v = \left\{
\begin{aligned}
\left( \sum_{i=1}^{m} \vert x_i  \vert_v^{2d} \right)^{s(v)/2} & \mbox{ if } v \mbox{ is real infinite}, \\
	 \left( \sum_{i=1}^{m} \vert x_i \vert_v^{d} \right)^{s(v)/2} & \mbox{ if } v \mbox{ is complex infinite},\\
\max \left\{ \vert x_1\vert_v, \ldots, \vert x_m \vert_v \right\} & \mbox{ if } v \mbox{ is finite}.
\end{aligned}
\right.
\end{equation*}
The (multiplicative) Weil height of $\xbf$ is defined as follows: $$H(\xbf)= \prod_{ v \in M_K}^{} \vert \xbf \vert_v.$$ By the product formula, this height is independent of the choice of the number field $K$ containing the coordinates of $\xbf$. For non-zero vector $\xbf=(x_1, \ldots, x_m) \in K^m$,  its homogeneous logarithm height is defined by $$h^{hom}(\xbf) = {\displaystyle \sum_{v \in M_K}} \max \{|x_1|_v, \cdots, |x_n|_v\}.$$  For $\alpha \in K$, we define the  multiplicative Weil height by $$H(\alpha):= \prod_{ v \in M_K}^{} \max \{1, \vert \alpha \vert_v\}$$ and the logarithm  Weil height  by  $$h(\alpha) = {\displaystyle \sum_{v \in M_K}} \log \max \{1, |\alpha|_v\}.$$ We now introduce the notion of "almost $S$-units" which generalizes the notion of $S$-units.

\begin{defin}\label{def}
Let $\eta \geq 0$ be a real number, and  let $S$ be a finite set of places of $K$ containing all archimedean places. An element $u \in K^*$ is called  almost $(S, \eta)$-unit if $$h_{S}(u) := {\displaystyle \sum_{v \notin S }^{}} \max\{0, \log |u|_v\} + \max \left \{ 0, \log \left|\frac{1}{u} \right|_v \right \} \leq \eta h(u).$$ In particular, $u$ is $S$-unit if and only if it is an almost $(S, 0)$-unit.  For $\textbf{u} = (u_1, \cdots, u_n) \in (K^*)^n$, we define $$\max \left \{ 0, \log \left|\frac{1}{\textbf{u}} \right|_v \right \} :=  \max \left \{ 0, \log \left|\frac{1}{u_1} \right|_v, \cdots, \log \left|\frac{1}{u_n} \right|_v \right \}.$$  We then set $$\mathbb{G}_m^n(K)_{S, \eta} := \{ \textbf{u} \in \mathbb{G}_m^n(K)~|~ h_{S}( \textbf{u}) \leq \eta h^{hom}( \textbf{u}, 1) \}.$$
\end{defin}

The height functions defined above satisfy the following standard properties (see, for instance, \cite[Chap.~III]{W} or \cite[Lemma~2.1]{AN}).

\begin{lem}\label{lemA}
For $\eta,\gamma\in \KK^*$, we have
\begin{compactitem}
\item[a)] ${H}(\eta) \geq 1$ and ${H}(\eta) = {H}(\eta^{-1})$,
\item[b)] ${H}(\eta_1 + \cdots + \eta_r) \leq r{H}(\eta_1)\cdots {H}(\eta_r)$,
\item[c)] ${H}(\eta^n) ={H}(\eta)^{\vert n\vert},$ for any $n \in \mathbb{Z}$,
\item[d)] ${H}(\eta\gamma) \leq  H(\eta)H(\gamma)$.
\end{compactitem}
\end{lem}

Let $S$ be a finite set of places of $K$ containing all archimedean places. Let $\mathcal{P}_1, \cdots, \mathcal{P}_t$ denote the prime ideals corresponding to the finite places in $S$ and set $$Q_S = N(\mathcal{P}_1\cdots \mathcal{P}_t),$$ where $N(\cdot)$ denotes the norm function on $K$. For $\alpha \in K^*$, the fractional ideal generated by $\alpha$ can be written as a product of two fractional ideals $\mathcal{A}_1\cdot \mathcal{A}_2$, where $\mathcal{A}_1$ is supported on the prime ideals corresponding to the finite places in $S$, and $\mathcal{A}_2$ is supported on the remaining prime ideals. We define the $S$-norm of $\alpha$ by $N_S(\alpha) := N(\mathcal{A}_2)$. With our definition of absolute values, this gives  $$N_S(\alpha) = {\displaystyle \prod_{v \in S}} \vert \alpha \vert_v.$$ We denote by $R_K$ the regulator of $K$, and set $$c_4 := 87d!d\sqrt{d-1}\log d.$$ The following result is due to Győry and Yu \cite{GY} and is based on Dirichlet's approach to the proof of $S$-unit theorem.

\begin{lem}\label{lemB}
For every nonzero $\alpha \in \mathcal{O}_S$, there exists $\epsilon \in \mathcal{O}_S^*$ such that $$h(\epsilon \alpha) \leq \log N_S(\alpha) + c_4R_K + h_K \log Q_S$$ where $h_K$ and $R_K$ are respectively  the class number and regulator of $K$. If $S= S_{\infty}$, then we replace $h_K \log Q_S$ by $1$. 
\end{lem}

Let $s$ be  the cardinality  of $S$. Dirichlet-Chevalley-Weil $S$-unit theorem, that is , the extension of the classical Dirichlet unit theorem asserts that the group $\mathcal{O}_S^*$ has rank $s-1$.  Let $\{\epsilon_1, \cdots, \epsilon_{s-1}\}$ be a  fundamental system  of $\mathcal{O}_S^*$. Thus, every $\epsilon \in \mathcal{O}_S^*$ admits a unique representation $$\epsilon = \eta \epsilon_1^{b_1} \cdots, \epsilon_{s-1}^{b_{s-1}},$$ where $\eta$ is a root of unity and $b_1, \cdots, b_{s-1}$ are rational integers. Pick $s-1$ places $v_1, \cdots, v_{s-1}$ from $S$. The $S$-regulator is defined by $$R_S = {\rm det}(\log |\epsilon_i|_{v_j})_{i,j}.$$ It is well known that $R_S$ is nonzero and independent of the choice of the places  $v_1, \cdots, v_{s-1}$  and of the chosen fundamental system of $S$-unit, provided the absolute values and the regulator are normalized consistently. Let $S_{\infty}$ be the set of infinite places of $K$. When $S = S_{\infty}$,  the $S$-regulator coincides with the regulator $R_K$ of $K$. Let $\mathcal{P}_1, \cdots, \mathcal{P}_t$  be the prime ideals corresponding to the finite places in $S$. By \cite[Pag. 61]{evertse-gyory-2022}, we have
\begin{equation}\label{eq:i'}  
R_S \leq h_K R_K {\displaystyle \prod_{i=1}^{t}} \log N(\mathcal{P}_i),
\end{equation}
where $h_K$ and $R_K$ denote  the class number and regulator of $K$, respectively. For $s \geq 2$,  define $$c_5 = \dfrac{((s-1)!)^2}{2^{s-2}d^{s-1}}(53 d^2 \log 6d)^{s-2}.$$ The following result provides an upper bound for the heights of a fundamental system of $S$-units.

\begin{lem}\label{lemC}
There exists a fundamental system $\{\epsilon_1, \cdots, \epsilon_{s-1}\}$ of $S$-unit in $K$ such that 
 \begin{equation*}
\max_i h(\epsilon_i) \leq \left\{
\begin{aligned}
c_5 h_K R_K {\displaystyle \prod_{i=1}^{t}} \log N(\mathcal{P}_i) & \mbox{ if } S \neq  S_{\infty},\\
c_5 R_K \quad & \mbox{ if } S = S_{\infty}.	
\end{aligned}
\right.
\end{equation*} 
\end{lem}

\begin{proof}
This lemma follows from \cite[Lemma~1]{BG} combined with the inequality \eqref{eq:i'}. Its proof relies on Dirichlet's approach to the unit theorem and the theorem of Minkowski on the successive minima of symmetric convex bodies in the geometry of numbers.
\end{proof}

\begin{lem}\label{lemD}
Let $\alpha_1, \cdots, \alpha_n$ be the $S$-units and $\Gamma$ be a subgroup generated by  $\alpha_1, \cdots, \alpha_n$ which is torsion-free of rank $m \leq n$. Then there exist $(\gamma_1, \cdots, \gamma_m) \in \Gamma^m$ multiplicatively independent such that $$\Gamma = <\gamma_1, \cdots, \gamma_m>$$ with $$ h(\gamma_i) \leq \dfrac{\max\{2, n\}}{2} \max \{h(\alpha_1), \cdots, h(\alpha_n)\}.$$
\end{lem}

\begin{proof}
For $\epsilon \in \mathcal{O}_S^*$, define  $$\psi(\alpha) := (\log |\alpha|_v)_{v \in S} \in \mathbb{R}^s.$$ By Chevalley-Weil-Dirichlet theorem, the image of $\mathcal{O}_S^*$  under $\psi$  is a full lattice of the hyperplane  $$H := \{(x_1, \cdots, x_s) \in \mathbb{R}^s~:~ x_1 + \cdots + x_s = 0 \}.$$ Consequently, there exist $\alpha_{i_1}, \cdots, \alpha_{i_m} \in \{\alpha_1, \cdots, \alpha_n\}$ such that $\psi(\alpha_{i_1}), \cdots, \psi(\alpha_{i_m})$ are linearly independent over $\mathbb{R}$. Applying \cite[Lemma 8, Chap V]{C} to the  convex distance function $f : \mathbb{R}^s \to \mathbb{R}_{\geq 0}$, define by $$f(\xbf) = \frac{1}{2}(|x_1| + \cdots + |x_s|), \quad  \xbf = (x_1, \cdots, x_s) \in \mathbb{R}^s,$$  we obtain the desired result.
\end{proof}

The following effective version of Northcott's property is due to Evertse \cite{E1}.

\begin{lem}\label{lemE}
Let $\theta$ be a nonzero algebraic number, not necessarily in $K$ and let $C\geq 1$ be a constant. Then the number of $x \in K$ with $h(\theta x) < C$ is at most $5(2C^3)^d$.
\end{lem}

We also need the following lower bound for the height of a product of multiplicatively independent elements of a number field.

\begin{lem}\label{lemF}
Let $\gamma_1, \cdots, \gamma_s \in K^*$ be multiplicatively independent.  Let $(n_1, \cdots n_s) \in \mathbb{N}^s$, and  suppose that $$\gamma_0 = \gamma_1^{n_1} \cdots \gamma_s^{n_s}$$ for some $\gamma_0 \in K$. Then $$h(\gamma_0) \ge \dfrac{\max_i n_i}{c_6 {\displaystyle \prod_{i=1}^{s}} h(\gamma_i)},$$ where 
\begin{equation}\label{eqA}
c_6(s) =  \dfrac{58s!3^sd^{s+1} \log d}{s^s}.
\end{equation}
\end{lem}

\begin{proof}
The result follows by combining \cite[Corollary~2.3]{LM} with the method used by Evertse and Győry in \cite[Corollary~7.5.3]{evertse-gyory-2022}.
\end{proof}

Let $K$ be a number field of degree $d$ and $A$ be a positive real number. An $A$-ceilling is an infinite tuple $(A_v)_{v \in M_K}$ of positive real numbers such that $A_v$ belong to the value group of $\vert \cdot \vert_v$ on $K_v^*$ for $v \in M_K$, $A_v =1$ for almost all $v$ and ${\displaystyle \prod_{v \in M_K} } A_v = A.$ This notion, introduced by Mahler, plays a crucial role in the geometry of numbers.

\begin{lem}\label{lemG}
Let $\Delta_K$ be the determinant of the number field $K$; $A>1$ and $(A_v)_{v \in M_K}$ be an $A$-ceilling. Let $(a_v)_v$ be a family of elements of $K$ such that $| a_v |_v \leq 1$ for all but finitely many $v.$ If $A \geq (d/2)|\Delta_K|^{1/2}$, then there is an $x \in K$ with $$|x-a_v|_v \leq A_v \quad \mbox{for}~~ \mbox{all}~~ v \in M_K.$$
\end{lem}

\begin{proof}
 The result follows from a generalized form of Minkowski's theorem due to Bombieri and Vaaler; see \cite[Lemma~7]{E}.
\end{proof}

The following lemma is also well known (see, for instance, \cite[Lemma~2.4]{AN}).

\begin{lem}\label{lemH}
If $s \geq 1, ~ T> (4s^2)^s$ and $T > x/ (\log x)^s,$ then $$x < 2^sT (\log T)^s.$$
\end{lem}

Let  $K$ be a number field of degree $d$, and $S$ be a finite set of places of $K$  containing all archimedean places and all finite places corresponding to prime ideals that occur in the chosen integral ideal representatives of the ideal classes of the class group of $K$. Let $n \geq 1,$  and $f_1, \cdots, f_m, \in \mathcal{O}_S[x_1, \cdots,x_n]$   be polynomials of total degree at most $l$.  We denote by $(f_1, \cdots, f_m)$ the ideal generated by these polynomials in $\mathcal{O}_S[x_1, \cdots,x_n]$, and $h^{hom}(f_1, \cdots, f_m)$   the homogeneous logarithm  height of the vector formed by all their coefficients. 

\begin{lem}\label{lemI}
If $1 \in (f_1, \cdots, f_m)$, then there exist $g_1, \cdots, g_m \in \mathcal{O}_S[x_1, \cdots,x_n]$ with $$1= f_1g_1 + \cdots + f_mg_m$$ and $${\rm deg} h_i \leq d(3l)^{2^{3n}}(h^{hom}(f_1, \cdots, f_m) +1 )$$ for all $i \in \{1, \cdots, m\}$.
\end{lem}

\begin{proof}
The proof follows that of \cite[Proposition~6.2]{Asc}, with  $R = \mathcal{O}_S$.
\end{proof}

We shall require the following quantitative version of Schmidt's Subspace theorem due to Evertse \cite{E}.

\begin{lem}[Subspace theorem]\label{thm1}
Let $\{L_{1v}, \ldots, L_{rv}\}~(v \in S)$ be a linearly independent set of linear forms in $r$ variables with coefficients in $K$ such that $H(L_{iv}) \leq H$ for $i\in\{1,\ldots,r\},v\in S$. Let $0 <\delta < 1.$ Consider the inequality 
 \begin{equation} \label{eq:3}
\prod_{ v \in S}^{} \prod_{i=1}^{r} \dfrac{\vert L_{iv}(\xbf)\vert_v}{\vert \xbf \vert_v} < \prod_{ v \in S}^{} \vert {\rm det}(L_{1v}, \ldots, L_{rv}) \vert_v H(\xbf)^{-r-\delta} \mbox{ with } \xbf \in K^r.
 \end{equation}
 There are proper linear subspaces $T_1, \ldots, T_{t_1}$ of $K^r,$ with $$t_1 < \left( 2^{60r^2}\delta^{-7r}\right)^{s+1}(3 + \log\log H)$$ such that every solution $\xbf \in K^r$ of inequality \eqref{eq:3}  belongs to $$T_1 \cup \cdots \cup T_{t_1}.$$
 \end{lem}
The bound obtained by Evertse for the number of proper subspaces containing the solutions to \eqref{eq:3} is not the strongest currently available, but it suffices for our purposes. Subsequent improvements were obtained by Evertse and Schlickewei \cite{ES}, and later by Evertse and Ferretti \cite{EF}. Their approaches draw on the geometry of numbers to construct a suitable family of subspaces and estimate its cardinality using a generalized form of Roth's method. In particular, Evertse and Ferretti replaced Schmidt's auxiliary-polynomial construction with the approach of Faltings and Wüstholz \cite{FW}, which is based on Faltings' product theorem.

The following result, due to Schmidt \cite{S}, plays a crucial role in the proof of Theorem~\ref{thmM}(b). It also provides an answer to a folklore conjecture sometimes attributed to Mordell.

\begin{lem}\label{lemJ}
Suppose that $\{u_n\}_{n \in \mathbb{Z}}$ is a non-degenerate linear recurrence sequence of complex numbers, whose characteristic polynomial has $k$ distinct roots of  multiplicity at most $a$. Then the number of solutions $n \in Z$ of the equation $u_n=0$ can be bounded above by $$c(k, a) = \exp((7k^a)^{8k^a}).$$
\end{lem}

Amoroso and Viada \cite{AV} obtained the strongest bound currently available in the literature. Nevertheless, Schmidt's bound stated above suffices for our purposes.

We now recall a result of Ahlgren \cite{ahlgren-1999}, which provides a quantitative version of Laurent's theorem \cite{laurent-1989}. In this sense, it can be viewed as a quantitative generalization of the Skolem–Mahler–Lech theorem to multivariate recurrence sequences.

Consider an equation of the form
	\begin{equation}
		\label{eq:genmultrec}
		\sum_{l=1}^{k} P_l(\textbf{x}) \boldsymbol{\alpha}_l^{\textbf{x}} = 0
	\end{equation} 
	in variables $\textbf{x} = (x_1, \cdots, x_m) \in \ZZ^m$ where $P_l$ are polynomials with coefficients in a number field $K$, and
	\begin{equation*}
		\boldsymbol{\alpha}_l^{\textbf{x}} = \alpha_{l,1}^{x_1} \cdot \ldots \cdot \alpha_{l,m}^{x_m}
	\end{equation*}
	with $\alpha_{l,j} \in K^* $ for $ 1 \leq l \leq k $ and $ 1 \leq j \leq m $.
	We put $\Gamma = \set{1,\ldots,k}$ and denote by $\Pc$ a partition of the set $\Gamma$.
	A subset $\lambda \subseteq \Gamma$ occurring in the partition $\Pc$ will be considered as element of $\Pc$.
	Given $\Pc$, consider the system of equations
	\begin{equation}\label{eq:Z2}
		\sum_{l \in \lambda} P_l(\textbf{x}) \boldsymbol{\alpha}_l^{\textbf{x}} = 0 \qquad (\lambda \in \Pc).
	\end{equation} 
	A solution $\textbf{x}$ of \eqref{eq:Z2} is called $\Pc$-degenerate  if, in at least one equation of the system, a nonempty proper subsum vanishes. Otherwise  $\textbf{x}$ is  call $\Pc$-non-degenerate.
	Let $M(\Pc)$ denote the set of $\Pc$-non-degenerate solutions of the system \eqref{eq:Z2}.
	Let $G(\Pc)$ be the subgroup of $\ZZ^m$ consisting of $\textbf{x}$ such that
	\begin{equation*}
		\boldsymbol{\alpha}_l^{\textbf{x}} = \boldsymbol{\alpha}_t^{\textbf{x}}
	\end{equation*}
	whenever $l$ and $t$ belong to the same set $\lambda$ of $\Pc$.
	Let $d$ be the degree of $K$. For $l \in \Gamma$, denote by $\delta_l$ the total degree of $P_l$.
	Set
	\begin{equation*}
		A = \sum_{l \in \Gamma} \binom{m+\delta_l}{m} \qquad \text{and} \qquad D= \max(m, A).
	\end{equation*}
	We equip $\RR^m$  with the Euclidean norm $\abs{\cdot}$,  and define $\log^+ z := \max\setb{\log z, 1}$.
	Let $r$ be the rank of $G(\Pc)$, and $H$  be the real vector subspace spanned by  $G(\Pc)$.
	Schlickewei and Schmidt \cite{SS} proved the following result.
	
	\begin{lem}\label{thm2}
		\label{thm:schsch}
		If $G(\Pc) = \set{0}$, then
		\begin{equation*}
			\# M(\mathcal{P}) < 2^{35D^3}d^{6D^2}.
		\end{equation*}
	\end{lem}
	
	Ahlgren's result \cite{ahlgren-1999} addresses the case in which the group $ G(\Pc)$ is non-trivial.
	More precisely, define 
	\begin{equation*}
		\Sc = \set{\textbf{x} \in \ZZ^m : \abs{\textbf{x}^{\perp}} \leq c_0 \log^+ \abs{\textbf{x}^{H}}}
	\end{equation*}
	where $\textbf{x}^{\perp} \in H^{\perp}$  and $\textbf{x}^{H} \in H$ are the orthogonal projections of $\xbf$ onto $H$  and its orthogonal complement $ H^{\perp}$, respectively so that  $$\textbf{x} = \textbf{x}^{H} + \textbf{x}^{\perp}.$$
	Let $\mathrm{GL}_m(\ZZ)$ denote the group of $m \times m$  integer matrices with determinant $\pm 1$.
	Set $c_0= D^{12}2^{4D+12}d^6$.
	Ahlgren's result can then be stated as follows.
	\begin{lem}\label{thm3}
		\label{thm:ahlgr}
		With the notation above, we have
		\begin{equation*}
			M(\Pc) \subseteq \bigcup_{\psi \in \T} \psi(\Sc)
		\end{equation*}
		for $ \T \subseteq \GG$ with $ |\T| \leq 2^{36D^3}d^{6D^2}$ and $ \GG $ the group of transformations of $\RR^m$ having the form $\psi(\textbf{x}) = \theta (\textbf{x}) + \textbf{u}$ where $\theta \in \mathrm{GL}_m(\ZZ)$ fixes $H$ pointwise and $\textbf{u} \in \ZZ^m$.
	\end{lem}
	
\section{Effective version of Aschenbrenner's method.}\label{sec2}

We begin by introducing the notion of a lattice over a principal ideal domain (PID) and establishing some fundamental properties needed in the proof of Theorem~\ref{thmD}. These results generalize \cite[Theorem~I, Chap.~I]{C} and \cite[Lemma~8, Chap.~V]{C}.

\subsection{Explicit bases of lattices over a PID}

Let $R$ be a PID, let $L$ be its field of fractions, and let $V$ be an $n$-dimensional vector space over $L$. Let $M \subset V$ be an $R$-module that is torsion-free.  We call $M$ an $R$-lattice in $V$  if it is finitely generated and contains an $L$-basis of $V$.  By \cite[Corollary~2]{PS}, the module $M$ admits an $R$-basis consisting of $n$ elements.

\begin{lem}\label{lema}
Let $M$ and $N$ be two lattices over $R$ such that $N \subseteq M$.
\begin{compactitem}
\item[i)] If the vectors $a_1, \cdots, a_n$ for a basis of $N$, then there exists a basis $b_1, \cdots, b_n$ of $M$ such that 
\begin{equation}\label{eq:l}
\left\{
\begin{aligned}
a_1 &= \mu_{11}b_1, \\
a_2 &= \mu_{21}b_1 + \mu_{22}b_2,\\
     \vdots \\     
a_n &= \mu_{n1}b_1 + \mu_{n2}b_2 + \cdots + \mu_{nn}b_n,
\end{aligned}
\right.
\end{equation}
where $\mu_{ij} \in R$ for $i,j \in \{1, \cdots, n\}$ and $\mu_{ii} \neq 0$ for all $i$.
\item[ii)] Conversely, to every basis $\{b_1, \cdots, b_n\}$ of $M$ there exists a basis $\{a_1, \cdots, a_n\}$ of $N$ such that \eqref{eq:l} holds.
\end{compactitem}
\end{lem}

\begin{proof}
The case $R= \mathbb{Z}$ was proved by Cassel \cite{C}, using the fact that $\mathbb{Z}$ is an Euclidean domain. However Cassel's approach may not apply in  our setting, since  $R$ is not necessary a Euclidean domain. We first proof  $~~ii)$. Let $\{b_1, \cdots, b_n\}$ be a basis of $M$. For $i \in \{1, \cdots, n\}$, let $Pr_i : M \to R$ denote  the $i$-th coordinate projection. Thus, if  $b \in M$ is written as $b = \alpha_1b_1 + \cdots + \alpha_nb_n,$ then $Pr_i (b)= \alpha_i$. By the definition of an $R$-lattice,  there exists a nonzero element $\delta \in R$ such that $\delta M \subseteq N$. For each $i\in \{1, \cdots, n\}$, define
$$\Lambda_i = \left\{ <Pr_i(a)>, ~~a \in N ~~ \mbox{such}~\mbox{that} ~~ a = {\displaystyle \sum_{j=1}^{i}} \mu_{ij}b_j \right\},$$ where $< \alpha >$ denote the ideal of $R$ generated by $\alpha$. Since $\delta b_i \in N$, we have $\Lambda_i \neq \emptyset .$ By \cite[Corollary of Theorem 1, Pag. 20]{PS},  $\Lambda_i$ has a maximal element, which we denote by $\mathcal{M}_i = <Pr_i(a_i)>$. We claim that $\{a_1, \cdots, a_n\}$ is a basis of $N$, which yields the desired result. Indeed, by construction, the vectors $a_1, \cdots, a_n$ are linearly independent over $R$. Suppose, for instance, that there exists an element of $N$ that cannot be expressed as an $R$-linear combination of  $a_1, \cdots, a_n$. Let $\mathcal{M}$ denote the set of all such elements, and choose the smallest integer $k$ for which some  $c \in \mathcal{M}$ admits an expression 
\begin{equation}\label{eq:a}
c = c_1b_1 + \cdots + c_kb_k, \quad \mbox{with} \quad c_1, \cdots,c_k \in R.
\end{equation}
Put $d= \gcd (\mu_{kk}, c_k)$. By the Euclidean algorithm, there exist $r, t \in R$ such that $d= r\mu_{kk}+ tc_k$. Consequently
 $$
\begin{array}{lcl}
tc+ra_k  &= & (tc_1+r\mu_{k1})b_1 + (tc_2+r\mu_{k2})b_2 + \cdots + (tc_k+r\mu_{kk})b_k\\\\
   
                                    &=& (tc_1+r\mu_{k1})b_1 + (tc_2+r\mu_{k2})b_2 + \cdots + db_k.
 \end{array} 
 $$
Since $tc+ra_k \in N$, the ideal $<d>$ belongs to  $\Lambda_k$. Moreover $d$ divides $\mu_{kk}$, so  $ <\mu_{kk}> \subseteq <d>$. The maximality of  $\mathcal{M}_k = <\mu_{kk}>$  therefore implies  $<d> = <\mu_{kk}>$. It follows $c_k \in <\mu_{kk}>$, and hence  there exists $\gamma \in R$ such that $c_k = \gamma \mu_{kk}$. Subtracting  $\gamma a_k$ from $c$, we obtain $$ c -\gamma a_k = {\displaystyle \sum_{j=1}^{k-1}}tjb_j, ~~ \mbox{where}~~ t_1, \cdots, t_{k-1} \in R.$$ Since  $ c -\gamma a_k \in N$, the minimality of $k$ implies $c-\gamma a_k  \notin \mathcal{M}$. Thus   $c-\gamma a_k$ is an $R$-linear combination of $a_1, \cdots, a_n$. Since $a_k$ is also among those vectors, the same must hold for $c$, contradicting  $c \in \mathcal{M}$.This contradiction completes the proof of   $ii)$.\\

Let $\{a_1, \cdots, a_n\}$ be a fixed basis of $N$. Since $\delta M$ is a sub-lattices of $N$, by $ii)$,  there exist a basis $\{\delta b_1, \cdots, \delta b_n \}$ of $\delta M$ such that 
\begin{equation*}
\left\{
\begin{aligned}
\delta b_1 &= w_{11}a_1, \\
\delta b_2 &= w_{21}a_1 + w_{22}a_2,\\
     \vdots \\     
\delta b_n &= w_{n1}a_1 + w_{n2}a_2 + \cdots + w_{nn}a_n,
\end{aligned}
\right.
\end{equation*}
with $w_{ij} \in R$ and $w_{ii} \neq 0$. Clearly $\{b_1, \cdots, b_n\}$ is a basis of $M$. Since $a_i \in M$ for every $i \in \{1, \cdots, n\}$, each $a_i$ admits a unique expression as an $R$-linear combination of $b_1, \cdots, b_n$. Thus, the desired result follows from the fact that, every $m \in M$  can be uniquely written as $$m = m_1b_1 + \cdots + m_nb_n,\quad \mbox{with} \quad m_1, \cdots,m_n \in R.$$
\end{proof}

Let $K$ be a number field of degree $d$, and let $S$ be a finite set of places of $K$ containing all archimedean places and all prime ideals dividing a fixed collection of integral ideals representing the ideal classes element of the class group of $K$. By \cite[Proposition 5.3.6, Chap. V]{bombieri-gubler-2006},  the ring of $S$-integers  $\mathcal{O}_S$ is a principal ideal domain. Let $s$ be the cardinality of $S$. The following lemma may be viewed as a slight generalization of \cite[Lemma 8, Chap. V]{C}.

\begin{lem}\label{lemb}
Let $M$ be a lattice over $\mathcal{O}_S$ and $a_1, \cdots, a_n$ be the vectors of $M$ which are linearly independent over $\mathcal{O}_S$. Then there exists a basis $\{c_1, \cdots, c_n\}$ of $M$ such that 
\begin{equation*}
\left\{
\begin{aligned}
 c_1 &= \alpha_{11}a_1, \\
c_2 &= \alpha_{21}a_1 + \alpha_{22}a_2,\\
     \vdots \\     
 c_n &= \alpha_{n1}a_1 + \alpha_{n2}a_2 + \cdots + \alpha_{nn}a_n,
\end{aligned}
\right.
\end{equation*}
with; if $s \geq 2$ then $$  \log |\alpha_{ij}|_v  \leq {\rm max} \left \{ {\displaystyle \sum_{i=1}^{s-1}} \max \{0, \log |\epsilon_i|_v\}, ~\log |d \Delta_K|_v \right \},$$ for~ all $i,j \in \{1, \cdots, n\}$ and $ v \in S$, where  $\{\epsilon_1, \cdots, \epsilon_{s-1}\}$ is the fundamental system issue from  Lemma \ref{lemC}, and $\Delta_K$ denote the discriminant of $K$. Otherwise, we have  $$  \log |\alpha_{ij}|_v  \leq {\rm max} \left \{ 0, ~\log |d \Delta_K|_v \right \},$$ for~ all $i,j \in \{1, \cdots, n\}$ and $ v \in S$.
\end{lem}

\begin{proof}
By Lemma \ref{lema}, there exists a basis  $\{b_1, \cdots, b_n\}$ of $M$ satisfying \eqref{eq:l}. Solving \eqref{eq:l} for $b_i$, we obtain 
\begin{equation}\label{eq:b}
b_i =  \mu_{ii}^{-1}a_i + k_{i, i-1}a_{i-1} + \cdots + k_{i,1}a_1,
\end{equation}
where $k_{i, j} \in K$ for all $i,j \in \{1, 2, \cdots, n\}.$ If $\mu_{ii} \in \mathcal{O}_S^*$, we set $c_i=b_i$. Otherwise, we have $N_S(\mu_{ii})>1$. If $s=1$, then $| \mu_{ii}^{-1}|_v = N_S( \mu_{ii}^{-1}) < 1.$ We now define
\begin{equation*}
A_v = \left\{
\begin{aligned}
\left| d \Delta_K\right|_v & \mbox{ if } v \mbox{ is  infinite}, \\
1 ~~~~~~~~~& \mbox{ if } v \mbox{ is finite}.
\end{aligned}
\right.
\end{equation*}
It is readily verified that $$ {\displaystyle\prod_{v\in M_K}} A_v = d |\Delta_K| \geq \frac{d}{2} |\Delta_K|^{1/2}.$$ For fixed $i,j \in \{1, \cdots, n\}$ we apply Lemma \ref{lemG} to the Adele $(a_v)_v$ defines as follows:
\begin{equation*}
a_v = \left\{
\begin{aligned}
k_{i,j}~~~~~~~ & \mbox{ if } v \in S, \\
0 ~~~~~~~~~& \mbox{ if } v \notin S.
\end{aligned}
\right.
\end{equation*}
It follows that, for each $i,j$, there exists   $x_{i,j} \in \mathcal{O}_S$ such that $$| k_{i,j} - x_{i,j}|_v \leq A_v \quad \mbox{for}~~\mbox{all} \quad v \in S.$$ Since $\{b_1, b_2, \cdots, b_n\}$ is a basis of $M$,  the vectors $$c_i = b_i + x_{i, i-1}a_{i-1} + \cdots + x_{i,1}a_1, \quad \mbox{for} \quad i=1, \cdots,n$$ also form  a basis of $M$. Combining this expression  with \eqref{eq:b}, we obtain $$c_i = \mu_{ii}^{-1}a_i + (k_{i, i-1}- x_{i,i-1})a_{i-1} + \cdots + (k_{i,1}- x_{i,1})a_1.$$ Consequently,   $$| k_{i,j} - x_{i,j}|_v \leq  |d \Delta_K|_v$$ and  the desired result follows. We now assume that  $s \geq 2$. Let $\{\epsilon_1, \cdots, \epsilon_{s-1}\}$ be a fundamental system of $S$-unit as in Lemma \ref{lemC}. By the Dirichlet-Chevalley-Weil $S$-unit theorem, the vectors $$(\log |\epsilon_i|_v)_{v\in S} \quad \mbox{with} \quad i =1, \cdots, s-1$$ form a basis of the hyperplane  given by $$H := \{ (x_1, \cdots, x_s) \in \mathbb{R}^s ~:~ x_1+\cdots +x_s =0\}.$$ Since  $$(\log N_S(\mu_{ii}^{-1})^{-s(v)}|\mu_{ii}^{-1}|_v)_{v\in S} \in H, $$ there exist an $S$-unit $\epsilon $ and  real numbers $t_1, \cdots, t_{s-1} \in [0,1)$ such that 
$$\log \left( N_S(\mu_{ii}^{-1})^{-s(v)}|\epsilon\mu_{ii}^{-1}|_v\right) =  t_1\log |\epsilon_1|_v + \cdots + t_{s-1}\log |\epsilon_{s-1}|_v \quad \mbox{for}~~ \mbox{all} \quad v \in S.$$ Using  $N_S(\mu_{ii}^{-1}) < 1$, we deduce that $$ \log |\epsilon\mu_{ii}^{-1}|_v  \leq {\displaystyle \sum_{i=1}^{s-1}} \max \{0, \log |\epsilon_i|_v\}$$ for all $v \in S$. Applying the same argument as above, with $b_i$ replaced by $\epsilon b_i$, we obtain a basis of $M$ whose elements are given by  $$c_i = \epsilon \mu_{ii}^{-1}a_i + (k_{i, i-1}- x_{i,i-1})a_{i-1} + \cdots + (k_{i,1}- x_{i,1})a_1,$$ where  $$| k_{i,j} - x_{i,j}|_v \leq  |d \Delta_K|_v.$$
This completes the proof.\\
 \end{proof}

\subsection{Height bound of generators.}
Let $m,n$ be two positive integers and set $a=\max\{m, n\}$. Let $S$ be a finite set of places of $K$ containing all archimedean places and all nonarchimedean places  corresponding to the  prime ideals dividing a fixed collection of integral ideals representing the ideal classes element of the class group of $K$. We denote the cardinality of $S$ by $s$. Let $\mathcal{P}_1, \cdots, \mathcal{P}_t$ denote the prime ideals corresponding to the finite places in $S$. Let $A$ be an $m \times n$- matrix with entries  $a_{ij} \in \mathcal{O}_S$. We define $${\rm Sol_{\mathcal{O}_S}}(A) = \{X \in \mathcal{O}_S^n~~|~~ AX=0\}$$ and $h= h^{hom}(1, a_{ij}~|~ ~ (i,j) \in \{1, \cdots, m\} \times \{1, \cdots, n\})$. Clearly, ${\rm Sol_{\mathcal{O}_S}}(A)$ is a finite generated $\mathcal{O}_S$-module.  For  convenience, we restate the Theorem \ref{thmC}.

\begin{lem}\label{lemed}(Theorem \ref{thmC})
There exists a generator $\{w^{(1)}, \cdots, w^{(k)}\}$   of  $~~{\rm Sol_{\mathcal{O}_S}}(A)$ such that  the logarithm Weil height of any coordinate of the vectors $w^{(1)}$ is bounded above by $$ a \left( \frac{5}{2}a^2(\log a + h)+ c_4R_K + h_K \log Q_S \right),$$
where $R_K, h_K$ are respectively the regulator and the class number of $K$ and 
$$c_4:=87d!\sqrt{d-1}d\log d, \quad Q_S = N(\mathcal{P}_1\cdots \mathcal{P}_t).$$
 If $S= S_{\infty}$, then we replace $ h_K \log Q_S$ by $1$.
  
\end{lem}

\begin{proof}
Our technique provides a quantitative refinement of Aschenbrenner's in \cite{Asc} which was first introduced by Hermann in \cite{H}.
Let $r = {\rm rank_K}(A) \geq 1$ denote the rank of $A$, and let $\Delta$ be an $r \times r$-submatrix of $A$ such that $\delta = {\det}(\Delta) \neq 0$. Without loss of generality, we may assume that $\Delta$ occupies the upper-left corner of $A$, that is,  $\Delta = (a_{ij})_{1 \leq i,j \leq r}$. Clearly an element $(X_1, \cdots, X_n)$ of ${\rm Sol_{\mathcal{O}_S}}(A)$ is solution of the system
\begin{equation}\label{eq:A2}
\begin{pmatrix}
a_{11} & a_{12} & \cdots & a_{1n} \\
a_{21} &a_{22} & \cdots &a_{2n} \\
\vdots & \vdots & \ddots & \vdots\\
a_{r1} &a_{r2}& \cdots & a_{rn} \\
\end{pmatrix}  \begin{pmatrix}
X_1  \\
X_2  \\
\vdots  \\
X_n  \\
\end{pmatrix} =  \begin{pmatrix}
0  \\
0  \\
\vdots  \\
0  \\
\end{pmatrix}.
\end{equation}
By  Hadamard's inequality for the determinant, we have $$h(\delta) \leq \frac{r}{2} \log r +rh \leq r(\log r + h).$$ Multiplying both side of  \eqref{eq:A2} by the adjoint matrix of $\Delta$, we obtain
\begin{equation}\label{eq:d}
\begin{pmatrix}
\delta & 0 & \cdots &0& c_{1,r+1}& \cdots &c_{1,n} \\
0 & \delta &  \cdots &0& c_{2,r+1}& \cdots &c_{2,n} \\
\vdots & \vdots & \ddots & \vdots\\
0 &0 &\cdots& \delta &  c_{r,r+1}& \cdots &c_{r,n} \\
\end{pmatrix}  \begin{pmatrix}
X_1  \\
X_2  \\
\vdots  \\
X_n  \\
\end{pmatrix} =  \begin{pmatrix}
0  \\
0  \\
\vdots  \\
0  \\
\end{pmatrix}
\end{equation}
with 
\begin{equation*}\begin{split}
h(c_{ij})
&\leq h \left( {\displaystyle \sum_{k=1}^{r}} \delta_{ik}a_{kj}\right) \leq r \max_k h(\delta_{ik} a_{kj}) + \log r \\
&\leq r(r-1)( \log (r-1) + h ) + rh + \log r  \\
&\leq r^2( \log r + h ),
\end{split}\end{equation*}
where $\delta_{ij}$ denotes the determinant of the matrix obtained by deleting the 
$i$-th row and the $j$-th column of  $\Delta$. All the above inequalities follow from the definition of  $\Delta$
 and the Lemma \ref{lemA}. Notice that the following $n-r$ vectors are linearly independent solutions of \eqref{eq:d}.
\begin{equation}\label{eq:c}
v_{\Delta}^{(1)}= \begin{pmatrix}
-c_{1,r+1} \\
\vdots \\
 -c_{r,r+1} \\
 \delta \\
 0\\
\vdots  \\
0  \\
\end{pmatrix}, v_{\Delta}^{(2)}= \begin{pmatrix}
-c_{1,r+2} \\
\vdots \\
 -c_{r,r+2} \\
 0\\
 \delta \\
\vdots  \\
0  \\
\end{pmatrix}~~, \cdots, ~~v_{\Delta}^{(n-r)}= \begin{pmatrix}
-c_{1,n} \\
\vdots \\
 -c_{r,n} \\
 0\\
\vdots  \\
0  \\
\delta \\
\end{pmatrix}.
\end{equation}
It is not difficult to see that 
\begin{equation}\label{eq:e}
\delta \cdot {\rm Sol_{\mathcal{O}_S}}(A) \subseteq \mathcal{O}_S v^{(1)} + \cdots + \mathcal{O}_S v^{(n-r)}.
\end{equation}
If $\delta \in \mathcal{O}_S^*$, then these vectors form a basis of ${\rm Sol_{\mathcal{O}_S}}(A)$, and the desired result follows.    Let $\mathcal{P}$ be a prime ideal dividing $\delta$. Since $\mathcal{O}_S$ is a principal ideal domain, we can write $\mathcal{P} = <t>$ for some $t \in \mathcal{O}_S$. Let $$\mu:= \nu_{\mathcal{P}}(\delta'),$$  where the minimum is taken over  all $r\times r$ submatrices of $A$ with nonzero determinant.  Choose a $r\times r$-submatrix $\Delta_{\mathcal{P}}$ of $A$  attaining this minimum, and put $\delta' = {\rm det} \Delta_{\mathcal{P}}.$ Since up to a sign, each $c_{ij}$ is determinant of a $r\times r$ submatrix of $A$,  the minimality of $\mu$  implies that $$ t^{-\mu}v_{\Delta_{\mathcal{P}}}^{(i)} \in {\rm Sol_{\mathcal{O}_S}}(A) \quad \mbox{for}~~\mbox{all}\quad i= 1, \cdots, n,$$ where the vectors $v_{\Delta_{\mathcal{P}}}^{(i)}$ is defined analogously to $v_{\Delta}^{(i)}$. Applying Lemma  \ref{lemB}, we deduce that there exists $\epsilon \in \mathcal{O}_S^*$ such that $$h(\epsilon t^{\mu}) \leq \log N_S(t^{\mu}) + c_4R_K + h_K \log Q_S. $$ Since $t^{\mu}$ divides $\delta$ over $\mathcal{O}_S$, we obtain  $$\log N_S(t^{\mu}) \leq \log N_S(\delta) \leq h(\delta) \leq r(\log r +h).$$ Consequently  $$h(\epsilon^{-1} t^{-\mu}) \leq r(\log r +h) + c_4R_K + h_K \log Q_S.$$ Setting $$\delta_{\mathcal{P}} = \delta' t^{-\mu}, \quad b_{\mathcal{P}}^{(i)} = \epsilon^{-1} t^{-\mu} v_{\Delta_{\mathcal{P}}}^{(i)} \quad \mbox{for}~~\mbox{ all}  \quad i= 1, \cdots, n, $$ we obtain  
\begin{equation}\label{eq:f}
\delta_{\mathcal{P}} \cdot {\rm Sol_{\mathcal{O}_S}}(A) \subseteq \mathcal{O}_Sb_{\mathcal{P}}^{(1)} + \cdots + \mathcal{O}_S b_{\mathcal{P}}^{(n-r)}
\end{equation}
with
\begin{equation*}\begin{split}
h(b_{\mathcal{P}}^{(i)})
&\leq \log n + n \left( 2r^2(\log r + h)+ c_4R_K + h_K \log Q_S \right)\\
&\leq a \left( \frac{5}{2}a^2(\log a + h)+ c_4R_K + h_K \log Q_S \right),
\end{split}\end{equation*}
for all $i \in \{1, \cdots, n-r\}$. Let $\{\mathcal{P}_1, \cdots \mathcal{P}_M \}$ be the set of all prime ideals above $\delta$. Clearly $M \leq 2 h(\delta)$. By construction of  the $\delta_{\mathcal{P}_i}$'s as above, we get $${\rm gcd}(\delta, \delta_{\mathcal{P}_1}, \cdots , \delta_{\mathcal{P}_M}) = 1$$ which implies $$1 = \delta l_0 + \delta_{\mathcal{P}_1} l_1 + \cdots + \delta_{\mathcal{P}_M} l_M $$ for some $l_0, \cdots, l_M \in \mathcal{O}_S$. Together with  \eqref{eq:e} and \eqref{eq:f}, we deduce that the set $\{v_{\Delta}^{(1)}, \cdots, v_{\Delta}^{(n-r)},  b_{\mathcal{P}_1}^{(1)}, \cdots , b_{\mathcal{P}_M}^{(n-r)}\}$ forms a systems of generators of ${\rm Sol_{\mathcal{O}_S}}(A)$. This completes the proof.
\end{proof}

Let $m,n$ be two positive integers and $a=\max\{m, n\}$.  Let $S$ be a finite set of places of $K$ containing all archimedean places and all nonarchimedean places  corresponding to the  prime ideals dividing a fixed collection of integral ideals representing the ideal classes element of the class group of $K$. Let $\mathcal{P}_1, \cdots, \mathcal{P}_t$ denote the prime ideals corresponding to the finite places in $S$ and set $$Q_S = N(\mathcal{P}_1\cdots \mathcal{P}_t),$$ where $N(\cdot)$ denotes the norm function on $K$. Let $S_{\infty}$ be the set of all archimedean places of $K$. Finally, set $$c_4:=87d!\sqrt{d-1}d\log d.$$ We set $D= \mathcal{O}_S[X_1, \cdots, X_N]$. Let $A$ be an $m \times n$ matrix with entries  $a_{ij} \in D$. We denote by ${\rm deg}(A)$ the maximum of the total degrees of the entries $a_{ij}$ and set $$T:=c_4R_K + h_K \log Q_S.$$ Furthermore,  We define $${\rm Sol_{D}}(A) = \{X \in D^n~~|~~ AX=0\}$$ and put $h= h^{hom}(1, a_{ij}~|~ ~ (i,j) \in \{1, \cdots, m\} \times \{1, \cdots, n\})$. Clearly ${\rm Sol_{\mathcal{O}_S}}(A)$ is a finite generated $D$-module. 

\begin{pro}\label{lemee} 
The $D$-module  ${\rm Sol_{D}}(A)$ is generated by the vectors whose the total degree of each coordinates is bounded by $$\gamma := (2 a\cdot {\rm \deg}(A))^{2((N+1)^N-1)}$$  and the height  by $$C := \frac{5}{2} \left(a \binom{N+ \deg(A)+ \gamma}{N} \right)^3 \log \left(a\binom{N+ \deg(A)+ \gamma}{N} \right) \left(h +T +1 \right).$$ If $S= S_{\infty}$, then we replace $T$ by $c_4R_K +1$.
\end{pro}

\begin{proof}
By  \cite[Proposition 4.3]{Asc}, the $D$-module  ${\rm Sol_{D}}(A)$ is generated by the  vectors of total degree  $\gamma$. Let $(Y_1, \cdots, Y_n)$ be one of such vectors.  We may write $$Y_i = {\displaystyle \sum_{|\nu| \leq \gamma}} Y_{i, \nu}X^{\nu} \quad \mbox{with} \quad Y_{i, \nu} \in \mathcal{O}_S,$$ where $\nu= (\nu_1, \cdots, \nu_N) \in\mathbb{N}^m$, $|\nu| = \nu_1 + \cdots + \nu_N$ and $X^{\nu} = X_1^{\nu_1} \cdots X_N^{\nu_N}$. We also set  $$a_{ij} = {\displaystyle \sum_{|\mu| \leq \deg(A)}} a_{ij, \mu}X^{\mu} \quad \mbox{with} \quad  a_{ij, \mu} \in \mathcal{O}_S.$$
It is well known that a polynomial in $N$ variables of total degree at most ${\rm deg}(A)$ contains at most $\binom{N+ \deg(A)}{N}$ monomials. Since $(Y_1, \cdots, Y_n) \in {\rm Sol_{D}}(A)$, we deduce that the  coefficients $ Y_{i, \nu}$ satisfy the following system of $n\binom{N+ \deg(A)+ \gamma}{N}$ homogeneous equations 
\begin{equation}\label{eq:g}
{\displaystyle \sum_{j=1}^{n}\sum_{\mu + \nu = \lambda}} a_{ij, \mu}Y_{j, \nu} = 0 \quad (|\lambda | \leq \gamma + \deg(A))
\end{equation}
in $m\cdot \binom{N+ \gamma}{N}$  variables $Y_{j, \nu}$ in $\mathcal{O}_S$. Consequently, the system of homogeneous equations 
\begin{equation*}
\begin{pmatrix}
a_{11} & a_{12} & \cdots & a_{1n} \\
a_{21} &a_{22} & \cdots &a_{2n} \\
\vdots & \vdots & \ddots & \vdots\\
a_{m1} &a_{m2}& \cdots & a_{mn} \\
\end{pmatrix}  \begin{pmatrix}
Y_1  \\
Y_2  \\
\vdots  \\
Y_n  \\
\end{pmatrix} =  \begin{pmatrix}
0  \\
0  \\
\vdots  \\
0  \\
\end{pmatrix}.
\end{equation*}
may be replace by the homogeneous system \eqref{eq:g} of $$ n \cdot\binom{N+ \deg(A)+ \gamma}{N} \quad \mbox{equations}~~\mbox{in}\quad m\cdot \binom{N+ \gamma}{N}$$  variables  $Y_{j, \nu}$ in $\mathcal{O}_S$. Applying Lemma \ref{lemed}, we deduce that the height of the coordinates of a generator of  ${\rm Sol_{D}}(A)$ is bounded above by 
$$a \binom{N+ \deg(A)+ \gamma}{N} \left(\frac{5}{2}a^2\binom{N+ \deg(A)+ \gamma}{N}^2\left(\log \left(a\binom{N+ \deg(A)+ \gamma}{N} \right)+ h \right) + T\right).$$ This quantity is further bounded above by  $$\frac{5}{2} \left(a \binom{N+ \deg(A)+ \gamma}{N} \right)^3 \log \left(a\binom{N+ \deg(A)+ \gamma}{N} \right) \left(h +T +1 \right),$$ which yields the desired result.
\end{proof}

Let $b \in D^m$ and $U$ be the maximum of the total degrees of the entries $a_{ij}$ and the coordinates of $b$, and  put $h:= h^{hom}(1, a_{ij}, b_i~|~ ~ (i,j) \in \{1, \cdots, m\} \times \{1, \cdots, n\})$.

\begin{pro}\label{pro}(Theorem \ref{thmB})
If the system $AX = b$ has a solution in $D^n$, then it admits a solution whose coordinates have total degrees at most  $$C := \gamma_1 \cdot  3 \left(a \binom{N+ U+ \gamma}{N} \right)^3 \log \left(a\binom{N+ U+ \gamma}{N} \right) \left(h +T +1 \right),$$ where $$\gamma_1 = d 3^{2^{3N}}(2 (a+1) \cdot U)^{2^{3N+1}(N+1)^N} $$ and 
$$\gamma = (2 (a+1) \cdot U)^{2((N+1)^N-1)}$$ and $$T=c_4R_K + h_K \log Q_{S}.$$ If $S= S_{\infty}$, then we replace $T$ by $c_4R_K +1$.
\end{pro}

\begin{proof}
By Proposition \ref{lemee} there exist generators $Y^{(1)}, \cdots, Y^{(M)}$ for the $D$-module of solutions to the homogeneous system of linear equations $(A, -b)Y =0$, where $Y =(Y_1 \cdots, Y_{m+1})$ is a vector of $n+1$ unknowns $Y_1 \cdots, Y_{m+1}$ such that
\begin{equation}\label{eq:h}
{\rm deg}(Y^{(j)}) \leq (2 (a+1) \cdot U)^{2((N+1)^N-1)} 
\end{equation}
\begin{equation}\label{eq:i}
h(Y^{(j)}) \leq  \frac{5}{2} \left(a \binom{N+ \deg(A)+ \gamma}{N} \right)^3 \log \left(a\binom{N+ U+ \gamma}{N} \right) \left(h +T +1 \right)
\end{equation}
for all $j \in  \{1, \cdots, M\}$. Clearly, if the system $AY=b$ is solvable, then $$1 \in (Y_{m+1}^{(1)}, \cdots, Y_{m+1}^{(M)}).$$ Conversely,suppose that there exist $h_1, \cdots, h_{m+1} \in D$ such that  $$1 = h_1Y_{m+1}^{(1)} + \cdots + h_{m+1}Y_{m+1}^{(M)}.$$ Then the vector $y \in D^{m}$  defined by $$ \binom{y}{1} = h_1Y^{(1)} + \cdots + h_{m+1}Y^{(M)}$$ is a solution to the system  $AY = b$. From Lemma \ref{lemI}, we deduce there exist $h_1, \cdots, h_{m+1}$  such that $$1 = h_1Y_{m+1}^{(1)} + \cdots + h_{m+1}Y_{m+1}^{(M)}$$ with total degree whose does not exceed
$${\rm deg} h_i \leq d(3\max_j \{\deg (Y_{m+1}^{(j)})\})^{2^{3N}}(\max_j \{h(Y_{m+1}^{(j)})\} +1 )$$ for all $i=1, \cdots, M$. Together with \eqref{eq:h} and \eqref{eq:i}, we get the desired result.
\end{proof}

Given $f_1, \cdots, f_m \in \mathcal{O}_S[x_1, \cdots,x_n]$, we write  $ I_K = <f_1, \cdots, f_m>$  for the ideal they generate in $K[x_1, \cdots,x_n]$.  A natural problem is to determine an explicit set of generators for the contraction  $$<f_1, \cdots, f_m>  \bigcap  \mathcal{O}_S[x_1, \cdots,x_n]$$ that is, the ideal obtained by contracting $I_K$ to  $\mathcal{O}_S[x_1, \cdots,x_n]$. 

For illustration, consider the case  $K= \mathbb{Q},~~\mathcal{O}_S = \mathbb{Z}$ and $f =x,~~ g= x+1$. We denote by $(f, g)$ the ideal of  $\mathbb{Z}[x]$ generated by $f,g$. we have $$ <f, g> \bigcap \mathbb{Z}[x] = \mathbb{Z}[x] = (f,g).$$ On the other hand, if we consider $f= x+y, ~~g =x-y$, then  
$$ <f, g> \bigcap \mathbb{Z}[x,y] = (x, y) \neq (f,g).$$

Evertse and Győry \cite[Theorem 6.1.7]{evertse-gyory-2022} solved this problem in the case  $K= \mathbb{Q}$ and $\mathcal{O}_S = \mathbb{Z}$.  We extend their result, providing an explicit description of the contracted ideal in this more general setting. 

\begin{pro}\label{pro2}(Theorem \ref{thmD})
Let $f_1, \cdots, f_m \in \mathcal{O}_S[x_1, \cdots, x_N]$ be the polynomials of total degrees  at most $D$ and $h = h^{hom}(f_1, \cdots , f_m)$ be the homogeneous logarithm height of the vector formed by their coefficients. Then $<f_1, \cdots, f_m>  \bigcap  \mathcal{O}_S[x_1, \cdots,x_N]$ is an ideal of 
$\mathcal{O}_S[x_1, \cdots,x_N]$ generated by polynomials of total degrees at most  $$ c :=  (2 m\cdot D)^{2((N+1)^N-1)}$$ and  homogeneous logarithm height  at most $$C_2 := (h+1) \left (\binom{c+D+N}{N} |d \Delta_K|\right)+ (s-1)c_5 h_K R_K {\displaystyle \prod_{i=1}^{t}} \log N(\mathcal{P}_i),$$  with
 $$c_5 :=  \dfrac{((s-1)!)^2}{2^{s-2}d^{s-1}}(53 d^2 \log 6d)^{s-2}$$
and $\Delta_K$  is the discriminant of $K$.  If $S= S_{\infty}$, then we  replace $h_K {\displaystyle \prod_{i=1}^{t}} \log N(\mathcal{P}_i)$ by $1$.
\end{pro}

\begin{proof}
Clearly $<f_1, \cdots, f_m> $ is the union of the ideals $\mathcal{I}_a ~(a \in \mathbb{N})$, where $\mathcal{I}_a$ consist of all polynomials $g \in \mathcal{O}_S[x_1, \cdots, x_N]$ such that $ag$ belongs to the ideal of $\mathcal{O}_S[x_1, \cdots, x_N]$ generated by $f_1, \cdots, f_m$. Consider the equation $$h_1f_1 + \cdots + h_mf_m - ag =0 \quad \mbox{in} ~~ (h_1, \cdots, h_m, g) \in (\mathcal{O}_S[x_1, \cdots,x_N])^{m+1}$$ By Proposition \ref{lemee}, the set of solutions of this equation is a module over $\mathcal{O}_S[x_1, \cdots,x_N]$  generated by tuples of polynomials of total degree which does not exceed $$c:= (2 m\cdot D)^{2((N+1)^N-1)}.$$ Consequently, for every $a \in \mathbb{N}$, the ideal $\mathcal{I}_a$ is generated by polynomials of the form  $h_1f_1 + \cdots + h_mf_m,$ where $$h_1, \cdots, h_m \in a^{-1}\mathcal{O}_S[x_1, \cdots,x_N]\quad \mbox{and} \quad {\rm deg}(h_i) \leq c$$ for every $i=1, \cdots, m$. It follows that the ideal  $<f_1, \cdots, f_m>  \bigcap  \mathcal{O}_S[x_1, \cdots,x_N]$ is generated by polynomials $$h_1f_1 + \cdots + h_mf_m,$$ where
\begin{equation}\label{eq:1"}
\left\{
\begin{aligned}
h_1, \cdots, h_m \in K[x_1, \cdots,x_N], \\
h_1f_1 + \cdots + h_mf_m \in \mathcal{O}_S[x_1, \cdots, x_N],\\    
{\rm deg}(h_i), \cdots, {\rm deg}(h_i) \leq c.
\end{aligned}
\right.
\end{equation}
Let $V$ be a $K$-vector space  spanned by all polynomials  $h_1f_1 + \cdots + h_mf_m$, where $h_1, \cdots, h_m$ satisfy \eqref{eq:1"} except for  the condition $$h_1f_1 + \cdots + h_mf_m \in \mathcal{O}_S[x_1, \cdots, x_N].$$  Since  $V$  is contained in the vector space of polynomials $K[x_1, \cdots,x_N]$ of degree $\leq c+D$, its dimension $L$ satisfies $$ L \leq \binom{c+D+N}{N}.$$ Consequently,  $V$ is spanned by the polynomials 
$$x_1^{\nu_1}\cdots x_N^{\nu_N}f_i ~(i= 1, \cdots, m,~~ \nu_1 + \cdots + \nu_N \leq c).$$ We may therefore select a basis $b_1, \cdots, b_L$ of $V$ from this set. In particular, each $b_i \in \mathcal{O}_S[x_1, \cdots, x_N]$ and has logarithm height  bounded by $h$ for all $i=1, \cdots, L$.  Let $\mathcal{A}$ be the $\mathcal{O}_S$ - module  generated by polynomials $h_1f_1 + \cdots + h_mf_m$ where $h_1, \cdots, h_m$ satisfy \eqref{eq:1"}. Then $\mathcal{A} \subseteq V$ and $b_1, \cdots, b_L \in \mathcal{A}$. By a classical result in module theory,  $\mathcal{A}$ is a finite generated torsion-free  $\mathcal{O}_S$-module, and hence a lattice over $\mathcal{O}_S$. Applying Lemma \ref{lemb}, we obtain a basis $\{P_1, \cdots, P_L\}$ of $\mathcal{A}$ such that 
\begin{equation*}
\left\{
\begin{aligned}
 P_1 &= \alpha_{11}b_1, \\
P_2 &= \alpha_{21}b_1 + \alpha_{22}b_2,\\
     \vdots \\     
P_L &= \alpha_{L1}b_1 + \alpha_{L2}b_2 + \cdots + \alpha_{LL}b_L,
\end{aligned}
\right.
\end{equation*}
with $$ \log |\alpha_{ij}|_v  \leq {\rm max} \left \{ {\displaystyle \sum_{i=1}^{s-1}} \max \{0, \log |\epsilon_i|_v\}, ~\log |d \Delta_K|_v \right \},$$ for~ all $i,j \in \{1, \cdots, n\}$ and $ v \in S$, where $c_5$ is defined in the preliminaries section \ref{sec1}. Clearly the total degrees of $P_j$ is at most $D+c$ for all $j=1, \cdots, L$. We denote by $a_{\textbf{i},j}, ~b_{\textbf{i},j} ~(j= 1, \cdots, L; ~ \textbf{i} =(i_1, \cdots, i_N)~/~ i_1 + \cdots + i_N \leq c+d)$ respectively coefficients of polynomials $P_j$ and $b_j$. We have 
\begin{equation}\label{eq:3'}
\begin{split}
\log |a_{\textbf{i},j}|_v
& \leq  \log \binom{c+D+N}{N}^{s(v)} + {\rm max}_k\{\log |\alpha_{jk}|_v\} + {\rm max}_k\{\log |b_{\textbf{i},k}|_v\}\\
& \leq \log \binom{c+D+N}{N}^{s(v)} + \log |d \Delta_K|_v +  {\displaystyle \sum_{i=1}^{s-1}} \max \{0, \log |\epsilon_i|_v\} \\
& +  {\displaystyle \sum_{k=1}^{L}} {\rm max}\{ \log |b_{\textbf{i},k}|_v\}\\
\end{split}
\end{equation}
 for all $v \in S$. Since $P_j \in \mathcal{O}_S[x_1, \cdots, x_N]$ for $j= 1, \cdots L$, we deduce that $$\log |a_{\textbf{i},j}|_v \leq 0 \quad \mbox{for} ~~ \mbox{all} \quad v \notin S.$$  Together with \eqref{eq:3'} and Lemma \ref{lemC}, we obtain 
$$h^{hom}(P_j) \leq (h+1) \left (\binom{c+D+N}{N} |d \Delta_K|\right)+ (s-1)c_5 h_K R_K {\displaystyle \prod_{i=1}^{t}} \log N(\mathcal{P}_i) $$  if $S \neq S_{\infty}$ and 
$$h^{hom}(P_j) \leq (h+1)\left(\binom{c+D+N}{N} |d \Delta_K|\right )+ (s-1)c_5 R_K $$ otherwise. This holds for all $j= 1, \cdots ,  L$ which give us the desired result. 
\end{proof}

\section{Effective results on G.C.D for polynomials }\label{sec3}

Let $K$ be a number field of degree $d$. For fixed $u = (u_1, \cdots, u_n) \in K^n$ and $i=(i_1,\cdots, i_n) \in \mathbb{N}^n$; we denote $|i| = i_1 + \cdots + i_n$ and $u^i:= u_1^{i_1}\cdots u_n^{i_n}$. Let $m$ be a positive integer, we define $$K[x_1, \cdots, x_n]_m = \{f \in K[x_1, \cdots, x_n] ~/~ {\rm deg}(f) \leq m \}$$ and for $f, g \in K[x_1, \cdots, x_n]$; $$(f,g)_m = \{fp+gq ~/ ~ {\rm deg}(fp) \leq m, ~{\rm deg}(gq) \leq m\}$$
For the proof of our Theorem \ref{thmJ}, we will need the following result which is the inhomogeneous version of  \cite[Lemma 3.2]{X}.

\begin{lem}\label{lem1}
Let $f_1, f_2 \in K[x_1, \cdots, x_N]$ be coprime polynomial of total degree $d_1$ and $d_2$ respectively. Let $m> {\rm max}\{d_1, d_2\}$ be an integer and  $B$ be a basis of the $K$-vector space $K[x_1, \cdots, x_N]_m/(f_1,f_2)_m$. Then $$|B| \leq d_1 d_2 \binom{m+n-2}{n-2},$$ where $|B|$ denote the cardinality of the set $B$. 
\end{lem}

We prove Theorem \ref{thmJ} by decomposing the sum on its left-hand side into two contributions, corresponding to places $v \notin S$ and $v \in S$, respectively. A noneffective version of this result was established in \cite{X}. As explained in the introduction, our aim here is to obtain an effective version, in the sense of providing an explicit bound for the number and the total degree of exceptional hypersurfaces. We begin by considering the contribution from places $v \notin S$.

\begin{pro}\label{pro3}
Let $K$ be a number field of degree $d$ and $S$ be a finite set of places of $K$ containing all archimedean ones. We denote the cardinality of $S$ by $s$. Let $f, g \in K[x_1, \cdots, x_n]$ be coprime polynomials  of total degree at most $D$  and $ 0< \epsilon < 1 $ be a real number. For $$\eta < \dfrac{\epsilon^2}{72D^2n^4} $$ we consider the inequality 
\begin{equation}\label{eq:A3}
- {\displaystyle \sum_{v \in M_K  \setminus S}^{}} \log ^- {\rm max} \left \{ |f(u_1, \cdots, u_n)|_v ,~ |g(u_1, \cdots, u_n)|_v \right \} > \epsilon {\displaystyle \sum_{i=1}^{n}} h(u_i)
\end{equation}
with $(u_1, \cdots, u_n) \in \mathbb{G}_{m}^{n}(K)_{S, \eta}$. There are hypersurfaces $H_1, \cdots , H_{t_1(\epsilon)}$ of $K^n$ of total degree at most $m: = 12D^2n^3\epsilon^{-1}$ with $$t_1(\epsilon )  := {\rm exp}\left( 75snd(2m)^{2n}(1+\log \epsilon^{-1})\right)\left( \log (6m)+ 2h^{hom}(f,g, 1) \right)^{3nd}$$ such that any solutions $(u_1, \cdots, u_n) \in \mathbb{G}_{m}^{n}(K)_{S, \eta}$ of \eqref {eq:A3} belongs to $$H_1 \bigcup  H_2 \bigcup \cdots \bigcup H_{t_1(\epsilon)}.$$
\end{pro}

\begin{proof}
Our approach extends the method developed in \cite{X}, which was originally introduced by Levin \cite{L}. Let $(f,g)$ be an ideal of $K[x_1, \cdots, x_n]$ generated by $f$ and $g$. Assume that $$(f,g)_m \neq K[x_1, \cdots, x_n]_m \quad \mbox{with} \quad  m= \lfloor{12 D^2n^3\epsilon^{-1}}\rfloor.$$ Then the $K$-vector space $$V_m := K[x_1, \cdots, x_N]_m/(f, g)_m$$ is nonzero. Fix $u:= (u_1, \cdots, u_n) \in  \mathbb{G}_{m}^{n}(K)_{S, \eta}$. For each $v \in S$, we construct a basis $B_v$ of the vector space $V_m$ as follows. First, choose a monomial $$X^{\textbf{i}_1} := x_1^{i_1} \cdots  x_1^{i_n} \in K[x_1, \cdots, x_n]$$ such that $|u^{\textbf{i}_1}|_v$  is minimal among all monomials  $X^{\textbf{i}} \notin (f, g)_m$. Suppose that monomials $X^{\textbf{i}_1}, \cdots, X^{\textbf{i}_j}$  have been chosen and are linearly independent modulo $(f,g)_m$ but do not yet span $V_m$.We then choose a monomial $X^{\textbf{i}_{j+1}} \in K[x_1, \cdots, x_n]$ for which $|u^{\textbf{i}_{j+1}}|_v$ is minimal, subject to the condition that $X^{\textbf{i}_1}, \cdots, X^{\textbf{i}_{j+1}}$ are linearly independent modulo $(f,g)_m$. Continuing in this way, we obtain a basis$ B_v$ of $V_m$ consisting of the residue classes of monomials  $$X^{\textbf{i}_1}, \cdots, X^{\textbf{i}_{N'}} \quad \mbox{with} \quad N' = {\rm dim} V_m.$$ Since there are at most $N = \binom{m+n}{n}$ monomials of degree at most $m$, it follows  that the number of possible $S$-tuples of bases $(B_v)_v \in S$ is bounded above by $$A:= \left( \dfrac{N!}{(N-N')!}\right)^s \leq {\rm exp}\left( sn(2m)^n \log (2m)\right) \leq {\rm exp}\left( sn^{n+1}(2D^2)^n \log (2D^2n)\right).$$ Although this bound could be improved by means of Stirling's formula, the above estimate will suffice for our purposes.
We now arrange the multi-indices $\textbf{i}_1, \cdots, \textbf{i}_{N'}$ in lexicographic order and put $M= {\rm dim}(f, g)_m$. Since $(f,g)_m$ is an subspace of $K[x_1, \cdots, x_n]_m$ generated by the polynomials $$x_1^{i_1} \cdots  x_1^{i_n}f,  ~x_1^{i_1} \cdots  x_1^{i_n}g ~~(i_1 + \cdots + i_n \leq m-D),$$  we may choose a basis $b_1, \cdots, b_{M}$ of $(f, g)_m$ from this set. Fix a multi-index $|\textbf{i}_0| $ with $|\textbf{i}_0| \leq m$ and  let $I_v = \{\textbf{i}_1, \cdots, \textbf{i}_{N'}\}.$ By the definition of the basis $B_v$, there exist coefficients $u_{\textbf{i}_0, \textbf{i}_j}$ such that $$X^{\textbf{i}_0} + {\displaystyle \sum_{j=1}^{N'}} u_{\textbf{i}_0, \textbf{i}_j}X^{\textbf{i}_j} = ff_1 + gg_1,$$ where $$f_1 = {\displaystyle \sum_{|\textbf{l}| \leq m-D}^{}} y_{\textbf{l}}x_1^{l_1}\cdots x_n^{l_n} \quad \mbox{and} \quad g_1 = {\displaystyle \sum_{|\textbf{t}| \leq m-D}^{}} z_{\textbf{t}}x_1^{t_1}\cdots x_n^{t_n};$$ and $$f = {\displaystyle \sum_{|\textbf{j}| \leq D}^{}} a_{\textbf{j}}x_1^{j_1}\cdots x_n^{j_n} \quad \mbox{and} \quad g = {\displaystyle \sum_{|\textbf{j}| \leq D}^{}} b_{\textbf{j}}x_1^{j_1}\cdots x_n^{j_n}.$$ Our next goal is to obtain explicit upper bounds for the heights of the coefficients  $ u_{\textbf{i}_0, \textbf{i}_j},  ~y_{\textbf{l}}, ~z_{\textbf{t}}.$ To this end, fix a set $I= \{\textbf{i}_1, \cdots, \textbf{i}_{M}\}$. We are led to the following system of inhomogeneous linear equations:\\
 If $\textbf{i}_0 \notin I_v$,  then
 \begin{equation}\label{eq:A4}
\left\{
\begin{aligned}
{\displaystyle \sum_{\textbf{l} + \textbf{j} = \zeta}^{}} a_{\textbf{j}}y_{\textbf{l}} + {\displaystyle \sum_{\textbf{l} + \textbf{j} = \zeta}^{}} b_{\textbf{j}}z_{\textbf{l}} - u_{\textbf{i}_0, \zeta} = 0 ~~& \mbox{ if } ~~\zeta \in I_v, \\
	{\displaystyle \sum_{\textbf{l} + \textbf{j} = \zeta}^{}} a_{\textbf{j}}y_{\textbf{l}} + {\displaystyle \sum_{\textbf{l} + \textbf{j} = \zeta}^{}} b_{\textbf{j}}z_{\textbf{l}}  = 0 ~~& \mbox{ if } ~~\zeta \notin I_v \cup \{\textbf{i}_0\}, \\
{\displaystyle \sum_{\textbf{l} + \textbf{j} = \zeta}^{}} a_{\textbf{j}}y_{\textbf{l}} + {\displaystyle \sum_{\textbf{l} + \textbf{j} = \zeta}^{}} b_{\textbf{j}}z_{\textbf{l}}  = 1 ~~& \mbox{ if } ~~\zeta = \textbf{i}_0.
\end{aligned}
\right.
\end{equation}
Otherwise (i.e. $\textbf{i}_0 \in I_v$),
 \begin{equation}\label{eq:A5}
\left\{
\begin{aligned}
{\displaystyle \sum_{\textbf{l} + \textbf{j} = \zeta}^{}} a_{\textbf{j}}y_{\textbf{l}} + {\displaystyle \sum_{\textbf{l} + \textbf{j} = \zeta}^{}} b_{\textbf{j}}z_{\textbf{l}} - u_{\textbf{i}_0, \zeta} = 0 ~~& \mbox{ if } ~~\zeta \in I_v \setminus \{\textbf{i}_0\} , \\
	{\displaystyle \sum_{\textbf{l} + \textbf{j} = \zeta}^{}} a_{\textbf{j}}y_{\textbf{l}} + {\displaystyle \sum_{\textbf{l} + \textbf{j} = \zeta}^{}} b_{\textbf{j}}z_{\textbf{l}}  = 0 ~~& \mbox{ if } ~~\zeta \notin I_v, \\
{\displaystyle \sum_{\textbf{l} + \textbf{j} = \zeta}^{}} a_{\textbf{j}}y_{\textbf{l}} + {\displaystyle \sum_{\textbf{l} + \textbf{j} = \zeta}^{}} b_{\textbf{j}}z_{\textbf{l}}   - u_{\textbf{i}_0, \zeta} = 1 ~~& \mbox{ if } ~~\zeta = \textbf{i}_0.
\end{aligned}
\right.
\end{equation}
Thus, we obtain an inhomogeneous linear system $AY = b$, where $A$ is an  $r\times T$ matrix,  with $$r:= \binom{n+m}{n} \quad \mbox{and} \quad T:= 2\binom{n+m-D}{n}+ N'$$ and whose entries belong to the set $\{0, -1, a_{\textbf{j}}, b_{\textbf{j}}; ~~ |\textbf{j}| \leq D\}$. Since the system $AY = b$ is solution, we can explicitly construct a solution as follows. Let $\Delta$ be a $k\times k$-sub-matrix of $A$ with $\delta = {\rm det}(\Delta) \neq 0,$ where $k= {\rm rank}(A) \geq 1$. Multiplying the corresponding $k$ entries of $b$ by the adjugate matrix  $\Delta^{ad}$ of $\Delta$, we obtain a vector $d: = (d_1, \cdots, d_k)$. Then, up to a permutation of the coordinates corresponding to the position of our sub-matrix $\Delta$ in $A$, a solution of  equation \eqref{eq:A4} or \eqref{eq:A5} is given by 
 \begin{equation}\label{eq:O}
\left\{
\begin{aligned}
d_j/ \delta ~~ & \mbox{ if } ~~  \mbox{ for}~~ 1 \leq j \leq k, \\
	0 ~~ & \mbox{ if } ~~  \mbox{ for}~~ k < j \leq T.
\end{aligned}
\right.
\end{equation}
In particular, this provides explicit expressions for the coefficients  $ u_{\textbf{i}_0, \textbf{i}_j},  ~y_{\textbf{l}}, ~z_{\textbf{t}}$. Since $\{b_1, \cdots, b_M \}$ is a base of $(f, g)_m$,  for each multi-index  $\textbf{i}$, there exists a linear form $L_{\textbf{i}}^{v}$ over $K$ such that $$L_{\textbf{i}}^{v}(b_1, \cdots, b_M) = X^{\textbf{i}} + {\displaystyle \sum_{j=1}^{N'}} u_{\textbf{i}, \textbf{i}_j}X^{\textbf{i}_j}.$$ Moreover, the collection $\{L_{\textbf{i}}^{v}(b_1, \cdots, b_M) ~:~~|\textbf{i}| \leq m,~ \textbf{i} \notin I_v\}$ forms a basis for $(f,g)_m$. Let $$P = b(\textbf{u}):= (b_1(\textbf{u}), \cdots, b_M(\textbf{u})) \in K^M.$$ By enlarging the exceptional set of hypersurfaces, we may assume that $b(\textbf{u}) \neq 0$. Finally, by the minimality property in the construction of the monomials $X^{\textbf{i}_1}, \cdots, X^{\textbf{i}_{N'}}$, together with the triangle inequality, we obtain, for every multi-index $\textbf{i}$ satisfying $| \textbf{i}| \leq m, ~~ \textbf{i} \notin I_v$ the following estimate:
\begin{equation}\label{eq:p}
{\rm log} |L_{\textbf{i}}^{v}(P)|_v \leq {\rm log} | \textbf{u}^{\textbf{i}}|_v + c_v,
\end{equation}
where $$c_v \leq s(v) {\rm log}(N'+1) + {\displaystyle \sum_{j=1}^{N'}} {\displaystyle \sum_{|\textbf{i}| \leq m}^{}} \log \max \{1, | u_{\textbf{i}, \textbf{i}_j}|_v\}.$$ By convention we set $h(0)=0$. Notice that
\begin{equation}\label{eq:X}
\begin{split}
h(u_{\textbf{i}, \textbf{i}_j})
& ~\leq ~ h \left(\frac{d_j}{\delta} \right) ~ \leq~~ h(d_j) + h(\delta) \\
& ~\leq~ T \log T + 2T h(f, g, 1)~ \leq ~~3(2m)^n \left( \log (3(2m)^n )+ 2h^{hom}(f,g, 1) \right)\\
& ~\leq~ 3n(2m)^n \left( \log (6m)+ 2h^{hom}(f,g, 1 ) \right), 
\end{split}
\end{equation}
where we used Hadamard inequality,  Lemma \ref{lemA} and the fact that $m \geq n$ and  $\binom{n+m}{n} \leq (2m)^n$. Consequently
\begin{equation}\label{eq:q'}
\begin{split}
c := {\displaystyle \sum_{v \in S}^{}}c_v & ~\leq ~ \log(N+1) + (2m)^{2n} h(u_{\textbf{i}, \textbf{i}_j})\\
& ~\leq ~ n\log (6m) + 3n(2m)^{3n} \left( \log (6m)+ 2h^{hom}(f,g, 1) \right)  \\
& ~\leq~ 4n(2m)^{3n} \left( \log (6m)+ 2h^{hom}(f,g, 1) \right).
\end{split}
\end{equation}
On the other hand, for every multi-index  $\textbf{i}$ satisfying  $| \textbf{i}| \leq m, ~~ \textbf{i} \notin I_v$, it follows from \eqref{eq:p}  that $L_{\textbf{i}}^{v}$ is a linear form whose  coefficients $a_{\textbf{i}, \textbf{j}}$ are of the form given in \eqref{eq:O}. Consequently, by \eqref{eq:X}, we obtain  $$ h(a_{\textbf{i}, \textbf{j}}) \leq 3n(2m)^n \left( \log (6m)+ 2h^{hom}(f,g, 1 ) \right).$$ Applying Hadamard's inequality, we get $$\log |{\rm det}(L_{\textbf{i}}^{v},~| \textbf{i}| \leq m, ~~ \textbf{i} \notin I_v)|_v \leq  s(v)\frac{N}{2} \log N + {\displaystyle \sum_{j=1}^{N'}} {\displaystyle \sum_{|\textbf{i}| \leq m}^{}}\log  \max \{1, |a_{\textbf{i}_j, \textbf{i}}|_v\}.$$ Thus
\begin{equation}\label{eq:t}
\begin{split}
{\displaystyle \sum_{v \in S}^{}} \log |{\rm det}(L_{\textbf{i}}^{v},~| \textbf{i}| \leq m, ~~ \textbf{i} \notin I_v)|_v  & ~\leq ~\frac{N}{2} \log N  + (2m)^{2n} h(a_{\textbf{i}, \textbf{j}})\\
& ~\leq ~ n(2m)^{n} \log (6m) + 3n(2m)^{3n} \left( \log (6m)+ 2h^{hom}(f,g, 1) \right)  \\
& ~\leq~ 4n(2m)^{3n} \left( \log (6m)+ 2h^{hom}(f,g, 1) \right):= \log H.
\end{split}
\end{equation}
It is clear that $H(L_{\textbf{i}}^{v}) \leq H$ and 
\begin{equation}\label{eq:Z10}
h \left({\rm det}(L_{\textbf{i}}^{v},~| \textbf{i}| \leq m, ~~ \textbf{i} \notin I_v)\right) \leq \log H.
\end{equation}
Therefore we may now apply a quantitative $p$-adic subspace theorem  (see Lemma \ref{thm1}). To this end, we first estimate $${\displaystyle \sum_{v \in S}^{}} {\displaystyle \sum_{|\textbf{i}|\leq m,~\textbf{i} \notin I_v}^{}}\log \dfrac{|P|_v}{|L_{\textbf{i}}^{v}(P)|_v}.$$ Using \eqref{eq:p}, we get $$ - {\displaystyle \sum_{v \in S}^{}} {\displaystyle \sum_{|\textbf{i}|\leq m,~\textbf{i} \notin I_v}^{}}\log |L_{\textbf{i}}^{v}(P)|_v  \geq - {\displaystyle \sum_{v \in S}^{}} {\displaystyle \sum_{|\textbf{i}|\leq m,~\textbf{i} \notin I_v}^{}}\log | \textbf{u}^{\textbf{i}}|_v - Mc,$$ where $c$ is given by the relation \eqref{eq:q'}. By the product formula, we have  $${\displaystyle \sum_{v \in S}^{}} \log | \textbf{u}^{\textbf{i}}|_v = -{\displaystyle \sum_{v \in M_K\setminus S}^{}} \log | \textbf{u}^{\textbf{i}}|_v.$$ Therefore
\begin{equation}\label{eq:m}
\begin{split}
- {\displaystyle \sum_{v \in S}^{}} {\displaystyle \sum_{|\textbf{i}|\leq m,~\textbf{i} \notin I_v}^{}}\log | \textbf{u}^{\textbf{i}}|_v & ~= ~{\displaystyle \sum_{v \in S}^{}} {\displaystyle \sum_{|\textbf{i}|\leq m,~\textbf{i} \in I_v}^{}}\log | \textbf{u}^{\textbf{i}}|_v - {\displaystyle \sum_{v \in S}^{}} {\displaystyle \sum_{|\textbf{i}|\leq m}^{}}\log | \textbf{u}^{\textbf{i}}|_v\\
& ~= ~ {\displaystyle \sum_{v \in S}^{}} {\displaystyle \sum_{|\textbf{i}|\leq m,~\textbf{i} \in I_v}^{}}\log | \textbf{u}^{\textbf{i}}|_v  +  {\displaystyle \sum_{v \in M_K \setminus S}^{}} {\displaystyle \sum_{|\textbf{i}|\leq m}^{}}\log | \textbf{u}^{\textbf{i}}|_v.
\end{split}
\end{equation}
From Lemma \ref{lem1}, we deduce that 
\begin{equation}\label{eq:n}
-{\displaystyle \sum_{v \in S}^{}} {\displaystyle \sum_{|\textbf{i}|\leq m, \textbf{i} \in I_v}^{}}\log | \textbf{u}^{\textbf{i}}|_v \leq mD^2\binom{m+n-2}{n-2}{\displaystyle \sum_{i=1}^{n}} h(u_i).
\end{equation}
Moreover $$ {\displaystyle \sum_{v \in M_K \setminus S}^{}} {\displaystyle \sum_{|\textbf{i}|\leq m}^{}}\log | \textbf{u}^{\textbf{i}}|_v =  {\displaystyle \sum_{v \in M_K \setminus S}^{}} \left ( {\displaystyle \sum_{j=0}^{m-1}}(m-j)\binom{n-1+j}{n-1}  \right) {\displaystyle \sum_{i=1}^{n}} \log |u_i|_v.$$
From Hockey-Stick identity, we deduce that 
\begin{equation}\label{eq:Z}
{\displaystyle \sum_{j=0}^{m-1}}(m-j)\binom{n-1+j}{n-1} = \binom{n+m}{n+1}.
\end{equation}
Therefore
$$
 {\displaystyle \sum_{v \in M_K \setminus S}^{}}{\displaystyle \sum_{|\textbf{i}|\leq m}^{}}\log | \textbf{u}^{\textbf{i}}|_v \geq - \binom{n+m}{n+1} {\displaystyle \sum_{v \in M_K \setminus S}^{}} {\displaystyle \sum_{i=1}^{n}} \log \max \left \{ 1, \left| \frac{1}{u} \right|_v\right\}.
$$
Together with \eqref{eq:m} and \eqref{eq:n} we obtain 

\begin{equation}\label{eq:t}
\begin{split}
- {\displaystyle \sum_{v \in S}^{}}{\displaystyle \sum_{|\textbf{i}|\leq m, \textbf{i} \notin I_v}^{}}\log |  L_{\textbf{i}}^{v}(P)|_v   & ~\geq ~ -mD^2\binom{m+n-2}{n-2}{\displaystyle \sum_{i=1}^{n}} h(u_i)\\
& ~- ~ \binom{n+m}{n+1} {\displaystyle \sum_{v \in M_K \setminus S}^{}} {\displaystyle \sum_{i=1}^{n}} \log \max \left \{ 1, \left| \frac{1}{u_i} \right|_v\right\}  - Mc.
\end{split}
\end{equation}
On the other hand, 
\begin{equation}\label{eq:u}
\begin{split}
{\displaystyle \sum_{v \in S}^{}}{\displaystyle \sum_{|\textbf{i}|\leq m, \textbf{i} \notin I_v}^{}}\log |P|_v   & ~= ~ M\left (h^{hom}(P)- {\displaystyle \sum_{v \in M_K \setminus S}^{}} \log |P|_v \right ).
\end{split}
\end{equation}
Now we want to estimate $$ {\displaystyle \sum_{v \in M_K \setminus S}^{}} \log |P|_v  = {\displaystyle \sum_{v \in M_K \setminus S}^{}} \max_i \{\log |b_i(\textbf{u})|_v\}.$$ For $v \in  M_K \setminus S$ and $i \in \{1, \cdots, M \}$, if we assume that $$b_i(\textbf{u}) = u_1^{j_1}\cdots u_n^{j_n}f(\textbf{u}) \quad \mbox{with} \quad  j_1+ \cdots + j_n \leq m-D,$$ then 
\begin{equation*}
\begin{split}
 | b_i(\textbf{u})|_v  & ~\leq ~ \left( \max_j |u_j|_v \right)^{j_1+ \cdots + j_n} |f(\textbf{u})|_v \\
& ~ \leq~    (\max\{1, |\textbf{u}|_v\})^{m} \cdot \max _{\textbf{i}}\{|a_{\textbf{i}}|_v, |b_{\textbf{i}}|_v, 1\} \cdot \min \left \{ 1, \max\{|f(\textbf{u})|_v, |g(\textbf{u})|_v\} \right \}.
\end{split}
\end{equation*}
Otherwise i.e.  $$b_i(\textbf{u}) = u_1^{j_1}\cdots u_n^{j_n}g(\textbf{u}) \quad \mbox{with} \quad  j_1+ \cdots + j_n \leq m-D,$$ we get the same bound. We sum over  $v \in  M_K \setminus S$ and use the notation $$ \log^{-}\max\{|f(\textbf{u})|_v, |g(\textbf{u})|_v\} := \min \left \{ 0, \log \max\{|f(\textbf{u})|_v, |g(\textbf{u})|_v\} \right \}$$ to get 
\begin{equation}\label{eq:v}
\begin{split}
{\displaystyle \sum_{v \in M_K \setminus S}^{}} \log |P|_v  & \leq {\displaystyle \sum_{v \in M_K \setminus S}^{}}\log^{-}\max\{|f(\textbf{u})|_v, |g(\textbf{u})|_v\} \\
&~ + ~m {\displaystyle \sum_{v \in M_K \setminus S}^{}} \log \max \{1, |\textbf{u}|_v\} ~+~ h^{hom}(f,g,1).
\end{split}
\end{equation}
For $v \in S$, we have $$ | b_i(\textbf{u})|_v  \leq  \left(\binom{n+D}{n}\right)^{s(v)}  |\textbf{u}|_v^{m} \cdot \max _{\textbf{i}}\{1, |a_{\textbf{i}}|_v, |b_{\textbf{i}}|_v\} $$ and so
\begin{equation}\label{eq:x}
h^{hom}(P) \leq mh^{hom}(\textbf{u}) + \binom{n+D}{n} + h^{hom}(f,g, 1).
\end{equation}
From \eqref{eq:u} and  \eqref{eq:v}, we deduce 
\begin{equation}\label{eq:w'}
\begin{split}
{\displaystyle \sum_{v \in S}^{}} ~ {\displaystyle \sum_{|\textbf{i}|\leq m, \textbf{i} \notin I_v}^{}}\log |P|_v   & ~\geq~ M h^{hom}(P)- M{\displaystyle \sum_{v \in M_K \setminus S}^{}}\log^{-}\max\{|f(\textbf{u})|_v, |g(\textbf{u})|_v\} \\
& ~ - ~ mM {\displaystyle \sum_{v \in M_K \setminus S}^{}} ~{\displaystyle \sum_{i=1}^{n}}\log \max \{1, |u_i|_v\} ~-~ M h^{hom}(f,g, 1) .
\end{split}
\end{equation}
For fixed $(B_v)_{v \in S}$, we apply Lemma \ref{thm1} with $\delta = \epsilon /4$ and deduce that there are finitely many hypersurfaces of degree at most $m$ denoted  $H_1, \ldots, H_{l}$ of $K^n,$ with $$l < \left( 2^{67M^2} \left(  \dfrac{\epsilon}{4} \right) ^{7M}\right)^{s+1}(3 + \log\log H)$$ such that every  $(S, \eta)$-unit which is not in $H_1 \bigcup \cdots \bigcup H_l$ satisfy
\begin{equation*}\label{eq:w}
\begin{split}
{\displaystyle \sum_{v \in S}^{}} {\displaystyle \sum_{|\textbf{i}|\leq m,~\textbf{i} \notin I_v}^{}}\log \dfrac{|P|_v}{|L_{\textbf{i}}^{v}(P)|_v} &~\leq~ \left(M+ \dfrac{\epsilon}{4} \right)\left(h^{hom}(P)+\dfrac{\log M}{2} \right)\\
&~ - {\displaystyle \sum_{v \in S}^{}} \log |{\rm det}(L_{\textbf{i}}^{v},~| \textbf{i}| \leq m, ~~ \textbf{i} \notin I_v)|_v, 
\end{split}
\end{equation*}
where $P = b(\textbf{u})$. Together with  Lemma \ref{lemA} and \eqref{eq:Z10}, we get
\begin{equation}\label{eq:A6}
\begin{split}
{\displaystyle \sum_{v \in S}^{}}~ {\displaystyle \sum_{|\textbf{i}|\leq m,~\textbf{i} \notin I_v}^{}}\log \dfrac{|P|_v}{|L_{\textbf{i}}^{v}(P)|_v} &~\leq~ \left(M+  \dfrac{\epsilon}{4} \right)\left(h^{hom}(P)+\dfrac{\log M}{2} \right)\\
&~ + 4sn(2m)^{3n} \left( \log (6m)+ 2h^{hom}(f,g,1) \right).
\end{split}
\end{equation}
Combining \eqref{eq:t}, \eqref{eq:w'} and \eqref{eq:A6}, we get  
\begin{equation*}
\begin{split}
-M{\displaystyle \sum_{v \in M_K \setminus S}^{}}\log^{-}\max\{|f(\textbf{u})|_v, |g(\textbf{u})|_v\} ~\leq~ 
mD^2\binom{m+n-2}{n-2}{\displaystyle \sum_{i=1}^{n}} h(u_i) \\ ~+ ~ \binom{n+m}{n+1} {\displaystyle \sum_{v \in M_K \setminus S}^{}} ~{\displaystyle \sum_{i=1}^{n}} \log \max \left \{ 1, \left| \frac{1}{u_i} \right|_v\right\}\\
~ + ~ mM {\displaystyle \sum_{v \in M_K \setminus S}^{}}~ {\displaystyle \sum_{i=1}^{n}}\log \max \{1, |u_i|_v\}  ~+~  \dfrac{m\epsilon}{4}h^{hom}(\textbf{u}) \\
~ + ~11ns(2m)^{4n} \left( \log (6m)+ 2h^{hom}(f,g,1) \right),
\end{split}
\end{equation*}
where we have used relation \eqref{eq:x} and a bound of $c$ given in \eqref{eq:q'}. Notice that  $$\binom{n+m}{n+1} = \frac{m}{n+1}\binom{m+n}{n}.$$ From Lemma \ref{lem1}, we have $$M \geq  \binom{m+n}{n}-D^2\binom{m+n-2}{n-2} \geq \frac{1}{n} \binom{m+n}{n}.$$ It follows that 
\begin{equation*}
\begin{split}
-{\displaystyle \sum_{v \in M_K \setminus S}^{}}\log^{-}\max\{|f(\textbf{u})|_v, |g(\textbf{u})|_v\} ~\leq~ 
 \left( \dfrac{D^2n^3}{m+1} + \dfrac{\epsilon}{4} \right) {\displaystyle \sum_{i=1}^{n}} h(u_i) \\ ~+ ~ m \left(  {\displaystyle \sum_{i=1}^{n}}  {\displaystyle \sum_{v \in M_K \setminus S}^{}}\log \max \left \{ 1, \left| \frac{1}{u_i} \right|_v\right\} 
~ + ~  \log \max \{1, |u_i|_v\} \right)   \\
~ + ~11ns(2m)^{4n} \left( \log (6m)+ 2h^{hom}(f,g, 1) \right),
\end{split}
\end{equation*}
where we used the fact that $$ \binom{m+n-2}{n-2} \leq \dfrac{n^2}{(m+1)m} \binom{m+n}{n}  $$
Since $\textbf{u} \in \mathbb{G}_m^n(K)_{S, \eta},$ we deduce $${\displaystyle \sum_{v \in M_K \setminus S}^{}}\log \max \left \{ 1, \left| \frac{1}{u_i} \right|_v\right\} 
~ + ~  \log \max \{1, |u_i|_v\}  \leq \eta h(\textbf{u}) \leq  \eta {\displaystyle \sum_{i=1}^{n}} h(u_i).$$ Therefore 
\begin{equation*}
\begin{split}
-{\displaystyle \sum_{v \in M_K \setminus S}^{}}\log^{-}\max\{|f(\textbf{u})|_v, |g(\textbf{u})|_v\} ~\leq~ 
\left( \dfrac{D^2n^3}{m+1} + \dfrac{\epsilon}{4} \right){\displaystyle \sum_{i=1}^{n}} h(u_i) \\
~+ ~ mn \eta {\displaystyle \sum_{i=1}^{n}} h(u_i)  
~ + 11ns(2m)^{4n} \left( \log (6m)+ 2h^{hom}(f,g, 1) \right)\\
~ \leq~ \left( \dfrac{D^2n^3}{m+1} + \dfrac{\epsilon}{4} +  mn \eta\right) {\displaystyle \sum_{i=1}^{n}} h(u_i) + 11ns(2m)^{4n} \left( \log (6m)+ 2h^{hom}(f,g, 1) \right).
\end{split}
\end{equation*}
Using the fact that $m= 12 D^2n^3 \epsilon^{-1}$, we deduce that 
\begin{equation}\label{eq:A}
\left( \dfrac{D^2n^3}{m+1} + \dfrac{\epsilon}{4} +  mn \eta\right) \leq \frac{\epsilon}{2}
\end{equation}
 whenever $$\eta < \dfrac{\epsilon^2}{72D^2n^4}.$$ So
\begin{equation*}
\begin{split}
-{\displaystyle \sum_{v \in M_K \setminus S}^{}}\log^{-}\max\{|f(\textbf{u})|_v, |g(\textbf{u})|_v\} 
~ \leq~  \dfrac{\epsilon}{2}  {\displaystyle \sum_{i=1}^{n}} h(u_i) \\
+ 11ns(2m)^{4n} \left( \log (6m)+ 2h^{hom}(f,g, 1) \right).
\end{split}
\end{equation*}
From Lemma \ref{lemE}, there are at most 
\begin{equation}\label{eq:B1}
t_2:=  \left(44^{4nd}(ns)^{3nd}\epsilon^{-3nd}(2m)^{12n^2d} \right)\left( \log (6m)+ 2h^{hom}(f,g, 1) \right)^{3nd} 
\end{equation}
vectors of $K^n$ such that $${\displaystyle \sum_{i=1}^{n}} h(u_i) \leq 22\epsilon^{-1}ns(2m)^{4n} \left( \log (6m)+ 2h^{hom}(f,g, 1) \right).$$
Excepted those vectors, we obtain 
\begin{equation*}
\begin{split}
-{\displaystyle \sum_{v \in M_K \setminus S}^{}}\log^{-}\max\{|f(\textbf{u})|_v, |g(\textbf{u})|_v\} 
~ \leq~  \epsilon {\displaystyle \sum_{i=1}^{n}} h(u_i).
\end{split}
\end{equation*}
In conclusion, the choice of linear forms in the application of the quantitative Subspace Theorem above depends only on the choice of the monomial bases $(B_v)_{v \in S}$. Since the number of possible such choices is at most $${\rm exp}\left( sn(2m)^n \log (2m)\right),$$ we obtain the following conclusion. Suppose that  $(f,g)_m \neq K[x_1, \cdots, x_n]_m$, then there exist hypersurfaces $H_1, \cdots , H_{t_1}$ of $K^n$ each of total degree at most  $m :=  12D^2n^3\epsilon^{-1} $ such that
\begin{equation*}
\begin{split}
 {\rm exp}\left( sn(2m)^n \log (2m)\right)\left( 2^{74M^2} \epsilon ^{-7M}\right)^{s+1}(3 + \log\log H)  + \\ \left(44^{4nd}(ns)^{3nd}\epsilon^{-3nd}(2m)^{12n^2d} \right)\left( \log (6m)+ 2h^{hom}(f,g, 1) \right)^{3nd}\\
~\leq ~ {\rm exp}\left( 75snd(2m)^{2n}(1+\log \epsilon^{-1})\right)\left( \log (6m)+ 2h^{hom}(f,g,1) \right)^{3nd}=: t_1(\epsilon)
\end{split}
\end{equation*}
and every solutions $(u_1, \cdots, u_n) \in \mathbb{G}_{m}^{n}(K)_{S, \eta}$ of \eqref {eq:A3} belongs to $$H_1 \bigcup  H_2 \bigcup \cdots \bigcup H_{t_1}.$$
We now turn to the case $(f,g)_m = K[x_1, \cdots, x_n]_m$. In this situation, there exist polynomials $f,g \in  K[x_1, \cdots, x_n]$  such that   $$1 = ff_1 + gg_1,$$ with $$f_1 = {\displaystyle \sum_{|\textbf{l}| \leq m-D}^{}} y_{\textbf{l}}x_1^{l_1}\cdots x_n^{l_n} \quad \mbox{and} \quad g_1 = {\displaystyle \sum_{|\textbf{t}| \leq m-D}^{}} z_{\textbf{t}}x_1^{t_1}\cdots x_n^{t_n};$$ and $$f = {\displaystyle \sum_{|\textbf{j}| \leq D}^{}} a_{\textbf{j}}x_1^{j_1}\cdots x_n^{j_n} \quad \mbox{and} \quad g = {\displaystyle \sum_{|\textbf{j}| \leq D}^{}} b_{\textbf{j}}x_1^{j_1}\cdots x_n^{j_n}.$$ Our next objective is to obtain explicit height bounds for the coefficients $ y_{\textbf{l}}, ~z_{\textbf{t}}.$ As in the preceding case, we obtain an inhomogeneous linear system  $AY = b$ analogous to  \eqref{eq:A4} or \eqref{eq:A5}, where $A$ is an  $r_1 \times r_2$ matrix , with $$r_1 := \binom{n+m}{n} \quad \mbox{and} \quad r_2 := 2\binom{n+m-D}{n}$$ and whose entries belong to  $\{0, a_{\textbf{j}}, b_{\textbf{j}}; ~~ |\textbf{j}| \leq D\}$. Since this system is solvable, we can explicitly construct a solution as follows.  Let $\Delta$ be a $r_3\times r_3$-submatrix of $A$ with $\delta = {\rm det}(\Delta) \neq 0,$ where $r_3= {\rm rank}(A) \geq 1$. Multiplying the corresponding entries of $b$ by the adjugate matrix  $\Delta^{ad}$ of $\Delta$, we obtain a vector $d: = (d_1, \cdots, d_{r_3})$. Consequently, up to a permutation of the coordinates corresponding to the position of $\Delta$ in $A$, a solution of the system is given by 
 \begin{equation}\label{eq:o}
\left\{
\begin{aligned}
d_j/ \delta ~~ & \mbox{ if } ~~  \mbox{ for}~~ 1 \leq j \leq r_3, \\
	0 ~~ & \mbox{ if } ~~  \mbox{ for}~~ r_3 < j \leq r_2.
\end{aligned}
\right.
\end{equation}
This gives us an explicit expression of $ y_{\textbf{l}}, ~z_{\textbf{t}}$. Notice that for $\rho \in \{y_{\textbf{l}}, ~z_{\textbf{t}}\},$ we have 
\begin{equation*}
\begin{split}
h(\rho)
& ~\leq ~ h \left(\frac{d_j}{\delta} \right) ~ \leq~~ h(d_j) + h(\delta) \\
& ~\leq~ r_2 \log r_2 + 2r_2 h^{hom}(f, g, 1)~ \leq ~~3(2m)^n \left( \log (3(2m)^n )+ 2h^{hom}(f,g, 1) \right)\\
& ~\leq~ 3n(2m)^n \left( \log (6m)+ 2h^{hom}(f,g, 1) \right), 
\end{split}
\end{equation*}
where we used the Lemma \ref{lemA} and the fact that $m \geq n$ and  $\binom{n+m}{n} \leq (2m)^n$. 
For $v \in M_K \setminus S$ and $\textbf{u} \in \mathbb{G}_m^n(K)_{S, \eta}$, we have
$$1= |f(\textbf{u}) f_1(\textbf{u}) + g(\textbf{u}) g_1(\textbf{u})|_v \leq {\rm max} \{|f(\textbf{u})|_v, |g(\textbf{u})|_v \} \cdot {\rm max} \{|f_1(\textbf{u})|_v, |g_1(\textbf{u})|_v \}.$$ Therefore 
\begin{equation*}\label{eq:q}
\begin{split}
- \log^{-}{\rm max} \{|f(\textbf{u})|_v, |g(\textbf{u})|_v \} & ~\leq ~ - \log^{-}{\rm min} \{1/|f_1(\textbf{u})|_v, 1/|g_1(\textbf{u})|_v \}\\
& ~\leq ~ - {\rm min} \{-\log |f(\textbf{u})|_v, -\log |g(\textbf{u})|_v, 0\}  \\
& ~\leq~ {\rm max} \{\log |f(\textbf{u})|_v, \log |g(\textbf{u})|_v, 0\}\\
& ~\leq~ {\rm max}_{\textbf{i}}\{\log |y_{\textbf{i}}|_v, \log |z_{\textbf{i}}|_v, 0\}
+ m{\displaystyle \sum_{i=1}^{n}}{\rm max} \{0, \log |u_i|_v\}  \\
& ~\leq~ {\displaystyle \sum_{|\textbf{i}| \leq m-D }^{n}}{\rm max} \{0, \log |y_{\textbf{i}}|_v \} + {\displaystyle \sum_{|\textbf{i}| \leq m-D }^{}}{\rm max} \{0, \log |y_{\textbf{i}}|_v \}  \\
& + m{\displaystyle \sum_{i=1}^{n}}{\rm max} \{0, \log |u_i|_v\}.
\end{split}
\end{equation*}
Summing over $v \in M_K \setminus S$ we obtain
\begin{equation*}\label{eq:q4}
\begin{split}
- {\displaystyle \sum_{v \in M_K \setminus S}^{}} \log^{-}{\rm max} \{|f(\textbf{u})|_v, |g(\textbf{u})|_v \} & ~\leq ~ m{\displaystyle \sum_{i=1}^{n}} ~{\displaystyle \sum_{v \in M_K \setminus S}^{}} {\rm max} \{0, \log |u_i|_v\}\\
& ~+ ~  {\displaystyle \sum_{|\textbf{i}| \leq m-D }^{}} ~ {\displaystyle \sum_{v \in M_K \setminus S}^{}} {\rm max} \{0, \log |y_{\textbf{i}}|_v \} \\
&+   {\displaystyle \sum_{|\textbf{i}| \leq m-D }^{}}~ {\displaystyle \sum_{v \in M_K \setminus S}^{}} {\rm max} \{0, \log |z_{\textbf{i}}|_v \} \\
& ~\leq~ mn\eta {\displaystyle \sum_{i=1}^{n}}h(u_i)  + 2 \binom{n+m-D}{n} h(\rho)\\
& ~\leq~  mn\eta {\displaystyle \sum_{i=1}^{n}}h(u_i) + 2(2m)^n h(\rho)\\
& ~\leq~  mn\eta {\displaystyle \sum_{i=1}^{n}}h(u_i) + 6n(2m)^{2n} \left( \log (6m)+ 2h^{hom}(f,g,1) \right),\\
\end{split}
\end{equation*}
where for the second inequality, we used the fact that $${\displaystyle \sum_{v \in M_K \setminus S}^{}} {\rm max} \{0, \log |u_i|_v\} \leq \eta h(\textbf{u}) \leq \eta {\displaystyle \sum_{i=1}^{n}}h(u_i)$$ and for the last one the fact that $$h(\rho) \leq 3n(2m)^n \left( \log (6m)+ 2h^{hom}(f,g,1) \right).$$
We deduce  that $$ - {\displaystyle \sum_{v \in M_K \setminus S}^{}} \log^{-}{\rm max} \{|f(\textbf{u})|_v, |g(\textbf{u})|_v \} \leq \frac{\epsilon}{2}{\displaystyle \sum_{i=1}^{n}}h(u_i) + 6n(2m)^{2n} \left( \log (6m)+ 2h^{hom}(f,g,1) \right).$$  Since a bound of the exceptional vector given in \eqref{eq:B1} also contains those satisfying $${\displaystyle \sum_{i=1}^{n}}h(u_i) < 2 \epsilon^{-1} 6n(2m)^{2n} \left( \log (6m)+ 2h^{hom}(f,g,1) \right)$$ it follows that $$ 6n(2m)^{2n} \left( \log (6m)+ 2h^{hom}(f,g,1) \right) \leq \frac{\epsilon}{2}{\displaystyle \sum_{i=1}^{n}}h(u_i)$$ and so 
$$- {\displaystyle \sum_{v \in M_K \setminus S}^{}} \log^{-}{\rm max} \{|f(\textbf{u})|_v, |g(\textbf{u})|_v \} \leq \epsilon {\displaystyle \sum_{i=1}^{n}}h(u_i)$$ which give us the desired result.
\end{proof}

We now prove an effective version of the $S$-part inequality in Theorem \ref{thmJ}, in the sense that we obtain an explicit bound for the number of hypersurfaces whose union contains all exceptional solutions.

\begin{pro}\label{pro4}
Let $K$ be a number field of degree $d$ and $S$ be a finite set of places of $K$ containing all archimedean ones. We denote the cardinality of $S$ by $s$. Let $f \in K[x_1, \cdots, x_n]$ be a polynomials  of total degree $D$  that does not vanish at the origin $(0, \cdots, 0).$ Let $ 0< \epsilon < 1 $ be a real number. For $$\eta < \dfrac{\epsilon}{24Dn} $$ we Consider the inequality 
\begin{equation}\label{eq:A7}
- {\displaystyle \sum_{v \in  S}^{}} \log ^-  |f(u_1, \cdots, u_n)|_v  > \epsilon {\displaystyle \sum_{i=1}^{n}} h(u_i)
\end{equation}
with $(u_1, \cdots, u_n) \in \mathbb{G}_{m}^{n}(K)_{S, \eta}$. There are hypersurfaces $T_1, \cdots , T_{t_2(\epsilon)}$ of $K^n$ of total degree at most $m: = 4Dn$ with $$t_2(\epsilon)  :=  {\rm exp}\left( 125ds(2mD)^{n}(1+\log \epsilon^{-1})\right)\left( \max \{1, h^{hom}(f,1)\}\right)^{3nd}$$ such that any solution $(u_1, \cdots, u_n) \in \mathbb{G}_{m}^{n}(K)_{S, \eta}$ of \eqref {eq:A7} belongs to $$T_1 \bigcup  T_2 \bigcup \cdots \bigcup T_{t_2(\epsilon)}.$$
\end{pro}

\begin{proof}
Our approach provides an effective version of the method developed in \cite{X}, which was originally introduced by Levin \cite{L}. Without loss of generality, we may assume that $(u_1, \cdots, u_n) \in \mathbb{G}_{m}^{n}(K)_{S, \eta}$ does not lie on the hypersurface  $$H : = \{(x_1, \cdots, x_n) \in K^n~:~ f(x_1, \cdots, x_n)= 0\}.$$ At the end of the proof, we will add $1$ to our upper bound for the number of exceptional hypersurfaces to account for $H$, thereby obtaining the desired conclusion.\\
For each subset  $S_1 \subseteq S$, let $R_{S_1}$ be the set of points $(u_1, \cdots, u_n) \in \mathbb{G}_{m}^{n}(K)_{S, \eta}$ such that $$S_1 = \{v \in S ~:~ \log |f(u_1, \cdots, u_n)|_v <0 \}.$$ Then, for $(u_1, \cdots, u_n) \in R_{S_1},$ we have 
\begin{equation*}
\log^{-}|f(u_1, \cdots, u_n)|_v = \left\{ \begin{aligned}
\log|f(u_1, \cdots, u_n)|_v ~~~~~~~& \mbox{ if } ~~~~~~v \in S_1,\\ 0~~~~~~ & \mbox{ if }~~~~~~ v \in S \setminus S_1. 
\end{aligned}
\right.
\end{equation*}
Since there are at most $2^s$ possible subsets $S_1 \subseteq S$, we may treat each of them separately. We set $m:=4Dn\epsilon^{-1}$ and consider the morphism   $\psi ~:~ \mathbb{P}^n \rightarrow \mathbb{P}^N$, $\psi = (\psi_0, \cdots, \psi_N),$ where $N = \binom{n+mD}{n}-1$ and the coordinate functions $\psi_0, \cdots, \psi_N$ are given by the monomials of degree $mD$ in $K[x_0, \cdots, x_n]$ arranged in lexicographic order. Let $V_{mD}$ be the vector space of homogeneous polynomial of degree $mD$ and $B_{mD}$ be its standard monomial basis, consisting of all monomials of degree  $mD$. 
We write $$X^{\textbf{i}} = x_0^{i_0}\cdots x_n^{i_n}.$$  For $v \in S_1,$ we construct a basis of $V_{mD}$  as follows. For every multi-index $\textbf{i}$ with $|\textbf{i}| = mD,$ set
$k_{\textbf{i}}= \lfloor \dfrac{{\rm ord}_{x_0} X^{\textbf{i}}}{D} \rfloor$ and define $$B_v^{\textbf{i}} = \dfrac{X^{\textbf{i}}}{x_0^{k_{\textbf{i}}D}}F^{k_{\textbf{i}}}\quad \mbox{with}\quad  F= x_0^Df(x_1/x_0, \cdots, x_n/x_0).$$ Let  $L_v^{\textbf{i}} $ be the linear form over $K$ such that $$L_v^{\textbf{i}}( X^{\textbf{j}}, |\textbf{j}| = mD )= \dfrac{X^{\textbf{i}}}{x_0^{k_{\textbf{i}}D}} F^{k_{\textbf{i}}}.$$ Let $B_v$ be the set of all $L_v^{\textbf{i}}$. Since $f$ does not vanish at the origin, the coefficient of  $x_0^D$ in $F$ is nonzero. Consequently, the linear forms  $L_v^{\textbf{i}},$  with $|\textbf{i}| = mD$ are linearly independent over $K.$ 

For every $v \in S \setminus S_1$,  we  instead define $L_v^{\textbf{i}} $ be the linear form over $K$ such that $L_v^{\textbf{i}}( X^{\textbf{j}}, |\textbf{j}| = mD ) = X^{\textbf{i}}$. In this case $B_v = B_{mD}$. Therefore we deduce that $$\log | {\rm det} (L_v^{\textbf{i}}, ~ |\textbf{i}| = mD )|_v =0$$ for all  $v \in S \setminus S_1$. Assume that  $v \in S_1$. We  denote by $(a_{\textbf{i}},~ |\textbf{i}| = mD)$ be the coefficients of polynomial $F$ and  $(b_{\textbf{i},~ \textbf{j}}, |\textbf{j}| = mD)$ those of polynomial $L_v^{\textbf{i}}$. From  Gelfond's lemma (see \cite[Lemma 4.1.3]{evertse-gyory-2022}, we deduce that
$$\max_{\textbf{j}} |b_{\textbf{i},~ \textbf{j}}|_v \leq 2^{mnDs(v)} \max_{\textbf{i}} |a_{\textbf{i}}|_v^m.$$ 

Hence we obtain 
\begin{equation*}
\begin{split}
\log | {\rm det} (L_v^{\textbf{i}}, ~ |\textbf{i}| = mD )|_v &~\leq ~ Ns(v)\log N + N \log \max _{\textbf{j}} |b_{\textbf{i},~ \textbf{j}}|_v\\
&~\leq~ s(v)(N\log N + NmnD) + mN \log  \max_{\textbf{i}} |a_{\textbf{i}}|_v
\end{split}
\end{equation*}
By summing over all $v \in M_K$, we get
\begin{equation}\label{eq:C}
h \left( {\rm det} (L_v^{\textbf{i}}, ~ |\textbf{i}| = mD )\right)  \leq 3(N+1)mnD + (N+1)m h^{hom}(f,1). 
\end{equation}

We now  apply Lemma \ref{thm1} with $\delta = \epsilon /3$ and deduce that There exist hypersurfaces of degree at most $m$ denoted  $H_1, \ldots, H_{l}$ of $K^n,$ with $$l < \left( 2^{60(N+1)^2} \left(  \dfrac{\epsilon}{3} \right) ^{7(N+1)}\right)^{s+1}(6m^{1/2} + \log\log \max \{4, mh^{hom}(f,1)\})$$ such that every  $(S, \eta)$-unit which is not in $H_1 \bigcup \cdots \bigcup H_l$ satisfy
\begin{equation*}\label{eq:w}
\begin{split}
{\displaystyle \sum_{v \in S}^{}} {\displaystyle \sum_{Q \in B_v}^{}}\log \dfrac{|\psi (P)|_v}{|Q(P)|_v} &~\leq~ \left(N+1 + \dfrac{\epsilon}{3} \right)\left(h^{hom}(\psi(P))+\dfrac{\log (N+1)}{2} \right)\\
&~ - {\displaystyle \sum_{v \in S}^{}} \log |{\rm det}(L_{\textbf{i}}^{v},~| \textbf{i}| = mD)|_v, 
\end{split}
\end{equation*}
where $P = (1, u_1, \cdots, u_n)$. Together with Lemma \ref{lemA}a) and  \eqref{eq:C} give us 
\begin{equation}\label{eq:E}
\begin{split}
{\displaystyle \sum_{v \in S}^{}}~ {\displaystyle \sum_{Q \in B_{mD}}^{}}\log \dfrac{|\psi (P)|_v}{|Q(P)|_v} - {\displaystyle \sum_{v \in S'}^{}} ~{\displaystyle \sum_{|\textbf{i}| = mD}^{}}\log \dfrac{|B_v^{\textbf{i}}(P)|_v}{|X^{\textbf{i}}(P)|_v} &~\leq~\\
\left(N+1 + \dfrac{\epsilon}{3} \right)\left(h^{hom}(\psi(P))+\dfrac{\log (N+1)}{2} \right)  + 4(N+1)^2 \max\{h^{hom}(f, 1), 1\}.
\end{split}
\end{equation}
From now, we estimate $$- {\displaystyle \sum_{v \in S'}^{}} {\displaystyle \sum_{|\textbf{i}| \leq mD}^{}}\log \dfrac{|B_v^{\textbf{i}}(P)|_v}{|X^{\textbf{i}}(P)|_v}.$$ To this end, let $(u_1, \cdots, u_n) \in R_{S_1}$. Then we have 
\begin{equation}\label{eq:y}
\begin{split}
- {\displaystyle \sum_{v \in S'}^{}} \log \dfrac{|B_v^{\textbf{i}}(P)|_v}{|X^{\textbf{i}}(P)|_v} &~=~ -k_{\textbf{i}} {\displaystyle \sum_{v \in S'}^{}} \log |f(u_1, \cdots, u_n)|_v \\ 
&~ =~  -k_{\textbf{i}} {\displaystyle \sum_{v \in S}^{}} \log^{-} |f(u_1, \cdots, u_n)|_v.
\end{split}
\end{equation}
Letting $I= {\displaystyle \sum_{|\textbf{i}| = mD}^{}} k_{\textbf{i}} $. By summing over all $\textbf{i}$, we get $$- {\displaystyle \sum_{v \in S'}^{}} {\displaystyle \sum_{|\textbf{i}| = mD}^{}} \log \dfrac{|B_v^{\textbf{i}}(P)|_v}{|X^{\textbf{i}}(P)|_v} = -I {\displaystyle \sum_{v \in S}^{}} \log^{-} |f(u_1, \cdots, u_n)|_v.$$ It is clear that 
\begin{equation*}
\begin{split}
I= {\displaystyle \sum_{|\textbf{i}| = mD}^{}} k_{\textbf{i}} \geq  {\displaystyle \sum_{|\textbf{i}| = mD}^{}}  \dfrac{{\rm ord}_{x_0} X^{\textbf{i}}}{D} - \binom{n+mD}{n} \\
~\geq~ \frac{1}{D} \binom{n+mD}{n+1} -  \binom{n+mD}{n} \\  
~ \geq \left( \dfrac{m-n-1}{n+1}\right)  \binom{n+mD}{n}.
\end{split}
\end{equation*}
Therefore, we finally get:
\begin{equation}\label{eq:D1}
- {\displaystyle \sum_{v \in S'}^{}} {\displaystyle \sum_{|\textbf{i}| \leq mD}^{}} \log \dfrac{|B_v^{\textbf{i}}(P)|_v}{|X^{\textbf{i}}(P)|_v} \geq  - \left( \dfrac{m-n-1}{n+1}\right)  \binom{n+md}{n} {\displaystyle \sum_{v \in S}^{}} \log^{-} |f(u_1, \cdots, u_n)|_v.
\end{equation}
Now, we estimate the quantity $${\displaystyle \sum_{v \in S}^{}} {\displaystyle \sum_{Q \in B_{mD}}^{}}\log \dfrac{|\psi (P)|_v}{|Q(P)|_v} = {\displaystyle \sum_{v \in M_K }^{}} {\displaystyle \sum_{Q \in B_{mD}}^{}}\log \dfrac{|\psi (P)|_v}{|Q(P)|_v} - {\displaystyle \sum_{v \in M_K \setminus S}^{}} {\displaystyle \sum_{Q \in B_{mD}}^{}}\log \dfrac{|\psi (P)|_v}{|Q(P)|_v}.$$ In one hand, we have $${\displaystyle \sum_{v \in M_K \setminus S}^{}} \log |\psi(P)|_v \leq mD {\displaystyle \sum_{v \in M_K \setminus  S}^{}} \log \max_i \{1, |u_i|_v\}.$$ On the other hand, $${\displaystyle \sum_{v \in M_K \setminus S}^{}} \log \left| \dfrac{1}{\psi(P)} \right|_v \leq mD {\displaystyle \sum_{v \in M_K \setminus  S}^{}} \log \max_i \left \{1,  \left| \dfrac{1}{u_i} \right|_v \right \}.$$ All together with the fact that $(u_1, \cdots, u_n) \in \mathbb{G}_m^n(K)_{S, \eta}$, we obtain  $${\displaystyle \sum_{v \in M_K \setminus S}^{}} \log |\psi(P)|_v  + \log \left| \dfrac{1}{\psi(P)} \right|_v  \leq mD \eta h(P).$$ It is clear that
\begin{equation}\label{eq:D2}
h^{hom}(\psi(P)) = mD h^{hom}(P,1).
\end{equation}
From product formula, we have $${\displaystyle \sum_{v \in M_K}^{}} {\displaystyle \sum_{Q \in B_{mD}}^{}}\log |Q(P)|_v =0.$$ Therefore 
\begin{equation*}
\begin{split}
{\displaystyle \sum_{v \in M_K}^{}} {\displaystyle \sum_{Q \in B_{mD}}^{}}\log \dfrac{|\psi (P)|_v}{|Q(P)|_v} = (N+1) {\displaystyle \sum_{v \in M_K}^{}}\log |\psi(P)|_v  \\
 \geq (N+1)h^{hom}(\psi(P)) \geq mD (N+1)h^{hom}(P,1),
\end{split}
\end{equation*}
 where for the last inequality, we used relation \ref{eq:D2}. Moreover 
 \begin{equation*}
\begin{split}
 {\displaystyle \sum_{v \in   M_K \setminus S}^{}}~ {\displaystyle \sum_{Q \in B_{mD}}^{}}\log \dfrac{|\psi (P)|_v}{|Q(P)|_v} \leq (N+1)\left( {\displaystyle \sum_{v \in M_K \setminus S}^{}} \log |\psi(P)|_v  + \log \left| \dfrac{1}{\psi(P)} \right|_v \right)\\
 \leq mD(N+1)\eta h^{hom}(P,1).
 \end{split}
\end{equation*}
All together give us 
 $${\displaystyle \sum_{v \in S}^{}} ~{\displaystyle \sum_{Q \in B_{mD}}^{}}\log \dfrac{|\psi (P)|_v}{|Q(P)|_v}  \geq mD(N+1)(1- \eta)h^{hom}(P,1).$$
 Together with relations \eqref{eq:D1} and \eqref{eq:A7}, we deduce that 
 \begin{equation*}
\begin{split}
 - \left( \dfrac{m-n-1}{n+1}\right)  \binom{n+mD}{n} {\displaystyle \sum_{v \in S}^{}} \log^{-} |f(u_1, \cdots, u_n)|_v + mD(N+1)(1-\eta)h^{hom}(P,1) \\
  \leq ~\left(N+1+  \dfrac{\epsilon}{3} \right)\left(h^{hom}(\psi(P))+\dfrac{\log (N+1)}{2} \right) + 4(N+1)^2\max \{1, h^{hom}(f,1)\}
\end{split}
\end{equation*}
From \eqref{eq:D2} and the fact that $$h^{hom}(P,1) \leq {\displaystyle \sum_{i=1}^{n}} h(u_i),$$  we deduce  
 \begin{equation*}
\begin{split}
- \left( \dfrac{m-n-1}{n+1}\right) {\displaystyle \sum_{v \in S}^{}} \log^{-} |f(u_1, \cdots, u_n)|_v \leq  \left( mD\eta + \frac{\epsilon}{3} \right) {\displaystyle \sum_{i=1}^{n}} h(u_i) \\
+ 6(N+1)\max \{1, h^{hom}(f,1)\}.
\end{split}
\end{equation*}

Using the fact that $m= 4Dn,$ we deduce that 
\begin{equation}\label{eq:A9}
\left( \dfrac{n+1}{m-n-1}mD\eta + \dfrac{\epsilon}{3} \right) \leq \frac{\epsilon}{2}
\end{equation}
 whenever $$\eta < \dfrac{\epsilon}{24Dn}.$$ So
\begin{equation*}
\begin{split}
-  {\displaystyle \sum_{v \in S}^{}} \log^{-} |f(u_1, \cdots, u_n)|_v  
~ \leq~  \dfrac{\epsilon}{2}  {\displaystyle \sum_{i=1}^{n}} h(u_i) 
+ 6(2mD)^n\max \{1, h^{hom}(f,1)\}.
\end{split}
\end{equation*}
From Lemma \ref{lemE}, there are at most 
\begin{equation}\label{eq:B5}
t_3:=  5^n\left(24\epsilon^{-1}(2mD)^n \right)^{3nd}\left( \max \{1, h^{hom}(f,1)\}\right)^{3nd} 
\end{equation}
vectors of $K^n$ such that $${\displaystyle \sum_{i=1}^{n}} h(u_i) \leq 2\epsilon^{-1}6(2mD)^n\max \{1, h^{hom}(f,1)\}.$$
Excepted those vectors, we obtain 
\begin{equation*}
\begin{split}
-  {\displaystyle \sum_{v \in S}^{}} \log^{-} |f(u_1, \cdots, u_n)|_v 
~ \leq~  \epsilon {\displaystyle \sum_{i=1}^{n}} h(u_i).
\end{split}
\end{equation*}
Since the choice of linear form in the application of the  quantitative subspace theorem used so far only depend on the choice of the monomial basis  $S_1$ which are at most $2^s$, we conclude that there are hypersurfaces $T_1, \cdots , T_{t_2(\epsilon)}$ of $K^n$  each of total degree at most  $m := 4nD$ with
\begin{equation*}
\begin{split}
  \left( 2^{61(N+1)^2} \left(  \dfrac{\epsilon}{3} \right) ^{-7(N+1)}\right)^{s+1} \left(2m^{1/2} \max \left \{1, h^{hom}(f,1) \right \}\right)  + \\   5^n\left(24\epsilon^{-1}(2mD)^n \right)^{3nd}\left( \max \left \{1, h^{hom}(f,1) \right \}\right)^{3nd} \\
~\leq ~ {\rm exp}\left( 125ds(2mD)^{n}(1+\log \epsilon^{-1})\right)\left( \max \left \{1, h^{hom}(f,1)\right \}\right)^{3nd} =: t_2(\epsilon)
\end{split}
\end{equation*}
such that every  solution $(u_1, \cdots, u_n) \in \mathbb{G}_{m}^{n}(K)_{S, \eta}$ of \eqref {eq:A7} belongs to $$T_1 \bigcup  T_2 \bigcup \cdots \bigcup T_{t_2(\epsilon)}.$$ 
\end{proof}

\subsection{Proof of Theorem \ref{thmJ}}
By Proposition \ref{pro3} and Proposition \ref{pro4}, we deduce  that, for fixed $0 < \epsilon < 1$ and $$\eta < \dfrac{\epsilon^2}{72 D^2 N^4},$$ there exist at most $$t :=  2{\rm exp}\left( 75sNd(2m)^{2N}(1+\log \epsilon^{-1})\right)\left( \log (6m)+ 2h^{hom}(f,g, 1) \right)^{3Nd}$$ hypersurfaces  $H_1, \cdots, H_t$ of $\mathbb{G}_m^N(K)$ of degree at most $m := 12D^2N^3 \epsilon ^{-1}$ such that any solution $(u_1, \cdots, u_N) \in \mathbb{G}_{m}^{N}(K)_{S, \eta}$ of the inequality
\begin{equation*}\label{eq:k}
- {\displaystyle \sum_{v \in M_K }^{}} \log ^- {\rm max} \left \{ |f(u_1, \cdots, u_N)|_v ,~ |g(u_1, \cdots, u_N)|_v \right \} > \epsilon {\displaystyle \sum_{i=1}^{N}} h(u_i)
\end{equation*}
belongs to $$Z := H_1 \bigcup H_2 \bigcup \cdots \bigcup H_t.$$

\section{Greatest common divisor and multi-recurrences.}\label{sec4}
Let $K$ be a number field of degree $d$ and let $(G(\nbf))$ be a multi-recurrence defined over $K$ as 
\begin{equation} \label{eq:1}
G(n_1, \ldots, n_t) = \sum_{i=1}^{r} P_i(n_1, \ldots, n_t)\alpha_{i1}^{n_1}\cdots \alpha_{it}^{n_t}, 
 \end{equation}
 where $r,t$ are positive integers, the $P_i$ are polynomials in $t$ variables with coefficients given by the vector $\textbf{a}_i =(a_{i1}, \ldots, a_{i\ri}) \in K^{m_i}$, ordered lexicographically, for $i=1,\ldots,r$, and $n_1, \ldots, n_t$ are non-negative integers and $\alpha_{i,j} \in K$ for all $(i,j) \in \{1, \cdots, r\} \times \{1, \cdots, t\}$ is called multi-root of $(G(\nbf))$.
 Let $S$ the set containing all nonarchimedean places corresponding to  the prime ideals above the $\alpha_{ij}$'s and all the archimedean places over $K$, and we denote the cardinality of $S$ by $s$. Let $\Gamma = \{1, \ldots,r\}$ and $\G$ be the subset of $\mathbb{Z}^t$ consisting of $\textbf{z}$ such that $\bsa_l^{\textbf{z}} = \bsa_k^{\textbf{z}}$ for some $l,k \in \Gamma$ with $l \neq k$. We denote by $$A' = \sum_{i \in \Gamma}^{} \binom{t+\ri}{t} \mbox{ and } A= \max \{t, A'\}.$$
Notice that $A = \max\{t,r\}$ if all the polynomials $P_1, \ldots, P_r$ are constants, and $A=A'$ otherwise. We write $$\mathcal{S}_i = \{ \textbf{z} \in \mathbb{Z}^t: P_i(\textbf{z}) = 0\}$$ for $i=1,\ldots,r$. We set $$h := {\rm max}_{i,j} \left \{ h(a_{i,j}), ~ 1  \right \}$$ where $a_{i,j}$ range over the coefficients of the polynomials $P_1, \cdots, P_r.$  For  convenience, we restate Theorem \ref{thmK}.

\begin{proposition}\label{propA'}(Theorem \ref{thmK})
Let $D$  be a positive real number. Assume that $\G = \{0\}$ and  $$H_i := \{(u_1, \cdots, u_t) \in \mathbb{Z}^t ~: ~P_i(u_1, \cdots, u_t) = 0\} $$ has at most $D$ elements and the $P_i$ have total degrees at most $\sigma \geq 1$.  Let $|\cdot |_v$ be an absolute value on K such that, for some $i,$ we have  $|\alpha_{i,j}|_v \geq 1$ for all $j\in \{1, \cdots, t\}$. For $0< \epsilon <1,$ there exists at most $$ \rho(\epsilon):= 4 \left( 2^{60r^2}\delta^{-7r}\right)^{s+1}\frac{r!}{2} 2^{35A^3}d^{6A^2} +D^r + L^t $$ with $$\eta =  \dfrac{\epsilon(r-1)c_6h^r}{(2r +1)(r + 6r^2t c_6h^{r+1})} \quad \mbox{and}\quad \delta =  \dfrac{\epsilon(r-1)c_6h^r}{r + 6r^2t c_6h^{r+1}}$$ and $$ L:= \max\left\{ (2t\sigma)^t\exp\left((2t\sigma)^th \right), (4 \eta^{-1}c_6h^r(\sigma+1)(r-1))^{3/2}\right\}$$ and 
$$c_6 =  \dfrac{58t!3^td^{t+1} \log d}{t^t}$$
values of $\nbf$ such that
\begin{equation}\label{eq:D9}
\log |G(\nbf)|_v < -\epsilon |\nbf|,
\end{equation}
where $|\nbf| = {\rm max}_i\{n_i\}.$
\end{proposition}

\begin{proof}
We first consider the easier case $r=1$. Let  $\mathcal{N}$ denote the set of solution $\nbf$ of  \eqref{eq:D9} satisfying $P_1(\nbf) \neq 0$.  By hypothesis, it suffices to obtain an explicit upper bound for the solutions of $\log |P_1(\nbf)|_v < -\epsilon |\nbf|.$ By Lemma \ref{lemA}, we have  $$-\log |P_1(\nbf)|_v \leq h(P_1(\nbf)) \leq \sigma \log |\nbf| + (2t \sigma)^t h + \log \binom{t+ \sigma}{t}.$$ Consequently, 
$$-\sigma \log |\nbf| \leq - \epsilon |\nbf| +(2t \sigma)^t h + \log \binom{t+ \sigma}{t}.$$ Suppose first that $$\frac{\epsilon}{2} |\nbf| \geq  (2t \sigma)^t h + \log \binom{t+\sigma}{t}.$$ Then $$\log |\nbf| \geq \frac{\epsilon}{2\sigma} |\nbf|.$$  It follows from Lemma \ref{lemH}  that $$|\nbf| \leq \dfrac{8\sigma}{\epsilon} \log \dfrac{4\sigma}{\epsilon}.$$ Consequently the cardinality of $\mathcal{N}$ does not exceed $$ D+ \left( \dfrac{8\sigma}{\epsilon} \log \dfrac{4\sigma}{\epsilon} \right)^t+ \left(\frac{2}{\epsilon}\left(  (2t \sigma)^t h + t\log (2t \sigma) \right) \right)^t.$$
We may therefore assume $r>1$. For each $w \in S \setminus \{v\}$, we define $r$ linear forms $L_{1w},\cdots ,L_{rw}$ in the variables $\xbf =(x_1, \cdots,x_r)$ by
$$L_{iw}(\xbf):= x_i, \quad i=1, \cdots ,r.$$
At the place $v$, we instead define $$ L_{1v}(\xbf):= x_1+ \cdots + x_r, \quad L_{iv}(\xbf):= x_i, ~~ \mbox{for}~~all \quad i= 2, \cdots, r.$$ Clearly, the linear forms $L_{1w},\cdots ,L_{rw}$ are linearly independent over $K$, and $$|{\rm det}(L_{1w},\cdots ,L_{rw})|_w =1.$$
 Without loss of generality we assume $|\alpha_{1,j}|_v \geq 1$ for all $j\in \{1, \cdots, t\}$. We put 
$$\xbf_{\nbf} = (P_1(\nbf)\bsa_1^{\nbf}, \ldots, P_r(\nbf)\bsa_r^{\nbf}).$$  Furthermore, set $$\eta =  \dfrac{\epsilon(r-1)c_6h^r}{(2r +1)(r + 2r^2 c_6h^{r+1})}$$ and  suppose that $\nbf$ is solution of  \eqref{eq:D9} satisfying
\begin{equation}\label{eq:A'}
|\nbf| \geq L:= \max\left\{ (2t\sigma)^{2t}\exp\left(2(2t\sigma)^th \right), (4 \eta^{-1}c_6h^r(\sigma+1)(r-1))^{3/2}\right\}.
\end{equation}
Our aim is to prove that $${\displaystyle \sum_{w \in S}}{\displaystyle \sum_{i=1}^{r}} \log \dfrac{|L_{iw}(\xbf_{\nbf})|_w}{|\xbf_{\nbf} |_w} + \log \dfrac{|G(\nbf)|_v}{|P_1(\nbf)\bsa_1^{\nbf}|_v} < -(r+ \delta) h(\xbf_{\nbf})$$ for a suitable choice of $\delta$ so that we can apply the quantitative Subspace Theorem. We first claim  that $\xbf_{\nbf} \in \mathbb{G}_m^r(K)_{S, \eta}$. Indeed, on the one hand,  we have 
\begin{equation}\label{eq:E2}
\begin{split}
{\displaystyle \sum_{w \in M_K \setminus S}}  \max \left \{ 0, \log \left|\xbf_{\nbf} \right|_v  \right \} + \max \left \{ 0, \log \left|\frac{1}{\xbf_{\nbf}} \right|_w  \right \}\\
~ =~  {\displaystyle \sum_{w \in M_K \setminus S}}  \max_i \left \{ 0, \log \left|P_i(\nbf) \right|_v \right \} + \max_i \left \{ 0, \log \left|\frac{1}{P_i(\nbf)} \right|_w  \right \}\\
~ \leq~  2{\displaystyle \sum_{i=1}^{r}} h(P_i(\nbf))~\leq ~ 2\sigma \log |\nbf| + 2\binom{t+ \sigma}{t}h + 2\log \binom{t+ \sigma}{t},
\end{split}
\end{equation}
for the first inequality, we have used the fact that $\alpha_{i,j} \in \mathcal{O}_S^*$, while the last inequality follows from Lemma \ref{lemA}. On the other hand Lemma \ref{lemF} yields $$\dfrac{|\nbf|}{c_6h^r} \leq h^{hom}(\bsa_1^{\nbf}, \cdots, \bsa_r^{\nbf}) \leq {\displaystyle \sum_{i=2}^{r}} h \left( \dfrac{P_i(\nbf)\bsa_i^{\nbf}}{P_1(\nbf)\bsa_1^{\nbf}}\right)+ {\displaystyle \sum_{i=2}^{r}} h \left( \dfrac{P_1(\nbf)}{P_i(\nbf)}\right).$$ Furthermore, by the definition of the homogeneous height, we have $$h^{hom}(\xbf_{\nbf}) \geq h \left( \dfrac{P_i(\nbf)\bsa_i^{\nbf}}{P_1(\nbf)\bsa_1^{\nbf}}\right).$$ Moreover, Lemma \ref{lemA} gives $$ h(P_i(\nbf)) \leq \sigma \log |\nbf| + \binom{t+ \sigma}{t}h + \log \binom{t+ \sigma}{t}$$ for all $i \in \{1, \cdots, r\}.$ Combining these estimates, we obtain
\begin{equation}\label{eq:F}
\begin{split}
h^{hom}(\xbf_{\nbf}) &~\geq~ \dfrac{|\nbf|}{(r-1)c_6h^r} - 2\sigma \log |\nbf| - 2\binom{t+ \sigma}{t}h - 2\log \binom{t + \sigma}{t} \\
&~ \geq~  \dfrac{|\nbf|}{2(r-1)c_6h^r}.
\end{split}
\end{equation}
where  we have combined  Lemma \ref{lemH} and the relation \eqref{eq:A'} to get $$\dfrac{|\nbf|}{2(r-1)c_6h^r} \geq 2\sigma \log |\nbf| + 2\binom{t+ \sigma}{t}h + 2\log \binom{t + \sigma}{t} $$ which implies the last inequality. It is also clear from Lemma \ref{lemH} and the relation \eqref{eq:A'} that
\begin{equation}\label{eq:G}
\eta  \dfrac{|\nbf|}{2(r-1)c_6h^r} \geq 2\sigma \log |\nbf| + 2\binom{t+ \sigma}{t}h + 2\log \binom{t + \sigma}{t}.
\end{equation}
Together with \eqref{eq:E2} and \eqref{eq:F}, we deduce that $${\displaystyle \sum_{w \in M_K \setminus S}}  \max \left \{ 0, \log \left|\xbf_{\nbf} \right|_v  \right \} + \max \left \{ 0, \log \left|\frac{1}{\xbf_{\nbf}} \right|_w  \right \} \leq \eta  \dfrac{|\nbf|}{2(r-1)c_6h^r} \leq \eta h^{hom}(\xbf_{\nbf}).$$ Therefore the claim is proven. For $i = 1, \cdots, r$ we have 
\begin{equation}\label{eq:F2}
\begin{split}
{\displaystyle \sum_{w \in S}}  \log \dfrac{|P_i(\nbf) \bsa_{\nbf}|_w}{|\xbf_{\nbf} |_w} &~=~ -{\displaystyle \sum_{w \in M \setminus S}} \log |P_i(\nbf)|_w - h(\xbf_{\nbf}) +  {\displaystyle \sum_{w \in M \setminus S}}  \log |\xbf_{\nbf}|_w \\  
&~\leq ~ h(P_i(\nbf)) + (\eta -1)h(\xbf_{\nbf}).
\end{split}
\end{equation}
Therefore $${\displaystyle \sum_{w \in S}} {\displaystyle \sum_{i=1}^{r}}  \log \dfrac{|P_i(\nbf) \bsa_{\nbf}|_w}{|\xbf_{\nbf} |_w}  \leq {\displaystyle \sum_{i=1}^r} h(P_i(\nbf))  + r(\eta -1)h(\xbf_{\nbf}).$$ On the other hand $$ \log |P_1(\nbf) \bsa_1^{\nbf}|_v \geq \log |P_1(\nbf)|_v \geq -h(P_1({\nbf})),$$ where for the first inequality, we used the fact that $\min_{j} \{|\alpha_{1j}|_v\} \geq 1$. All together give us $${\displaystyle \sum_{w \in S}} {\displaystyle \sum_{i=1}^{r}}  \log \dfrac{|P_i(\nbf) \bsa_{\nbf}|_w}{|\xbf_{\nbf} |_w}  + \log \dfrac{|G(\nbf)|_v}{|P_1(\nbf)\bsa_1^{\nbf}|_v} \leq  \eta  \dfrac{(r+1)|\nbf|}{2(r-1)c_6h^r} + r(\eta -1)h(\xbf_{\nbf}) - \epsilon |\nbf|,$$ where we used the relation \eqref{eq:G}. Therefore, we deduce from \eqref{eq:F} that  
\begin{equation}\label{eq:F3}
\begin{split}
 {\displaystyle \sum_{w \in S}} {\displaystyle \sum_{i=1}^{r}}  \log \dfrac{|P_i(\nbf) \bsa_{\nbf}|_w}{|\xbf_{\nbf} |_w}  + \log \dfrac{|G(\nbf)|_v}{|P_1(\nbf)\bsa_1^{\nbf}|_v} ~\leq~  \eta (r+1)h(\xbf_{\nbf}) +  r(\eta -1) h(\xbf_{\nbf}) - \epsilon |\nbf|.
\end{split}
\end{equation}
Since $$h(\xbf_{\nbf}) \leq \frac{\log r}{2} + h^{hom}(\xbf_{\nbf}) \leq  \frac{\log r}{2} + {\displaystyle \sum_{i=1}^{r}} h \left( \dfrac{P_i(\nbf)\bsa_i^{\nbf}}{P_1(\nbf)\bsa_1^{\nbf}}\right) \leq \left ( \dfrac{ r + 6r^2t c_6h^{r+1}}{2(r-1)c_6h^r} \right) |\nbf|,$$ it follows from \eqref{eq:F3} that 
$$ {\displaystyle \sum_{w \in S}} {\displaystyle \sum_{i=1}^{r}}  \log \dfrac{|P_i(\nbf) \bsa_{\nbf}|_w}{|\xbf_{\nbf} |_w}  + \log \dfrac{|G(\nbf)|_v}{|P_1(\nbf)\bsa_1^{\nbf}|_v} < - \left( 2\dfrac{\epsilon(r-1)c_6h^r}{r + 6r^2t c_6h^{r+1}} - \eta (2r+1) +r \right) h(\xbf_{\nbf}).$$ Using the fact that $$\eta =  \dfrac{\epsilon(r-1)c_6h^r}{(2r +1)(r + 6r^2t c_6h^{r+1})},$$ we deduce that $${\displaystyle \sum_{w \in S}} {\displaystyle \sum_{i=1}^{r}}  \log \dfrac{|P_i(\nbf) \bsa_{\nbf}|_w}{|\xbf_{\nbf} |_w}  + \log \dfrac{|G(\nbf)|_v}{|P_1(\nbf)\bsa_1^{\nbf}|_v} < -(\delta + r) h(\xbf_{\nbf})$$ with $$ \delta =  \dfrac{\epsilon(r-1)c_6h^r}{r + 6r^2t c_6h^{r+1}}.$$ Moreover $$ {\displaystyle \sum_{w \in S}} \log |{\rm det}(L_{1w}, \cdots, L_{rw})|_w =0.$$ Applying Lemma \ref{thm1} with this choice of $\delta$, we deduce that there are at most  $$4 \left( 2^{60r^2}\delta^{-7r}\right)^{s+1}$$ polynomial-exponential Diophantine equations of the shape 
\begin{equation} \label{eq:1A}
\sum_{i \in I}^{} c_iP_i(\nbf)\bsa_i^{\nbf}=0
\end{equation}
with $(c_1, \ldots, c_r) \in K^r / \{\textbf{0}\}$ and $I \subseteq \Gamma$ such that every solution $\nbf$ of \eqref{eq:D9} satisfies at least one of these equations.  Fix an equation of the form \eqref{eq:1A}, and let $\mathcal{P}$ be a partition of $I$. The set $\lambda \subseteq I$ occurring in the partition $\mathcal{P}$ will be considered elements of $\mathcal{P}$. To each such partition, we associate the system of equations  
\begin{equation} \label{eq:12}
\sum_{i \in \lambda}^{} c_iP_i(\nbf)\bsa_i^{\nbf}=0~(\lambda \in \mathcal{P}).
\end{equation} 
A solution $\nbf$ of \eqref{eq:12} is called $\mathcal{P}$-degenerate if a subsum of one of the equations of the system vanishes. Otherwise, we will say $\nbf$ is $\mathcal{P}$-non-degenerate. Let $M(\mathcal{P})$ be the set consisting of $\mathcal{P}$-non-degenerate solutions of the system \eqref{eq:12}. Clearly, every solution of \eqref{eq:1A} lies in $M(\mathcal{P})$ for some partition $\mathcal{P}$. By \cite[Lemma 2.5]{Le}, the set of partitions of $I$ which only contain subsets with at least two elements has a cardinality, which does not exceed $r! / 2.$ Hence, by using the fact that $\G = \{\textbf{0}\}$ and applying Lemma \ref{thm2}, we deduce that the number of solutions of \eqref{eq:1A} is at most  $$\frac{r!}{2} 2^{35A^3}d^{6A^2}.$$ Combining this estimate with the preceding upper bound for the number of polynomial-exponential Diophantine equations of the form \eqref{eq:1A}, we conclude that, for $ |\nbf| \geq L$, the inequality \eqref{eq:D9} has at most  $$T_1:= 4 \left( 2^{60r^2}\delta^{-7r}\right)^{s+1}\frac{r!}{2} 2^{35A^3}d^{6A^2} +D^m$$ for all $|\nbf| \geq L$. Finally, taking into account the solutions satisfying $| \nbf | < L$, we obtain the following upper bound for the total number of solutions of \eqref{eq:D9}:  $$l:= 4 \left( 2^{60r^2}\delta^{-7r}\right)^{s+1}\frac{r!}{2} 2^{35A^3}d^{6A^2} +D^m + L^t.$$ This completes the proof of the desired result.
\end{proof}

Let $\Gamma$ be a torsion-free multiplicative subgroup of $K^*$ of rank $m \geq 1.$  It is readily verified that the set of multi-recurrences of the form \eqref{eq:A1}, defined over $K$ and having all their  multi-roots in $\Gamma$, is closed under addition and multiplication. It therefore forms a ring, which we denote by $R_{\Gamma}$. Let $\{u_1, \cdots,u_m\}$ be a basis of $\Gamma$.  It is straightforward to verify that $R_{\Gamma}$ is isomorphic to the Laurent polynomial ring $$\mathbb{K}= K[x_1, \cdots, x_t, T_{1,1}^{\pm 1}, \cdots,  T_{m,t}^{\pm 1} ]$$ with the isomorphism induced by mapping $x_i$ to the multi-recurrence $F(\nbf)=n_i$ and $T_{i,j}$ to the multi-recurrence $F(\nbf)=u_i^{n_j}$. Hence we say that two elements $F$ and $G$ of $R_{\Gamma}$ are coprime if their corresponding polynomials are coprime in $\mathbb{K}$.\\

Let $K$ be a number field of degree $d$ and  $(F(\nbf)), ~(G(\nbf))$ be two multi-recurrences defined over $K$ by 
$$
F(n_1, \ldots, n_t) = \sum_{i=1}^{r} Q_i(n_1, \ldots, n_t)c_{i1}^{n_1}\cdots c_{it}^{n_t}, 
$$
and
$$
G(n_1, \ldots, n_t) = \sum_{i=1}^{r} P_i(n_1, \ldots, n_t)\alpha_{i1}^{n_1}\cdots \alpha_{it}^{n_t}, 
$$
 where $r$ and $ t$ are positive integers,  $P_i, Q_i \in K[x_1, \cdots,x_t]$ for $i=1,\ldots,r$, and $n_1, \ldots, n_t$ are non-negative integers and $\alpha_{i,j}, \beta_{i,j} \in K$ for all $(i,j) \in \{1, \cdots, r\} \times \{1, \cdots, t\}$. $\alpha_{i,j}$ and $ c_{i,j}$ are called "multi-roots" of multi-recurrences  $(F(\nbf)), ~(G(\nbf))$ respectively.
 
 Let $S$ be a set  containing all the archimedean places of $K$ and the nonarchimedean places corresponding to the prime ideals lying above  $\alpha_{ij}$ and $c_{i,j} $. Let $s$ be the number of elements of $S$. Put $\Gamma = \{1, \ldots,r\}$ and $\G$ be the subset of $\mathbb{Z}^t$ consisting of those  $\textbf{z}$ for which at least one of the following equalities holds for some distinct  $l,k \in \Gamma$ :  
 $$\bsa_l^{\textbf{z}} = \bsa_k^{\textbf{z}}, \quad
 \textbf{c}_l^{\textbf{z}} = \textbf{c}_k^{\textbf{z}}, \quad
 \bsa_l^{\textbf{z}} = \textbf{c}_k^{\textbf{z}}.$$ Define $$A' = \sum_{i \in \Gamma}^{} \binom{t+\ri}{t} \mbox{ and } A= \max \{t, A'\}.$$
Note that $A = \max\{t,r\}$ if all the polynomials $P_1, \ldots, P_r$ are constants, whereas $A=A'$ otherwise. Finally, set $$h := h^{hom}(\textbf{V}, 1)$$ where $\textbf{V}$ is the vector consisting of the coefficients of $P_1, \cdots, P_r$ and $Q_1, \cdots, Q_r$.   For  convenience, we restate Theorem \ref{thmL}.

\begin{thm}\label{thmA'}(Theorem \ref{thmL}). 
Let $\epsilon >0$ and $$S_0 = \{v \in M_K~ :~ \max \{|\alpha_{i,j}|_v, |c_{i,j}|_v \} <1\}.$$ Assume  $\G = \{0\}$.   
 We also suppose that $F$ and $G$ are  coprime and simple (i.e. $P_i$, $Q_i$ are constants) and the inequality
 \begin{equation}\label{eq:d'}
 {\displaystyle \sum_{M_K \setminus S_0}} -\log^{-} \max \{|F(\nbf)|_v, |G(\nbf)|_v\} > \epsilon |\nbf | 
 \end{equation}
 has infinitely many solutions. There exist finitely many  proper subspaces $H_1, \cdots, H_l$ of $\mathbb{R}^t$; constants $C_1, \cdots C_l >0$ and finitely many tuples  $(\theta, \textbf{u}) \in \mathrm{GL}_t(\ZZ) \times \mathbb{Z}^t$  such that all but finitely many of solutions  satisfy  $$\nbf = \theta (\textbf{x}) + \textbf{u} \quad  \mbox{where} \quad \textbf{x} \in \Sc_i$$ with
\begin{equation*}
		\Sc_i = \set{\textbf{x} \in \ZZ^t : \abs{\textbf{x}^{\perp}} \leq C_i \log^+ \abs{\textbf{x}^{H_i}}},
\end{equation*}
for some $i \in \{1, \cdots, l \}$,  where $\textbf{x}^{\perp} $ lies in the orthogonal of $H$ and $\textbf{x}^{H_i} \in H_i$ with $\textbf{x} = \textbf{x}^{H_i} + \textbf{x}^{\perp}$.
\end{thm}

\begin{proof}
By definition of $S$, we have $S_0 \subseteq S$. In the view of Proposition \ref{propA'} and the inequality
\begin{equation*}
\begin{split}
{\displaystyle \sum_{M_K \setminus S_0}} -\log^{-} \max \{|F(\nbf)|_v, |G(\nbf)|_v\}  \leq  {\displaystyle \sum_{M_K \setminus S}} -\log^{-} \max \{|F(\nbf)|_v, |G(\nbf)|_v\} + \\
{\displaystyle \sum_{S \setminus S_0}} -\log^{-} \max \{|F(\nbf)|_v, |G(\nbf)|_v\}, 
\end{split}
\end{equation*}
we deduce that, with at most finitely many exceptions, every solution of inequality \eqref{eq:d'} also satisfies $${\displaystyle \sum_{M_K \setminus S}} -\log^{-} \max \{|F(\nbf)|_v, |G(\nbf)|_v\} > \frac{\epsilon}{2} |\nbf|.$$ Let $\{u_1, \cdots, u_N\}$ be  a basis of the group  $\mathcal{A}$ generated by $\alpha_{ij}$ and $c_{i,j}$. Let $$ f, g \in  \mathbb{K}_1 := K[ T_{1,1}^{\pm 1}, \cdots, T_{N,t}^{\pm1} ]$$
be the Laurent polynomials corresponding to $F$ and $G$ respectively.
We write $$f(T_{1,1}, \cdots, T_{N,t}) = T_{1,1}^{i_{1,1}}\cdots T_{N,t}^{i_{N,t}}  f_0(T_{1,1}, \cdots, T_{N,t})$$ and $$g(T_{1,1}, \cdots,  T_{N,t}) = T_{1,1}^{j_{1,1}}\cdots T_{N,t}^{j_{N,t}}  g_0(T_{1,1}, \cdots,  T_{N,t}),$$ where $i_{1,1}, \cdots, i_{N,t}, j_{1,1}, \cdots, j_{N,t} \in \mathbb{Z}$ and  $f_0, g_0 \in K[ T_{1,1}, \cdots, T_{N,t}]$ are such that no variable $T_{i,j}$ divides $f_0g_0,$ for any $(i, j) \in \{1, \cdots, N\} \times \{1, \cdots, t\}$. Let $F_0, G_0$ denote the simple multi-recurrence sequences associated with $f_0, g_0$ respectively. Since $f, g$ are coprime and each $T_{i,j}$ is unit in $\mathbb{K}_1$, the polynomials $f_0$ and $g_0$ are also coprime. Clearly, $u_1, \cdots, u_N \in \mathcal{O}_S^*$. Consequently, $${\displaystyle \sum_{M_K \setminus S}} -\log^{-} \max \{|F(\nbf)|_v, |G(\nbf)|_v\} = {\displaystyle \sum_{M_K \setminus S}} -\log^{-} \max \{|F_0(\nbf)|_v, |G_0(\nbf)|_v\}.$$ It therefore suffices to prove the theorem under the assumption that inequality  \eqref{eq:d'} is replaced by  $${\displaystyle \sum_{M_K \setminus S}} -\log^{-} \max \{|F_0(\nbf)|_v, |G_0(\nbf)|_v\} > \frac{\epsilon}{2} |\nbf|.$$ Note that any $\nbf$ satisfying the inequality 
\begin{equation}\label{eq:T1}
{\displaystyle \sum_{M_K \setminus S}} -\log^{-} \max \{|F_0(\nbf)|_v, |G_0(\nbf)|_v\} > \frac{\epsilon}{2} |\nbf|
\end{equation}
also satisfies $${\displaystyle \sum_{M_K \setminus S}} -\log^{-} \max \{|F_0(\nbf)|_v, |G_0(\nbf)|_v\} > \dfrac{\epsilon}{2t \sum_{i=1}^N h(u_i)} (h(u_1^{n_1})+ \cdots + h(u_N^{n_t})).$$  We set $$P_{\nbf} := (u_1^{n_1},\cdots, u_N^{n_t}) \in \mathbb{G}_m^{Nt}(K)_{S, 0}.$$ By Lemma \ref{lemF} and Lemma \ref{lemD}, the total degrees of $f_0$ and $g_0$ are bounded above by $$D := Nt c_6 \left( \max_{i,j} \{h(\alpha_{i,j}), h(c_{i,j})\} \right)^{N+1}.$$ Applying Proposition \ref{pro3}  with $\epsilon$ replaced by $ \epsilon/ t \sum_{i=1}^N h(u_i)$), we obtain finitely many hypersurfaces $H_1, \cdots, H_{t_1}$  such that, for all but finitely many solutions $\nbf$ of \eqref{eq:T1}, we have $$ (u_1^{n_1},\cdots, u_N^{n_t}) \in H_1 \bigcup \cdots \bigcup H_{t_1}.$$ Consequently, there are at most $t_1$ polynomial-exponential Diophantine equations, at least one of which is satisfied by every such $\nbf$. Fix
\begin{equation} \label{eq:1B}
\sum_{\textbf{i} \in I}^{} d_{\textbf{i}}u_1^{i_1n_1}\cdots u_N^{i_{Nt}n_t} = 0,
\end{equation}
where  $I$ is a subset of  $\{\textbf{i} \in \mathbb{N}^{Nt}~:~ i_1 + \cdots + i_{Nt} \leq C\}$ for some constant $C>0$ and $d_{\textbf{i}} \in K$.  Let $\mathcal{P}$ be a partition of the set $I$, whose constituent subsets  will be regarded as elements of $\mathcal{P}$. Associated with  $\mathcal{P}$, we consider the system of equations
\begin{equation} \label{eq:12'}
\sum_{\textbf{i} \in \lambda}^{}  d_{\textbf{i}}u_1^{i_1n_1}\cdots u_N^{i_Nn_t} =0~(\lambda \in \mathcal{P}).
\end{equation} 
 Let $M(\mathcal{P})$ be the set consisting of $\mathcal{P}$-non-degenerate solutions of the system \eqref{eq:12'}. Clearly, every solution of \eqref{eq:1B} belongs to $M(\mathcal{P})$ for some partition $\mathcal{P}$.  Let $\abs{\cdot}$ be the euclidean norm on $\RR^m$ and define $\log^+ z = \max\setb{\log z, 1}$. Let $G(\Pc)$ be the subgroup of $\ZZ^t$ consisting of $\textbf{x}$ such that
	\begin{equation*}
		\boldsymbol{\gamma}_{\textbf{i}}^{\textbf{x}} = \boldsymbol{\gamma}_{\textbf{j}}^{\textbf{x}}
	\end{equation*}
for all $ \textbf{i}, \textbf{j} \in \lambda$, where $$\boldsymbol{\gamma}_{\textbf{i}}^{\textbf{x}} = (u_1^{i_1}\cdots u_N^{i_N})^{x_1} \cdots (u_1^{i_{(t-1)N + 1}}\cdots u_N^{i_{tN}})^{x_t}.$$ 
Let $H$ be the subspace of $\RR^m$ spanned by $G(\Pc)$. Note that $G(\Pc)$ is the set of integer solutions of a system of homogeneous linear equations in $t$ variables with integer coefficients. Since the system  of equation \eqref{eq:12'} has infinitely many solutions, Lemma \ref{thm2} implies that $G(\Pc) \neq \{0\}$. Moreover, by  a generalization of Siegel's lemma due to Bombieri and Vaaler \cite{BV}, we conclude that $H$ is a proper subspace. 
We define 
\begin{equation*}
		\Sc = \set{\textbf{x} \in \ZZ^t : \abs{\textbf{x}^{\perp}} \leq c_0 \log^+ \abs{\textbf{x}^{H}}},
\end{equation*}
	with $c_0= \max\{t, |I|\}^{12}2^{4 \max\{t, |I|\}+12}d^6$ where $\textbf{x}^{\perp} $ lies in the orthogonal of $H$ and $\textbf{x}^{H} \in H$ with $\textbf{x} = \textbf{x}^{H} + \textbf{x}^{\perp}$. The desired result now follows from Lemma\ref{thm3}. 

\end{proof}

For convenience we restate Theorem \ref{thmM}.

\begin{thm}\label{thmD"}(Theorem \ref{thmM}).
Let $\epsilon >0$ and $$S_0 = \{v \in M_K~ :~ \max_{ijkl} \{|\alpha_{i,j}|_v, |c_{k,l}|_v \} <1\}.$$   
\begin{compactitem}
 \item[a)] If $F$ and $G$ are  simple (i.e. $P_i$, $Q_i$ are constants) and the $\alpha_{i,j}$ and $c_{i,j}$  are multiplicatively independent.  Then, the number of solutions of the inequality 
 \begin{equation}\label{eq:A"}
 {\displaystyle \sum_{M_K \setminus S_0}} -\log^{-} \max \{|F(\nbf)|_v, |G(\nbf)|_v\} > \epsilon |\nbf |
 \end{equation}
does not exceed $$ T_1:= t_1 \left( \dfrac{\epsilon}{2t \left(\sum_{i,j} h(\alpha_{i,j}) + h(c_{i,j}) \right)}\right) \frac{|I|!}{2} 2^{35|I|^3}d^{6|I|^2} + s \rho\left( \frac{\epsilon}{2s}\right)$$
with $$|I| = \binom{2rt^2+96r^3t^8\epsilon^{-1}}{2rt^2},$$ where, the functions $\rho(\cdot)$ is defined in Theorem \ref{thmK}, with parameter ($D=0$) and the function $t_1(\cdot)$ is defined in Proposition \ref{pro3}, with the parameters ($h^{hom}(f,g,1)=h,~n = 2rt^2,~ D=t, ~ m= 96t^8r^3$).
  \item[b)] If $t=1$ (i.e. $F, G$ are linear recurrence sequences) and $F,~G$ are  coprime and $\G = \{0\}$, then the number of solutions of inequality \begin{equation}\label{eq:l_4}
  {\displaystyle \sum_{M_K \setminus S_0}} -\log^{-} \max \{|F(n)|_v, |G(n)|_v\} > \epsilon n
  \end{equation}
does not exceed $$ T_2(\epsilon) :=  L + t_1 \left( \dfrac{\epsilon}{3r^2 \max_{i,j}\{h(c_i), h(\alpha_j), 1\}}\right) \exp((7C^m)^{8C^m}) + s \rho\left( \frac{\epsilon}{2s}\right)$$
with $$D= \max \{\deg P_i, \deg Q_i\} + 2c_7 \left(r\max_{i,j} \{h(\alpha_{i}), h(c_j)\} \right)^{2r+1};~ L=  \left(\dfrac{288 D^2(2r+1)^4d^4}{ \epsilon^2} \right)^{3/2}$$ and $$ m= 12 D^2(2r+1)^3\epsilon^{-1}3r^2 \max_{i,j}\{h(c_i), h(\alpha_j), 1\},~C= \binom{2r+1+m}{2r+1}$$ and $$c_7 =  \dfrac{58(2r+1)!(3d)^{2r+2} \log d}{(2r+1)^{2r+1}}.$$ Here, the functions $\rho(\cdot)$ is defined in Theorem \ref{thmK}, with parameters ($D= \max \{\deg P_i, \deg Q_i\} + 1$ and $t=1$) and the function $t_1(\cdot)$ is  defined in Proposition \ref{pro3}, with the parameters ($h^{hom}(f,g,1)=h,~n = 2r+1$ and , $ D, ~ m$ as defined above).
\end{compactitem}
\end{thm}

\begin{proof}
We first prove assertion $a)$. 
 By definition of $S$, we have $S_0 \subseteq S$. In particular
\begin{equation*}
\begin{split}
{\displaystyle \sum_{M_K \setminus S_0}} -\log^{-} \max \{|F(\nbf)|_v, |G(\nbf)|_v\}  \leq  {\displaystyle \sum_{M_K \setminus S}} -\log^{-} \max \{|F(\nbf)|_v, |G(\nbf)|_v\} + \\
{\displaystyle \sum_{S \setminus S_0}} -\log^{-} \max \{|F(\nbf)|_v, |G(\nbf)|_v\}, 
\end{split}
\end{equation*}
It follows from  Theorem \ref{thmK} that the number of solutions of  $${\displaystyle \sum_{S \setminus S_0}} -\log^{-} \max \{|F(\nbf)|_v, |G(\nbf)|_v\} > \frac{\epsilon}{2} |\nbf|$$ is at most $$s \rho\left( \frac{\epsilon}{2s}\right),$$ where $\rho(\cdot)$ is the function defined in  Theorem \ref{thmK}, with $D=0$. It remains to estimate the number of solutions of $${\displaystyle \sum_{M_K \setminus S}} -\log^{-} \max \{|F(\nbf)|_v, |G(\nbf)|_v\} > \frac{\epsilon}{2} |\nbf|.$$ Since $\alpha_{i,j}'$s and $c_{i,j}'$s  are multiplicatively independent, the group  $\mathcal{A}$ generated by them has rank   $N:= 2rt$. Let $$ f, g \in  \mathbb{K}_1 := K[ T_{1,1}, \cdots, T_{N,t}]$$
be the Laurent polynomials corresponding to $F$ and $G$ respectively. Observe that 
$$f(T_{1,1}, \cdots, T_{N,t})= f(T_{1,1}, \cdots, T_{r,t})$$ and $$g(T_{1,1}, \cdots, T_{N,t})= g(T_{r+1,1}, \cdots, T_{N,t})$$
and any solution $\nbf$ of
\begin{equation}\label{eq:1C}
{\displaystyle \sum_{M_K \setminus S}} -\log^{-} \max \{|F(\nbf)|_v, |G(\nbf)|_v\} > \frac{\epsilon}{2} |\nbf|
\end{equation}
also satisfies   $${\displaystyle \sum_{M_K \setminus S}} -\log^{-} \max \{|F(\nbf)|_v, |G(\nbf)|_v\} > \dfrac{\epsilon}{2t (\sum_{i,j} h(\alpha_{i,j}) + h(c_{i,j}))} (h(\alpha_{1,1}^{n_1})+ \cdots + h(c_{r,t}^{n_t})).$$  We set $$P_{\nbf} := (\alpha_{1,1}^{n_1},\cdots,\alpha_{r,t}^{n_t}, c_{1,1}^{n_1}, \cdots, c_{r,t}^{n_t}) \in \mathbb{G}_m^{2rt^2}(K)_{S, 0}.$$ It is easy to see that the total degrees of the polynomials $f$ and $g$ are bounded above by $D :=t.$ Applying  Proposition \ref{pro3}, with $\epsilon$ replaced by $ \epsilon/2t \left(\sum_{i,j} h(\alpha_{i,j}) + h(\alpha_{i,j}) \right)$), we obtain  at most $$t_1 \left( \dfrac{\epsilon}{2t \left(\sum_{i,j} h(\alpha_{i,j}) + h(c_{i,j}) \right)}\right)$$ hypersurfaces $H_1, \cdots, H_{t_1}$, where $t_1(\cdot)$  is defined in Proposition \ref{pro3},   such that,  every solution $\nbf$ of \eqref{eq:1C} satisfies $$ (\alpha_{1,1}^{n_1},\cdots,\alpha_{r,t}^{n_t}, c_{1,1}^{n_1}, \cdots, c_{r,t}^{n_t})  \in H_1 \bigcup \cdots \bigcup H_{t_1}.$$ Consequently, there are at most   $$t_1 \left( \dfrac{\epsilon}{2t \left(\sum_{i,j} h(\alpha_{i,j}) + h(c_{i,j}) \right)}\right)$$ polynomial-exponential Diophantine equations, at least one of which is satisfied by every such $\nbf$. Fix
\begin{equation} \label{eq:13'}
\sum_{\textbf{i} \in I}^{} d_{\textbf{i}} \alpha_{1,1}^{i_1n_1}  \cdots c_{r,t}^{i_{N}n_t} = 0,
\end{equation}
where  $I$ is a subset of  $\{\textbf{i} \in \mathbb{N}^{t}~:~ i_1 + \cdots + i_{2rt^2} \leq 96r^3t^8\epsilon^{-1}\}$ and $d_{\textbf{i}} \in K$. Clearly $|I| = \binom{N+96r^3t^8\epsilon^{-1}}{N}$.  Let $\mathcal{P}$ be a partition of the set $I$. The set $\lambda \subseteq I$ occurring in the partition $\mathcal{P}$ will be considered elements of $\mathcal{P}$. Given $\mathcal{P}$, consider the system of equations
\begin{equation} \label{eq:14'}
\sum_{\textbf{i} \in \lambda}^{}  d_{\textbf{i}} \alpha_{1,1}^{i_1n_1}  \cdots c_{r,t}^{i_{N}n_t} =0~(\lambda \in \mathcal{P}).
\end{equation} 
 Let $M(\mathcal{P})$ be the set consisting of $\mathcal{P}$-non-degenerate solutions of the system \eqref{eq:14'}. It is clear that every solution of \eqref{eq:13'} lies in $M(\mathcal{P})$ for some partition $\mathcal{P}$. By \cite[Lemma 2.5]{Le}, the set of partitions of $I$ which only contain subsets with at least two elements has a cardinality, which does not exceed $|I|! / 2.$ Hence, by using the fact that $\G = \{\textbf{0}\}$ and Lemma \ref{thm2}, it follows that the number of solutions of equation \eqref{eq:13'} does not exceed $$\frac{|I|!}{2} 2^{35|I|^3}d^{6|I|^2}.$$ By combining this bound with the upper bound above on the number of polynomial-exponential Diophantine equations of the shape \eqref{eq:13'}, we deduce that the inequality \eqref{eq:1C}  has at most $$ t_1 \left( \dfrac{\epsilon}{t \left(\sum_{i,j} h(\alpha_{i,j}) + h(c_{i,j}) \right)}\right) \frac{|I|!}{2} 2^{35|I|^3}d^{6|I|^2}$$ solutions.  Finally, the number of solutions of the inequality \eqref{eq:A"}  does not exceed $$t_1 \left( \dfrac{\epsilon}{2t \left(\sum_{i,j} h(\alpha_{i,j}) + h(c_{i,j}) \right)}\right) \frac{|I|!}{e} 2^{35|I|^3}d^{6|I|^2} + s \rho\left( \frac{\epsilon}{2s}\right),$$ which give us the desired result.\\
 
 We now prove assertion $b)$.  Since  $S_0 \subseteq S$, we have 
\begin{equation*}
\begin{split}
{\displaystyle \sum_{M_K \setminus S_0}} -\log^{-} \max \{|F(n)|_v, |G(n)|_v\}  \leq  {\displaystyle \sum_{M_K \setminus S}} -\log^{-} \max \{|F(n)|_v, |G(n)|_v\} + \\
{\displaystyle \sum_{S \setminus S_0}} -\log^{-} \max \{|F(n)|_v, |G(n)|_v\}.
\end{split}
\end{equation*}
 By Proposition \ref{propA'}, the number of solutions of  $${\displaystyle \sum_{S \setminus S_0}} -\log^{-} \max \{|F(n)|_v, |G(n)|_v\} > \frac{\epsilon}{2} n$$ does not exceed $$s \rho\left( \frac{\epsilon}{2s}\right),$$ where $\rho(\cdot)$ is the function defined in  Theorem \ref{thmK} with $D= \max \{\deg P_i, \deg Q_i\} + 1$ and $t=1$.  Let $u_1, \cdots, u_{N-1}$ be generators for the group $\mathcal{A}$ generated by the $\alpha_{i,1}$ and $c_{i,1}$. Clearly $N \leq 2r+1$. Let $ f, g \in  \mathbb{K}_1 := K[X, T_{1}^{\pm 1}, \cdots, T_{N-1}^{\pm 1}]$
denote the Laurent polynomials corresponding to $F$ and $G$ respectively. We write $$f(X, T_{1}, \cdots, T_{N-1}) = T_{1}^{i_{1}}\cdots T_{N-1}^{i_{N-1}}  f_0(X, T_{1}, \cdots, T_{N-1})$$ and $$g(X, T_{1}, \cdots,  T_{N-1}) = T_{1}^{j_{1}}\cdots T_{N-1}^{j_{N-1}}  g_0(X, T_{1}, \cdots,  T_{N-1}),$$ where $i_{1}, \cdots, i_{N-1}, j_{1}, \cdots, j_{N-1} \in \mathbb{Z}$ and  $f_0, g_0 \in K[ X, T_{1}, \cdots, T_{N-1}]$ with $T_{i}$ not dividing $f_0g_0$ for $i \in \{1, \cdots, N-1\} $. Let $F_0, G_0$ be the corresponding linear recurrence sequences of $f_0, g_0$ respectively. Since $f, g$ are coprime and each $T_{i}$ is a unit in $\mathbb{K}_1$, the polynomials $f_0$ and $g_0$ are also coprime. Since $u_1, \cdots, u_{N-1} \in \mathcal{O}_S^*$, we have $${\displaystyle \sum_{M_K \setminus S}} -\log^{-} \max \{|F(n)|_v, |G(n)|_v\} = {\displaystyle \sum_{M_K \setminus S}} -\log^{-} \max \{|F_0(n)|_v, |G_0(n)|_v\}.$$ It therefore remains to estimate the number of solutions of
\begin{equation}\label{eq:e'}
{\displaystyle \sum_{M_K \setminus S}} -\log^{-} \max \{|F_0(n)|_v, |G_0(n)|_v\} > \frac{\epsilon}{2} n
\end{equation}
By Lemma \ref{lemD} and Lemma \ref{lemF}, we deduce that the total degrees of polynomials $f_0$ and $g_0$ are bounded above by  $$D := 2\max \{\deg P_i, \deg Q_i\} + 4c_7 \left(r\max_{i,j} \{h(\alpha_{i}), h(c_j)\} \right)^{2r+1},$$ where $c_7 =c_6(2r+1)$ with $c_6(\cdot)$ defined in \ref{eqA}.  We now consider the sufficiently large solutions $n$ of  \eqref{eq:e'}, namely those satisfying  $$n \geq  \left(\dfrac{288 D^2(2r+1)^4d^4}{ \epsilon^2} \right)^{3/2}.$$ Put $$\eta = \dfrac{\epsilon^2}{144 D^2(2r+1)^4}.$$ By Dobrolowski's  lower bound for the absolute Weil height, together with Lemma \ref{lemG}, we obtain $$P_{n} := (n, u_1^{n},\cdots, u_N^{n}) \in \mathbb{G}_m^{N+1}(K)_{S, \eta}.$$
 Moreover, every sufficiently large solution $n$ of \eqref{eq:e'} satisfies $${\displaystyle \sum_{M_K \setminus S}} -\log^{-} \max \{|F_0(n)|_v, |G_0(n)|_v\} > \dfrac{\epsilon}{1+2\sum_{i=1}^{N-1} h(u_i)} (h(u_1^{n})+ \cdots + h(u_{N-1}^{n})).$$ 
Applying Proposition \ref{pro3}, with $\epsilon$ replaced by $ \epsilon /\left( 1+ 2\sum_{i=1}^{N-1} h(u_i) \right)$, it follows that there exist at most $t_1$ hypersurfaces $H_1, \cdots, H_{t_1}$, each of total degree at most  
$$m= 12 D^2(2r+1)^3\epsilon^{-1}3r^2 \max_{i,j}\{h(c_i), h(\alpha_j), 1\}$$
such that, every sufficiently large solution $n$ of \eqref{eq:e'} satisfies $$ (n, u_1^{n},\cdots, u_{N-1}^{n}) \in H_1 \bigcup \cdots \bigcup H_{t_1}.$$ Consequently,  there are at most $t_1(\epsilon/ \left( 3r^2 \max_{i,j}\{h(c_i), h(\alpha_j), 1\} )\right)$  polynomial-exponential Diophantine equations, at least one of which is satisfied by every sufficiently large  such solution  $\nbf$.  Such polynomial-exponential Diophantine equations have the form
\begin{equation} \label{eq:11'}
\sum_{\textbf{i} \in I}^{} P_{\textbf{i}}(n)u_1^{i_1n}\cdots u_{N-1}^{i_{N-1}n} = 0,
\end{equation}
where  $I$ is a subset of  $\{\textbf{i} \in \mathbb{N}^{N-1}~:~ i_1 + \cdots + i_{N-1} \leq m \}$ and $d_{\textbf{i}} \in K$. By Lemma \ref{lemJ},  the number of solution in $n$ of inequality \eqref{eq:e'} does not exceed $$  L + t_1 \left( \dfrac{\epsilon}{3r^2 \max_{i,j}\{h(c_i), h(\alpha_j), 1\}}\right) \exp((7C^m)^{8C^m})$$ with $$L=  \left(\dfrac{288 D^2(2r+1)^4d^4}{ \epsilon^2} \right)^{3/2}; \quad C= \binom{2r+1+m}{2r+1}$$ where the functions $\rho(\cdot)$ and $t_1(\cdot)$ are defined in Theorem \ref{thmK} and Proposition \ref{pro3} respectively. This completes the proof.
\end{proof}

In this final part of the paper, we establish an effective version of the Hadamard quotient theorem (Theorem \ref{thmN}). For the reader's convenience, we first recall the statement of the result.
 Let $K$ be a number field of degree $d$ and  $(F_n)_n$ and $(G_n)_n$  be two linear recurrence sequences with Binet representation 
$$
G_n = \sum_{i=1}^{r} Q_i(n)c_{i}^{n}, \quad \mbox{and}\quad F_n = \sum_{i=1}^{r} P_i(n)\alpha_{i}^{n},
$$
where $P_i, Q_i \in K[x]$ have coefficients  $a_{ij}, b_{kl}$ respectively,  and  $\alpha_{i}, c_i\in K$ for all $i \in \{1, \cdots, r\} $. Let $S$ be the set containing all nonarchimedean places corresponding to the prime ideals above $\alpha_{i}$ and $c_{i} $ and all the archimedean places of $K$. We denote the cardinality of $S$ by $s$. Finally,  set $h_1= {\rm max}_{i,j, k, l}\{ h( a_{ij}),  h( b_{kl}),1\}$.

\begin{pro}(Theorem \ref{thmN}.)
Assume that the roots  of the characteristic polynomials of $(F_n)_n$ and  $(G_n)_n$ generate together a torsion-free multiplicative group denoted $\mathcal{A}$. We also suppose that  the set $\mathcal{R}$ of $n$ such that $G(n) \neq 0$ and $F(n)/G(n) \in \mathcal{O}_S$ contains at least $$C = T_2 \left( \frac{1}{4s(s+1)d^4}\right) + \rho\left( \frac{1}{8sd^4}\right) + (s-1) \rho\left( \frac{1}{4s^2d^4}\right) +1$$ elements, where  $T_2(\cdot)$  is the function defined in Theorem \ref{thmM} b), with $r$ replaced by $l$ defined in the Lemma \ref{lemA"} and $\rho(\cdot)$ is the function defined in Lemma \ref{lemA"}. Then there exists a nonzero polynomial $P(X) \in K[X]$ such that both sequences $(P(n)F_n/G_n)_n$ and $(G_n/Q(n))_n$ are linear recurrence sequences.
\end{pro}

\begin{proof}
We distinguish two cases according to the number of roots of the characteristic polynomial of $(G_n)_n$. If this polynomial has only one root, then it is readily seen that both   $Q_1(n)F_n/G_n$ and $G_n/Q_1(n)$ are linear recurrence sequences. We may therefore assume henceforth that the characteristic polynomial of the sequence $(G_n)_n$ has at least two roots.  Let $u_1, \cdots, u_{N-1}$ be a basis of the group $\mathcal{A}$. Clearly $N \leq 2r+1$. Let $ f, g \in  \mathbb{K}_1 := K[X, T_{1}, \cdots, T_{N-1}, T_{1}^{-1}, \cdots, T_{N-1}^{-1}]$
be the Laurent polynomials corresponding to the sequences $F_n$ and $G_n$, respectively. We write $$f(X, T_{1}, \cdots, T_{N-1}) = T_{1}^{i_{1}}\cdots T_{N-1}^{i_{N-1}}  f_0(X, T_{1}, \cdots, T_{N-1})$$ and $$g(X, T_{1}, \cdots,  T_{N-1}) = T_{1}^{j_{1}}\cdots T_{N-1}^{j_{N-1}}  g_0(X, T_{1}, \cdots,  T_{N-1}),$$ where $i_{1}, \cdots, i_{N-1}, j_{1}, \cdots, j_{N-1} \in \mathbb{Z}$ and  $f_0, g_0 \in K[ X, T_{1}, \cdots, T_{N-1}]$ are such that no  $T_{i}$ divides $f_0g_0$ for $i \in \{1, \cdots, N-1\}. $ Let $F_0, G_0$ be the corresponding linear recurrence sequences of $f_0, g_0$, respectively. Let $G$ be the greatest common divisor of  $f_0$, and $ g_0$, normalized to be monic, and write $$ f_0 = Gf_1, \quad g_0 = G g_1,$$  $f_1, g_1 \in K[ X, T_{1}, \cdots, T_{N-1}]$ are coprime.  By Lemma \ref{lemD} and Lemma \ref{lemF}, the total degrees of  $f_0$ and $g_0$ are bounded above by  $$D := 2\max \{\deg P_i, \deg Q_i\} + 4c_7 \left(r\max_{i,j} \{h(\alpha_{i}), h(c_j)\} \right)^{2r+1}.$$  It follows from  Gelfond's lemma (see \cite[Lemma 4.1.3]{evertse-gyory-2022}) that  $$ h^{hom}(f_1,g_1,1) \leq (2r+1 +D)^{2r+3}h_1.$$ Let $F_1(n)$ and $G_1(n)$ be two linear recurrences corresponding to $f_1, g_1$, respectively.

In other word $$F_1(n) = f_1(n, u_1^n, \cdots, u_{N-1}^n) \quad \mbox{and} \quad  G_1(n) = g_1(n, u_1^n, \cdots, u_{N-1}^n).$$ It is readily verified that    $$\mathcal{R}= \left \{ n \in \mathbb{N}~~:~~ \frac{F_1(n)}{G_1(n)} \in \mathcal{O}_S \right \}$$ and $${\displaystyle \sum_{M_K \setminus S}} -\log |G_1(n)|_v  \leq {\displaystyle \sum_{M_K \setminus S}} -\log^{-} \max \{|F_1(n)|_v, |G_1(n)|_v\}.$$ Since the group  $\mathcal{A}$  is torsion-free, no quotient of two distinct elements of  $\mathcal{A}$  is a root of unity. It follows from Kronecker's theorem (see \cite[Theorem 1.5.9]{bombieri-gubler-2006}) that there exists a place $w$ of $K$ at which the roots of the characteristic polynomial of the recurrence sequence $G_1(n)$ do not all have the same absolute value.
Observe that the roots of the characteristic polynomial of  $G_1(n)$ are of the form $$u_1^{m_1}\cdots u_{N-1}^{m_{N-1}}, \quad \mbox{with}\quad \sum_{i=1}^{N-1}m_i \leq D.$$

Let $G_2(n)$ be the linear recurrence sequence obtained from $G_1(n)$ by dividing its Binet representation  by a suitable root, chosen so that all the characteristic roots of  $G_2(n)$ have $w$-adic absolute value at least 1. We denote by $\beta_1, \cdots \beta_r$ the roots of the characteristic polynomial of $G_2(n)$. Hence, we obtain $${\displaystyle \sum_{\nu \in S}} \log |G_2(n)|_v  = {\displaystyle \sum_{\nu \in  M_K \setminus S}} -\log |G_2(n)|_v  \leq {\displaystyle \sum_{\nu \in M_K \setminus S}} -\log^{-} \max \{|F_1(n)|_{\nu}, |G_1(n)|_{\nu}\}.$$ Without loss of generality we put $$|\beta_1|_v =\max_{\nu \in S} |\beta_1|_{\nu} >1  \quad \mbox{and} \quad  \epsilon = \frac{1}{2} \log \left|\beta_1\right |_v \geq \frac{1}{4sd^4}.$$

\begin{lem}\label{lemA"}
 The equation 
\begin{equation}\label{eq:T}
\log |G_1(n)|_v < \epsilon n,
\end{equation}
has  at most $$ \rho(1/ (8sd^4)):= 4 \left( 2^{60l^2}\delta(1/ (8sd^4) )^{-7l}\right)^{s+1}\exp((7l^{\sigma})^{8l^{\sigma}}) +\sigma^l + L(1/ (8sd^4)) $$ with $$\eta(x) =  \dfrac{x}{21l^3 h},~\sigma = \max\{\deg P_i, \deg Q_i \},~ h = (2r+1 +D)^{2r+3}h_1,~ \delta(x) =  \dfrac{x}{7l^2 h} $$ and $$ L(x):= \max\left\{ (2\sigma)\exp\left((2\sigma)h \right), (4 \eta(x)^{-1}c_6(2r) h^l(\sigma+1)(l-1))^{3/2}\right\},~ l= (2r+1 +D)^{2r+1}$$ and $$D := 2\max \{\deg P_i, \deg Q_i\} + 4c_7 \left(r\max_{i,j} \{h(\alpha_{i}), h(c_j)\} \right)^{2r+1}$$
solutions in $n.$
\end{lem}

\begin{proof}
The proof of this lemma is analogous to that of Theorem \ref{thmK}, and is therefore omitted. 
\end{proof}
By  Theorem \ref{thmM} b), the number of solutions of the inequality $${\displaystyle \sum_{M_K \setminus S}} -\log^{-} \max \{|F_1(n)|_v, |G_1(n)|_v\} > \frac{\epsilon}{s+1} n$$ is atmost   $$ T_2 \left( \frac{1}{4s(s+1)d^4}\right)$$ is the function defined in Theorem \ref{thmM} b), with $r$ replaced by $l$ in the Lemma \ref{lemA"}, and with  the function $\rho(\cdot)$ as defined in the Lemma \ref{lemA"}. It follows from the hypothesis that the inequality $${\displaystyle \sum_{\nu \in S}} \log |G_2(n)|_v \leq \epsilon/(s+1) n$$ has at least $$ \rho \left(\frac{1}{8sd^4}\right) + (s-1) \rho \left(\frac{1}{4s^2d^4} \right)$$ solutions. On the other hand, by Theorem \ref{thmK}, for every  $\nu \in S$, the inequality  $$\log |G_2(n)|_{\nu}  < -  \epsilon/s n $$ has at most $\rho(1/ (4s^2d^4)$ solutions. Consequently,   $$\log |G_2(n)|_{v}  <  \epsilon n $$ has at least $$\rho \left(\frac{1}{8sd^4}\right) +1$$ solutions. This contradicts Lemma \ref{lemA"}, and hence completes the proof of the proposition.
\end{proof}

\section*{Acknowledgement}

I would like to thank Umberto Zannier for pointing out his paper \cite{CZ3} and  Robert Tichy, Martin Widmer and Daniel Smertnig for many hepful remarks and suggestions on an earlier version of this article. 

The author was funded  by the Austrian Science Fund (FWF), Grant-DOI\\
10.55776/ESP157132.

\subsection*{Use of AI tools.} GPT-6 was used exclusively for language editing and for deriving the identity \eqref{eq:Z}, which follows directly from the hockey-stick identity.



\begin{thebibliography}{0}
\bibitem{A} L. M. Adleman, C. Pomerance, R. S. Rumely,  \emph{On distinguishing prime numbers from
composite numbers.}  Ann. Math. \textbf{117} (1),  173 -- 206 (1983).

\bibitem{AAGT} M. Agrawal, S. Akshay, B. Genest, and  P. S. Thiagarajan, \emph{ Approximate verification of the symbolic dynamics of Markov chains},  \textit{J. ACM} \textbf{62}(1) (2015) 2:1 -- 2:34.

\bibitem{AR} N. Ailon, Z. Rudnick,  \emph{Torsion points on curves and common divisors of $a^k-1$ and $b^k-1$.}  Acta Arith., . \textbf{113} (1),  31 -- 38 (2004).

\bibitem{ahlgren-1999} S. Ahlgren,
			\emph{The set of solutions of a polynomial-exponential equation,}
			\textit{Acta Arith.} \textbf{87} (1999), 189--207.

  \bibitem{AV} F. Amoroso and E. Viada,
			\emph{On the zeros of linear recurrence sequences.}
			\textit{Acta Arith.} \textbf{147 }(4) (2011), 387 -- 396.                                           
            
\bibitem{Asc} M. Aschenbrenner \emph{Ideal membership in polynomial ring over the integers}, J. Amer. Math. Soc.  \textbf{17(2)} (2004),  407 -- 441.
\bibitem{Ay} C. Ayoub,  \emph{On constructing bases for ideals in polynomial rings over the integers}, J. Numb. Theory  \textbf{17} (1983),  no. 2, 204 -- 225.

\bibitem{Ba}  W. Baur,   \emph{Rekursive Algebren mit Kettenbedingungen.}  Z. Math. Logik Grundlagen Math., \textbf{20}(1974),  37 -- 46.

\bibitem{BM1} J. Berstel, and M. Mignotte  \emph{Deux propri\'et\'es d\'ecidables des suites r\'ecurrentes lin\'eaires}.  Bull. Soc. Math. France  \textbf{104} (1976).

\bibitem{BS} F. Beukers and C. Smyth, F, \emph{Cyclotomic points on curves.} Proceedings of the Milennial Conference in Urbana-Champaign, 2002.

\bibitem{BV} E. Bombieri and J.Vaaler, \emph{On Siegel's lemma}, Invent. Math.,\textbf{73} (1983), 11 -- 32.

\bibitem{bombieri-gubler-2006} E. Bombieri and W. Gubler,
			\emph{Heights in Diophantine geometry,}  Cambridge university press,  2006.

\bibitem{Bo} I. Borosh , M. Flahive, D. Rubin, and B. Treybig, \emph{A sharp bound for solutions of
linear Diophantine equations}, Proc. Amer. Math. Soc., \textbf{105} (1989), 844 -- 846. 

 \bibitem{BG}Y. Bugeaud, K. Gy\H{o}ry and K. Yu,  \emph{Bounds for the solutions of unit equations},  Acta Arith.,  \textbf{74} (1996),  67 -- 80.

 \bibitem{BCZ}Y. Bugeaud, P. Corvaja and U. Zannier,  \emph{An upper bound for the G. C. D. of $a^n-1$ and $b^n-1$},  Math. Z.,  \textbf{243}(1),  79 -- 84 (2003).

 \bibitem{C} J. W. S. Cassels, \emph{An introduction to the geometry of numbers},  Springer,  1959.
 
\bibitem{CZ3} P. Corvaja and U. Zannier,  \emph{On the number of integral points on algebraic curves}, J. reine angew. Math.,  \textbf{565} (2003),  27 -- 42.

\bibitem{CZ1} P. Corvaja and U. Zannier,  \emph{On the greatest prime factor of $(ab+1)(ac+1)$},  Proc. Amer. Math. Soc.,   \textbf{131}(6),  1705 -- 1709 (2003).

\bibitem{CZ2} P. Corvaja and U. Zannier,  \emph{A lower bound for the height of a rational function at $S$-unit points},  Monatsh. Math.,  \textbf{144}(3),  203 -- 224 (2005).

\bibitem{CZ5} P. Corvaja and U. Zannier,  \emph{Finiteness of integral value for the ratio of two linear recurrence},  Invent. math.,  \textbf{149}(3),  431 -- 451 (2002).

\bibitem{CZa} P. Corvaja and U. Zannier, \emph{Some cases of Vojta’s conjecture on integral points
over function fields}, J. Algebraic Geom.  \textbf{17} (2008), p. 295--333. Addendum in
Asian Journal of Math.,  \textbf{14}(2010), p. 581 -- 584

\bibitem{CZ4} P. Corvaja and U. Zannier, \emph{On the maximal order of a torsion point on a curve in $\mathbb{G}_m^n$}, Rend. Lincei Mat. Appl.  \textbf{19} (2008): 73 -- 78, 2010.

\bibitem{Ed}  H. M. Edwards, \emph{Kronecker's views on the foundations of mathematics.} The History of Modern
Mathematics,  Vol. I (Poughkeepsie, NY, 1989), Academic Press, Boston, MA, 1989, P. 67--77.

\bibitem{Ev}  T. Evans  \emph{Some connections between residual finiteness, finite embeddability and the word
problem}, J. London Math. Soc.  \textbf{(2)1} (1969),  399 -- 403.

\bibitem{E1} J.-H. Evertse, \textit{On sums of $S$-units and linear recurrences}. Compos. Math. \textbf{53(2)} (1984), 225--244.

\bibitem{E} J.-H. Evertse, \textit{An improvement of the quantitative Subspace Theorem}, Compos. Math. \textbf{101(3)} (1996), 225--311.

\bibitem{EF} J.-H. Evertse and  R. G. Ferretti, \textit{An further improvement of the quantitative Subspace Theorem}, Ann. of Math., \textbf{177} (2023), 513--590.

\bibitem{ES} J.-H. Evertse and H. P. Schlickewei, \textit{A quantitative version of the absolute  Subspace Theorem}, J. Reine. Angew. Math., \textbf{548} (2002), 21--127.

\bibitem{evertse-gyory-2022} 
			\textsc{J.-H. Evertse and K. Gy\H{o}ry},
			Effective results and methods for Diophantine equations over finitely generated domains,
			\textit{London Math. Soc. Lecture Note Ser.} \textbf{475},
			Cambridge University Press, Cambridge, 2022.

            \bibitem{FW} G. Faltings and G. Wutholz, \textit{Diophantine approximation on projective space}, Invent. Math., \textbf{116} (1994), 109-- 138.

\bibitem{FH1} C. Fuchs and S. Heintze, \emph{On the growth of multi-recurrences}. Arch. Math. (Basel) \textbf{119} (2022), 489--494. 

\bibitem{FH} C. Fuchs and  S. Heintze, \emph{On the growth of linear recurrences in function fields}.  Bull. Austr. Math. Soc. \textbf{104}(1) (2021). 

\bibitem{FN} C. Fuchs  and  A. Noubissie,    \emph{On  Growth of Multi-recurrences.}  Journ. Number Theory \textbf{293} (2027),  105 -- 127.

\bibitem{GM} G. Gallo and B. Mishra,  T. Evans  \emph{A solution to Kronecker's problem}, Appl. Algebra in Engrg. Comm.
Comput. \textbf{5} (1994), no. 6,  343 -- 370.

\bibitem{GW} N. Grieve and J. Wang,    \emph{Greatest common divisors with moving targets and consequences for linear recurrence sequences.}  Trans. American Math. Soc.,  \textbf{373} (2020),  8095 -- 8126.

\bibitem{GY}K. Gy\H{o}ry and K. Yu,  \emph{Bounds for the solutions of $S$-unit equations and decomposable form equations},  Acta Arith.,  \textbf{126} (2006),  9 -- 41.

 \bibitem{HN}  K. Hentzelt and E. Noether, \emph{ Zur Theorie der Polynomideale und Resultanten}. Math. Ann. \textbf{88}(1923), 53--79.
            
\bibitem{H} G. Hermann,  \emph{die frage der endlich vielen schritte in der theorie der polynomiadeale},  Math. Ann.,  \textbf{95} (1926),  736-- 788.

\bibitem{HL} S. Hern\'andez, F. Luca,  \emph{On the largest prime factor of $(ab + 1)(ac + 1)(bc + 1)$}, Bol. Soc. Mat. Mex., (3)\textbf{9} (2), 235 -- 244  (2003).

\bibitem{K}  J. Kollar, \emph{Sharp effective Nullstellensatz}.J. Amer. Math. Soc. \textbf{1}(1988),  963 -- 975.

\bibitem{L1} S. Lang, \emph{Fundamentals of Diophantine Geometry},  Springer-Verlag, New York, 1983.

\bibitem{La} S. Lang, \emph{Division points on curves}, Annali di matematica pura ed applicata (IV) LXX (1965) pp. 229 -- 234.

\bibitem{laurent-1989}
			\textsc{M. Laurent},
			Equation exponentielles polynomes et suite recurrentes lineaires II,
			\textit{J. Number Theory} \textbf{31} (1989), 24-53.

\bibitem{Le} L. Leroux, \emph{Computing the torsion points of a variety defined by lacunary polynomials}. Math. Comp. \textbf{81(279)} (2012), 1587--1607.

\bibitem{L} A. Levin, \emph{Greatest common divisors and Vojta’s conjecture for blowups of algebraic tori}. Invent. math. \textbf{215}, 493 -- 533 (2019). 

 \bibitem{LM}  T. Loher and D. Masser,  \emph{Uniformly counting points of bounded height,} Acta Arith., \textbf{111}(2004) : 277 -- 297.

\bibitem{LV} B. H.  Loxton and  and A. J. Van der Poorten, \emph{On the growth of linear recurrences }.   Math. Proc. Camb. Phil. Soc. (1977), \textbf{81}, 369.

\bibitem{AN} F. Luca and A. Noubissie,  \emph{Linear combinations of factorial and $S$-units in a ternary recurrence sequence with a double root}.  Periodica Mathematica Hungarica, \textbf{86}, 422--441  (2023).

\bibitem{LOW1} F. Luca and J. Ouaknine and J. Worrell , \emph{Universal Skolem set with positive lower density}.  Proceedings of MFCS 22, LIPIcs 241, 2022.

\bibitem{LOW2} F. Luca and J. Ouaknine and J. Worrell , \emph{Universal Skolem Set }.  Proceedings of LICS 21,  2021.

\bibitem{LOW} F. Luca, J. Ouaknine, and J. Worrell,  Algebraic model checking for discrete linear dynamical Systems   \textit{FORMATS 2022}, Lecture note in Computer Sciences, vol. 13465, springer, 22, pp. 3-- 15.

\bibitem{LMNOW} F. Luca,  J. Maynard, A. Noubissie, J. Ouaknine,  J. Worrell,   \emph{Skolem meets Bateman-Horn},   To appear in TheoretiCS Journal.//arxiv.org/abs/2308.01152v2 

\bibitem{Mas} D. W. Masser, \emph{A note on Siegel's Lemma}, The Rocky Mountain Journal of Mathematics, \textbf{26} (1996) No. 3, 1057 -- 1068. 

\bibitem{Ma} G. Maurin, \emph{Courbes alg\'ebriques et \'equations multiplicatives}, Math. Ann., \textbf{341} (2008): 789 -- 824. 

\bibitem{MM} E. Mayr and A. Meyer, \emph{ The complexity of the word problems for commutative semigroups
and polynomial ideals}, Adv. Math. \textbf{46} (1982),  no. 3, 305--329.




\bibitem{MST} M. Mignotte, T. Shorey, and R. Tijdeman  \emph{The distance between terms of an algebraic recurrence sequence}.  J. fur die reine und angewandte Math., \textbf{349} (1984) . 


\bibitem{N}  A. Noubissie,    \emph{A Quantitative Growth of  Linear Recurrences.}  Journal of Aus. Math. Soc., \textbf{120(1)},  107 -- 140  (2026).

\bibitem{OS} A. Ostafe and I. Shparlinski, \emph{On the Skolem problem and some related questions for parametric families of linear recurrence sequences},Canada Journal of Math., \textbf{74(3)} : 773 -- 792, 2022.

\bibitem{OW15} J. Ouaknine, and J. Worrell,   \emph{On linear recurrence sequences and loop termination},  \textit{ACM SIGLOG}, New \textbf{2} (2015) no. 2, 4--13.

\bibitem{OW3} J. Ouaknine and J. Worrell . 2014c.  \emph{Ultimate Positivity Problem  is Decidable for Simple Linear Recurrence Sequences}.  In Proc. ICALP (LNCS), Vol. 8573. Springer. 


  
\bibitem{R} B. Renschuch, \emph{ Beitr\H{a}ge zur konstruktiven Theorie der Polynomideale.} XVII/1. Zur Hent
zelt/Noether/Hermannschen Theorie der endlich vielen Schritte, Wiss. Z. P?dagog. Hochsch.
"Karl Liebknecht" Potsdam  \textbf{24}(1980), no. 1 87 -- 99.

\bibitem{RT} D. Roy and J.Thunder, \emph{A note on Siegel's Lemma over number fields}, Monatsh. Math.,\textbf{120} (1995), 307 -- 318.

\bibitem{PS} P. Samuel, \emph{Algebraic theory of numbers}, Hermann publisher, Paris, 1967. 

\bibitem{Se}  A. Seidenberg,  \emph{ Constructions in algebra}. Trans. Amer. Math. Soc. , \textbf{197} (1974), 273--313.

\bibitem{SS} H. P. Schlickewei and W. M. Schmidt , \emph{The number of solution of the polynomials-exponential equations}. Compos. Math. \textbf{120} (2000), 193--225.

\bibitem{VS} H. P. Schlickewei and A. J. Van der Poorten, \emph{The growth conditions for recurrence sequences}. Macquarie Math. Report \textbf{82-0041} (August 1982), Macquarie University, Australia, 2109. 

\bibitem{S} W. M. Schmidt, \emph{The zero multiplicity of linear recurrence sequences}. Acta Math. \textbf{182} (1999), 243--282.

\bibitem{Si} H. Simmons,  \emph{The solution of a decision problem for several classes of rings}. Pacific J. Math.,, \textbf{34}, No. 2, 1970.

\bibitem{So}  M. Sombra,   \emph{A sparse effective Nullstellensatz},  Adv. in Appl. Math., \textbf{22} (1999), 271--295.


\bibitem{V}  N. K. Vereshchagin,   \emph{The problem of appearance of a zero in a linear recurrence sequence (in Russian)}.  Mat. Zametki   \textbf{38, 2} (1985).
 
\bibitem{Wa}  S. S. Wainer,   \emph{ A classification of the ordinal recursive functions}, Arch. Math. Logik Grundla
genforsch. , \textbf{13} (1970), 136--153.



\bibitem{W} M. Waldschmidt, \emph{Diophantine approximation on linear algebraic groups}. Vol. \textbf{326} of Grundlehren der Mathematischen Wissenschaften [Fundamental Principles of Mathematical Sciences]. Springer-Verlag, Berlin, 2000. Transcendence properties of the exponential function in several variables.
\bibitem{X} Z. Xiao, \emph{Greatest common divisors for polynomials in almost units and applications to linear recurrence sequences,} Math. Z. \textbf{306(4)} (2024), Paper No. 61, 42 p.
 
\end{thebibliography}
\bibliographystyle{plain}

\end{document}